\documentclass[12pt,a4paper]{article}
\usepackage[latin1]{inputenc}
\usepackage[T1]{fontenc}
\usepackage{amsmath,amssymb}
\usepackage{geometry}
\newtheorem{theorem}{Theorem}
\newtheorem{definition}{Definition}
\newtheorem{proposition}{Proposition}
\newtheorem{corollary}{Corollary}
\newtheorem{lemma}{Lemma}
\newtheorem{remark}{Remark}
\newtheorem{example}{Example}
\newcommand\Var{\mathrm {Var}}
\usepackage{hyperref}

\newcommand\Cov{\mathrm{Cov}}
\newcommand\Card{\mathrm{Card}}
\newcommand\argmin{\mathrm{argmin}}
\newcommand\MSE{\mathrm{MSE}}
\newcommand\IMSE{\mathrm{IMSE}}
\renewcommand{\textbf}[1]{\begingroup\bfseries\mathversion{bold}#1\endgroup}
\usepackage{xcolor}
\usepackage{tikz,tkz-tab}
\usepackage{dsfont}
\usepackage{mathrsfs}
\usepackage{hyperref}
\usepackage{multirow}
\usepackage{float}
\date{}
\makeatletter
\def\hlinewd#1{%
	\noalign{\ifnum0=`}\fi\hrule \@height #1 %
	\futurelet\reserved@a\@xhline}
\makeatother
\title{The reproducing kernel Hilbert space approach in nonparametric regression problems with correlated observations}
\begin{document}
\maketitle
\begin{center}
	D. BENELMADANI, 
	 K. BENHENNI  and S. LOUHICHI\\[0.2cm]
{\it
		Laboratoire Jean Kuntzmann (CNRS 5224), Université Grenoble Alpes, France.
	}
{\footnotesize djihad.benelmadani@univ-grenoble-alpes.fr, karim.benhenni@univ-grenoble-alpes.fr, sana.louhichi@univ-grenoble-alpes.fr.}
\end{center}
\textbf{Abstract:}
In this chapter we investigate  the problem of estimating the regression function in models with correlated observations. The data is obtained from several experimental units, each of them forms a  time series. 
Using the properties of the Reproducing Kernel Hilbert spaces, we construct a new estimator based on the inverse of the autocovariance matrix of the observations. We give the asymptotic expressions of its bias and its variance. In addition, we give a theoretical comparison between this new estimator and the popular one  proposed by Gasser and Müller, we show that the proposed estimator has an asymptotically smaller variance then the classical one. Finally, we conduct a simulation study to investigate the performance of the proposed estimator and to compare it to the  Gasser and Müller's estimator in a finite sample set.\\[0.3cm]
\section{Introduction}

One of the situations that statisticians encounter in their studies is the  estimation of a whole function based on partial observations of this function. For instance, in pharmacokinetics one wishes to estimate the concentration-time of some injected medicine in the organism, based on the observations of the concentration from blood tests over a period of time.  In statistical terms, one wants to estimate a function, say $g$, relating two random variables: the explanatory variable $X$ and the response variable $Y$, without any parametric restrictions on the function $g$. The statistical model often used is the following: $Y_i=g(X_i)+\varepsilon_i$ where $(X_i,Y_i)_{1 \leq i \leq n}$ are $n$ independent replicates of $(X,Y)$ and  $\{\varepsilon_i, i=1,\cdots,n\}$ are centered random variables (called errors).

The most intensively treated model has been the one in which $(\varepsilon_i)_{1 \leq i \leq n}$ are independent errors and $(X_i)_{1 \leq i \leq n}$ are fixed within some domain. We mention the works of Priestly and Chao (1972) \cite{Priestly and Chao}, Benedetti (1977) \cite{Benedetti} and Gasser and Müller (1979) \cite{Gasser and Müller 1979} among others.  However, the independence of the observations  is not always a realistic assumption. For instance, the growth curve models are usually used in the case of longitudinal data, where the same experimental unit is being observed on multiple points of time. As a real life example, the heights observed on the same child are correlated.
The temperature observations measured along the day are also correlated. For this, we  focus, in this paper, on the nonparametric kernel estimation problem where the observations are correlated.

In the current chapter, we consider a situation where the data is generated from  $m$ experimental units each of them having $n$  measurements of the response. For this data, we consider the so-called fixed design regression model with repeated measurements given by,
\begin{equation}
Y_j(t_i)=g(t_i)+\varepsilon_j(t_i)~~\text{for}~i=1,\cdots,n~\text{and}~j=1,\cdots,m,\label{model}
\end{equation}
where $\{\varepsilon_j,j=1,\cdots, m\}$ is a sequence of i.i.d. centered error processes with the same distribution as a process $\varepsilon$. The non correlation of the errors $\{\varepsilon_j,j=1,\cdots, m\}$ is a natural assumption since it is equivalent to assuming that the experimental units (in general individuals) are independent.

This model is usually used in the growth curve analysis and dose response problems, see for instance, the work of Azzalini (1984) \cite{Azzalini}. It  has also been considered by Müller (1984) \cite{Müller} with $m=1$, where he supposed that the observations are asymptotically uncorrelated when the number of observations tends to infinity, i.e.,  $\Cov (\varepsilon(s),\varepsilon(t))=O(1/n)$ for $s \ne t$, which is not a realistic assumption, for instance, in the growth curve analysis and temperature.

The correlated observations case was considered by Hart and Wherly (1986) \cite{Hart and Wherly}, who investigated the estimation of $g$ in Model \eqref{model} where $\varepsilon$ is a stationary error process. Using the kernel estimator proposed by Gasser and Müller, see Gasser and Müller (1979) \cite{Gasser and Müller 1979},  they proved the consistency in $\mathbb L^2$ space of this estimator, when the number of experimental units $m$ tends to infinity, but not when $n$ tends to infinity as in the case of independent observations.

The assumption of stationarity made on the observations is however restrictive. In the previous pharmacokinetics example for instance, it is clear that the concentration of the medicine will be high at the beginning then decreases with time. For this, we shall investigate the estimation of $g$ in Model \eqref{model} where $\varepsilon$ is not necessarily a stationary error process. This case was partially investigated by Benhenni and Rachdi (2016) \cite{Benhenni and Rachdi} and Ferreira et al. (1997) \cite{Ferreira}, where the Gasser and Müller's estimator was used.

In this chapter, we propose a new estimator for the regression function $g$ in Model \eqref{model}. This estimator, which is also a linear kernel estimator, is based on the inverse of the autocovariance matrix of the observations, that we assume known and invertible.

The proposed estimator was inspired by the work of Sacks and Ylvisaker (1966, 19868, 1970) \cite{Sacks and Ylvisaker,Sacks and Ylvisaker II,Sacks and Ylvisaker III} but in a different context than ours. They considered the parametric model: $Y(t)=\beta f(t)+\varepsilon(t)$ where $\beta$ is an unknown real parameter and  $f$ is a known function belonging to the Reproducing Kernel Hilbert Space associated to the autocovariance function of the error process $\varepsilon$, denoted by RKHS($R$). They also assumed that the autocovariance matrix is known and invertible. It is worth noting  that the Reproducing Kernel Hilbert Spaces have been used in several domains, for instance, in Statistics by  Sacks and Ylvisaker (1966) \cite{Sacks and Ylvisaker} and more recently  by Dette et al. (2016) \cite{Dette 2016}, in Mathematical Analysis in Schwartz (1964) \cite{Schwartz} and in Signal Processing in Ramsay and Silverman (2005) \cite{Ramsay}.

We also give the asymptotic statistical performance of the proposed estimator and we compare it to the classical Gasser and Müller's estimator (GM estimator), proving, in particular, that the proposed estimator outperforms the GM estimator, in the sense that it has  an asymptotically smaller variance, wheras they both are asymptotically unbiased. This can be argued by the fact that, in statistics in general, the best linear estimator (or optimal predictor) is based on the inverse of the autocovariance matrix, see for instance, Benhenni and Cambanis  (1992) \cite{Benhenni and Cambanis}, whereas  the GM estimator does not  take into account this correlation requirement. In addition, the GM estimator is an approximation of an integral and, as known in statistics, the best linear approximation of an integral is based on some projection property.

This chapter is organized as follows. In section \ref{section2projection}, we construct our proposed estimator for the function $g$ in Model \eqref{model} where $\varepsilon$ is a centered, second order error process with a continuous autocovariance function $R$.  It is constructed  through the following function defined, for $x\in [0,1]$, by, \begin{equation}
f_{x,h}(t)=\int_{0}^{1} R(s,t)\varphi_{x,h}(t) \;ds~~\text{where}~~ \varphi_{x,h}(t)=\frac{1}{h}K\Big(\frac{x-s}{h}\Big)~~~\text{for}~~ t \in [0,1],
\label{f_xh}
\end{equation}
where $K$ is a Kernel and  $h=h(n)$ is a bandwidth.

We shall see that this function belongs to the RKHS($R$). This allows us to use the properties of this space to control the variance of the proposed estimator.  These properties were introduced by Parzen (1959) \cite{Parzen} to solve various problems in statistical inference on time series.
We also give, in this section, the analytical expressions of this estimator for the generalised Wiener process and the Ornstein-Uhlenbeck process, since the analytical expression of the inverse of the autocovariance  matrix can be derived  for this class of processes.

In Section \ref{localasymptoticresultsprojection}, we derive the asymptotic performances of this estimator. We give an asymptotic expression of the weights of this linear estimator, which is used to derive the asymptotic expression of its bias. The properties of the RKHS($R$) not only allow us to obtain the asymptotic expression of the variance, but also to find the optimal rate of convergence of the residual variance. After obtaining the asymptotic expression of the Integrated Mean Squared Error (IMSE), we derive the asymptotic optimal bandwidth with respect to the IMSE criterion.  Moreover, we prove the asymptotic normality of the proposed estimator.

In Section \ref{comparisonGM}, we give a theoretical comparison between the new estimator and  the Gasser and Müller's estimator. We prove that the proposed estimator has, asymptotically, a smaller variance than that of  Gasser and Müller. Moreover, the proposed estimator has an asymptotically smaller IMSE, for instance, in the case of a Wiener process $\varepsilon$.

In Section \ref{simulationprojection}, we conduct a simulation study in order to investigate the performance of the proposed estimator in a finite sample set, then we compare it with the Gasser and Müller's estimator for different values of the number of experimental units and different values of the sample size.  Since the classical cross-validation criterion is shown to be inefficient in the presence of correlation (see for instance, Altman (1990) \cite{Altman 1990}, Chiu (1989) \cite{Chiu 1989} and Hart (1991, 1994) \cite{Hart 1991, Hart 1994}), we use the optimal bandwidth that minimizes the exact IMSE, obtained using the Conjugated Gradient Algorithm. The results of this simulation study confirm our theoretical statements given in Section 3 and Section 4.

Finally, the supplementary materials section is dedicated to the proofs of the theoretical results, in addition to an appendix about the RKHS($R$) and some technical details.
\section{Construction of the estimator using the RKHS approach}\label{section2projection}
\noindent  We consider  Model \eqref{model} where $g$ is the unknown regression function on $[0,1]$ and $\{\varepsilon_j(t), t \in [0,1] \}_j$ is a sequence of error processes. We assume that $g \in C^2([0,1])$ and that  $(\varepsilon_j)_j$ are i.i.d. processes with the same distribution as a centered second order process $\varepsilon$. We denote by $R$  its autocovariance function, assumed to be known, continuous and forms a non singular matrix when restricted to $T\times T$ for any finite set $T \subset [0,1]$.

\subsection{Projection estimator}
In this section, we shall give the definition of the new proposed estimator for the regression function $g$ in Model \eqref{model}. This estimator (see Definition \ref{definitionestimator} below) is constructed using the function $f_{x,h}$ given by \eqref{f_xh}  for $x \in  [0,1]$, $h \in  ]0,1[$ and $K$ is a first order kernel\footnote{The kernel $K$ satisfies: $\int_{-1}^1 K(t)dt=1,$ $\int_{-1}^1 t K(t)dt=0$ and $\int_{-1}^1 t^2 K(t)dt <+ \infty.$} of support $[-1,1]$ belonging to $C^1$.
This function is  well known in time series analysis and has been used by several authors. We mention, among others, the works of  \cite{Belouni} and \cite{Sacks and Ylvisaker} for linear regression models with correlated errors. It is mainly used due to its  belonging to the Reproducing Kernel Hilbert Space associated to the autocovariance function $R$ (RKHS(R)) (see Appendix 1 for more details). This space is spanned by the functions $\{R(\cdot,t_i)_{1 \leq i \leq n}\}$ forming a closed subspace on which an orthogonal projection of the function $f_{x,h}$ is feasible. We shall call the estimator obtained by this approach, the projection estimator.\\
The proposed estimator, which is a kernel estimator, is linear in the observations $\overline{Y}(t_i)$ and is given by the following definition.

\begin{definition}\label{definitionestimator}
The projection estimator of the regression function $g$ in Model \eqref{model} based on the observations  $(t_i,Y_j(t_i))_{\underset{1 \leq j \leq m}{1 \leq i \leq n}}$ is given for any $ x \in [0,1]$ by,
\begin{equation}
\hat{g}^{pro}_n (x) = \sum_{i=1}^n m_{x,h} (t_i) \overline{Y} (t_i),
\label{estimator}
\end{equation}
where $\overline{Y}(t_i)=\frac{1}{m}\sum_{j=1}^{m} Y_j(t_i)$ and the weights $(m_{x,h}(t_i))_{1 \leq i \leq n}$ are being determined, letting $T_n=(t_i)_{1 \leq i \leq n}$, by,
\begin{equation}
m_{{x,h}_{|T_n}}'={{f_{x,h}}_{|T_n}}' R_{|T_n}^{-1},
\label{defintionmxh}
\end{equation}
\noindent with $f_{{x,h}_{|T_n}}:=(f_{x,h}(t_1),\dots,f_{x,h}(t_n))'$,  $R_{|T_n}:=(R(t_i,t_j))_{1 \leq i,j \leq n}$, $R_{|T_n}^{-1}$ the inverse of $R_{|T_n}$ and  ${m_{x,h}}_{|T_n}:=(m_{x,h}(t_1),\dots,m_{x,h}(t_n))'$,  where $v'$ denotes the transpose of a vector v.
\end{definition}
\begin{remark} In order to motivate the proposed estimator, consider the regression model using $m$ continuous experimental units, i.e., 
 \begin{equation}\label{modelcontinue}
Y_j(t)=g(t)+\varepsilon_j(t)~\text{for}~t \in [0,1]~~\text{and } j=1,\cdots,m.
\end{equation}
A continuous kernel estimator of $g$ in Model $\eqref{modelcontinue}$ is given for any $x \in [0,1]$ by, \begin{equation}
\hat{g}_{[0,1]}(x)=\int_0^1 \varphi_{x,h}(t)\overline{Y}(t)~dt~~~\text{with}~~\overline{Y}(t)=\frac{1}{m}\sum_{j=1}^{m}Y_j(t), \label{continuousestimator}
\end{equation}
where $\varphi_{x,h}(t)=\frac{1}{h}K\Big( \frac{x-t}{h}\Big)$ for a kernel $K$ and a bandwidth $h$. We refer the reader to the works of 
Blanke and Bosq (2008) \cite{Blanke and Bosq} and Didi and Louani (2013) \cite{Didi and Louani} for more details on the Kernel estimation of the regression function based on continuous observations.\\
Since in practice we only have access to discrete observations,  then a linear approximation of the continuous kernel estimator should be of the form:$$\hat{g}_n (x) = \sum_{i=1}^{n} W_{x,h}(t_i)\overline{Y}(t_i).$$
Now let, $$f_{n,x}(t)=\sum_{i=1}^{n} W_{x,h}(t_i)R(t_i,t)~~~\text{for}~~t\in[0,1].$$
Then the Mean Squared Error (MSE) of  approximation can be written as: $$\mathbb E \big(\hat{g}_{[0,1]}(x)- \hat{g}_n(x) \big)^2=||f_{x,h}-f_{n,x} ||^2,$$
where $f_{x,h}$ is given by \eqref{f_xh} and $||\cdot||$ is the norm of the RKHS(R)(see the Appendix for more details). Then the best linear predictor $\hat{g}^{\text{pro}}_n (x)$ of $\hat{g}_{[0,1]}(x)$ satisfies:
$$ \underset{{W_{x,h}}_{|T_n}}{\inf} \mathbb E \big(\hat{g}_{[0,1]}(x)- \hat{g}_n(x) \big)^2=||f_{x,h}-P_{|T_n}f_{x,h} ||^2,$$
where $P_{|T_n}f_{x,h}$ is the orthogonal projection of $f_{x,h}$ on the subspace of RKHS spanned by the function $\{R(\cdot,t_i), i=1,\cdots,n\}.$ The optimal coefficients $(W_{x,h}^*(t_i))_{1 \leq i \leq n}$ can then be derived by using the fact that $P_{|T_n}f_{x,h}(t_i)=f_{x,h}(t_i)$ for $i=1,\cdots,n$ (see Equation \eqref{appendixPtn=fxh} in the Appendix) and this yields ${W_{{x,h}_{|T_n}}^*}'={{f_{x,h}}_{|T_n}}' R_{|T_n}^{-1}$.
\end{remark}
For some classical error processes, such as the Wiener and the  Ornstein-Uhlenbeck processes, the estimator \eqref{estimator} has a simplified expression as shown in the following proposition.

\begin{proposition}\label{winnerestimatorproposition} Consider the regression model \eqref{model} where $\varepsilon$ is of autocovariance function $R(s,t)=\int_{0}^{\min(s,t)}u^\beta\;du$ for a positive constant $\beta$. Let $t_0=0,~ t_{n+1}=1$. Set $\overline{Y}(t_0)=0$ and $\overline{Y}(t_{n+1})=\overline{Y}(t_n)$. For any $x \in [0,1]$,
	the projection estimator \eqref{estimator} can be written as follows:
		\begin{align}
		\hat{g}^{pro}_n (x) &= \frac{1}{\beta +1}\bigg(\sum_{i=1}^{n+1}\overline{Y} (t_i) \int_{t_{i-1}}^{t_i}\varphi_{x,h}(s)ds  \nonumber \\
	   &~~~~~~	+ \sum_{i=0}^{n-1} \frac{\overline{Y} (t_{i+1})-\overline{Y} (t_i)}{t_{i+1}^{\beta+1}-t_i^{\beta+1}} \int_{t_i}^{t_{i+1}}(s^{\beta+1}-t_{i+1}^{\beta+1})\varphi_{x,h}(s)ds\bigg).
     	\label{winnerestimator}
        \end{align}
\end{proposition}
\begin{remark}
Taking $\beta=0$ in the previous proposition gives the expression of the projection estimator \eqref{estimator} in the case where $\varepsilon$ is the classical standard Wiener error process.
\end{remark}

\begin{proposition}
	If the error process $\varepsilon$  in Model \eqref{model} is the Ornstein-Uhlenbeck process with $R(s,t)=e^{-|t-s|}$ then  for any $x \in [0,1]$,
	\begin{align*}
	&\hat{g}^{pro}_n (x) = \sum_{i=2}^{n-1}\overline{Y}(t_i) \int_{t_{i-1}}^{t_{i+1}} e^{|s-t_i|}\varphi_{x,h}(s)\;ds +\overline{Y}(t_1)\int_{0}^{t_2} e^{s-t_1}\varphi_{x,h}(s)\;ds\\
	&  +\overline{Y}(t_{n})\int_{t_{n-1}}^{1} e^{t_n-s}\varphi_{x,h}(s)\;ds -\sum_{i=1}^{n-1} \frac{e^{t_{i+1}}\overline{Y}(t_{i+1})-e^{t_i}\overline{Y}(t_i)}{1-e^{-2(t_{i+1}-t_i)}} \int_{t_i}^{t_{i+1}} e^{-s}\varphi_{x,h}(s)\;ds\\
	&+\sum_{i=1}^{n-1} \frac{e^{-t_{i+1}}\overline{Y}(t_{i+1})-e^{-t_i}\overline{Y}(t_i)}{1-e^{-2(t_{i+1}-t_i)}} \int_{t_i}^{t_{i+1}} e^{s}\varphi_{x,h}(s)\;ds,
	\end{align*}
	where $\varphi_{x,h}$ is defined in the previous proposition.
	\label{ornsteinestimator}
\end{proposition}
\begin{remark}
As the previous propositions show, the expression of $m_{{x,h}_{|T_n}}$ is known analytically for error processes of practical interest. For more complicated error processes, numerical methods can be used. 
For more general error processes, we will give an asymptotic simplified expression of the weights of the projection estimator (see Lemma \ref{mxhlemma} below).
\end{remark}
\subsection{Assumptions and comments}
In order to derive our asymptotic results, the following assumptions on  the autocovariance function $R$ and the Kernel $K$ are required.
\begin{itemize}
	\item[(A)] $R$ is continuous on the entire unit square and has left and right derivatives up to order two at the diagonal (i.e. when $s=t$), i.e.,  \[
	R^{(0,1)}(t,t^-)=\lim\limits_{s\uparrow t}^{}\frac{\partial R(t,s)}{   \partial s} \quad  \text{and} \quad  R^{(0,1)}(t,t^+)=\lim\limits_{s\downarrow t}\frac{\partial  R(t,s)}{   \partial s},
	\]exist and are continuous. In a similar way we define $R^{(0,2)}(t,t^-)$ and $R^{(0,2)}(t,t^+)$.\\
	Off the diagonal (i.e. when $s \ne t$ in the unit square), $R$ has continuous derivatives up to order two.
\end{itemize}
For $t \in ]0,1[$, let $\alpha(t)=R^{(0,1)}(t,t^-)-R^{(0,1)}(t,t^+)$. Assumption (A) gives the following lemma concerning the jump function $\alpha$.
\begin{lemma}\label{alphapos}
If Assumption $(A)$ is satisfied then the jump function $\alpha$ is a positive function.
\end{lemma}
To obtain our asymptotic results, we shall give next a stronger assumption on the jump function $\alpha$.
\begin{itemize}
\item[(B)]  We assume that $\alpha$ is Lipschitz on $]0,1[$,  $\underset{0< t < 1}{\inf} \alpha (t) = \alpha_0 >0$ and  $\underset{0< t < 1}{\sup} \alpha(t) = \alpha_1 < \infty$.
\end{itemize}
Assumptions (A) and (B) are classical regularity conditions and were used in several works, see for instance, Sacks and Ylvisaker (1966), Su and Cambanis (1993) \cite{Su Cambanis} and most recently Belouni and Benhenni (2015) \cite{Belouni}. 
\begin{itemize}
\item[(C)] For each $t \in [0,1]$, $R^{(0,2)}(.,t^+)$ is in the Reproducing Kernel Hilbert space associated to $R$, denoted by RKHS($R$), equipped with the norm $||\cdot||$.
In addition, $\underset{0 \leq t \leq 1}{\sup}||R^{(0,2)}(.,t^+)||< \infty$  (see the Appendix  for more details).
\end{itemize}
Assumption (C), which is more restrictive than (B) as indicated by Sacks and Ylvisaker (1966) \cite{Sacks and Ylvisaker}, is necessary to evaluate the weights of the projection estimator (see Lemma \ref{mxhlemma} below).
\begin{itemize}
\item[(D)] $K$ is an even function and $K'$ is a Lipschitz function on $[-1,1]$.
\end{itemize}
Examples of autocovariance functions which satisfy  Assumptions $(A)$, $(B)$ and $(C)$ are given below.
\begin{example}~~~~~
	\begin{enumerate}
		\item The autocovariance function $R(s,t) = \sigma^2 min (s,t)$ of the Wiener process, has a constant jump function  $\alpha(t) = \sigma^2 $ and $R^{(i,j)}(s,t) = 0$ for all integers $i,j$ such that $i+j=2$ and $s \ne t$.
		
		\item The autocovariance function $R(s,t)= \sigma^2 e^{-\lambda|s-t|}$ of the stationary Ornstein-Uhnelbeck process  with $\sigma > 0$ and $\lambda > 0$. For this process the jump function is $\alpha(t)=2\sigma^2\lambda$ and $R^{(0,2)}(s,t)=\sigma^2\lambda^2e^{-\lambda|s-t|}$.
		
		\item Another general class of autocovariance functions was given by Sacks and Ylvisaker (1966) \cite {Sacks and Ylvisaker} and has the form,
		\[R(s,t) = \int_{0}^{1/|t-s|} (1-\mu |t-s|)p(\mu)\;d\mu,\]
		where $p$ is a probability density and $p'$ its derivative are such that,
		\[ \underset{\mu \to \infty}{\lim}{\mu^3 p(\mu)} < \infty, ~~~\text{and}~~ \int_{a}^{\infty} (\mu p'(\mu)+3p(\mu))^2)\mu^6 d\mu < \infty,\]
		for some $a$. We have $\alpha(t)=2\int_0^\infty u p(u)~du.$
	\end{enumerate}
\end{example}
\section{Local asymptotic results}\label{localasymptoticresultsprojection}
Let $T_n =(t_{i,n})_{1 \leq i \leq n}$ for $n \geq 1$, be a fixed sequence of designs with $T_n \in D_n$, where, \[D_n = \{(s_1,s_2,\dots,s_n)\;:\; 0 \leq s_1 < s_2 <\dots <s_n \leq 1\}.\] Set $t_{0,n}=0, t_{n+1,n}=1$,  $d_{j,n}= t_{j+1,n}-t_{j,n}$ and let for $x \in [0,1]$, $h=h(n)$,
 \[I_{x,h} =\{i=1,\cdots,n~:~ [t_{i-1,n},t_{i+1,n}]\cap]x-h,x+h[ \ne \emptyset\}.\]
Denote by $N_{T_n} =  \Card (I_{x,h}). $
Recall that $[x-h,x+h]$ is the support of the function $\varphi_{x,h}$.
To obtain the asymptotic results, we require that the sequence $(T_n)_{n \geq 1}$ satisfies the next assumption.
\begin{itemize}
\item[(E)] $\underset{n \to \infty}{\lim} ~\underset{0\leq j \leq n}{\sup}~d_{j,n} = 0$,  $\underset{n \to \infty}{\lim} \Big(\frac{1}{h}\; \underset{0\leq j \leq n}{\sup}~d_{j,n}\Big) = 0$,  $\underset{n \to \infty}{\lim} \Big( N_{T_n} \frac{1}{h^2} ~\underset{0\leq j \leq n}{\sup}\;d_{j,n}^2\Big) = 0$ and

 $\underset{n \to \infty}{\lim \sup}\Big(N_{T_n}^2 \frac{1}{h^2} ~\underset{0\leq j \leq n}{\sup}\;d_{j,n}^2\Big) < \infty $ \textcolor{red}{}.
 \end{itemize}
A simple sequence of designs that verifies Assumption $(E)$ was presented by Sacks and Ylvisaker (1970) \cite{Sacks and Ylvisaker III} as follows.
\begin{definition} \label{regulardesignexample} Let $F$ be a distribution function of some  density function $f$ such that $\underset{0 < t < 1}{ \sup} f(t) < \infty$ and $\underset{0 < t < 1}{ \inf} f(t) > 0$. The so-called regular sequence of designs generated by $f$ is defined by,
	\begin{equation*}
	T_n = \bigg \{t_{i,n}=F^{-1}\bigg(\frac{i}{n}\bigg), i=1,\dots,n \bigg\}.
	\end{equation*}
\end{definition}
	In the sequel, the density $f$ is assumed to be at least in $C^2([0, 1])$.
	This sequence of designs verifies the following Lemma  (see for instance  Benelmadani et al. (2018a) for its proof).
	\begin{lemma} \label{lemmaregulardesign}
		Let $ (T_n)_{n \geq 1}$ be a regular sequence of designs generated by some density function.  For $x \in ]0,1[$ and $h>0$,  suppose that  $T_n \cap [x-h,x+h] \ne \emptyset$ and that $nh \geq 1$. Then,
		\begin{equation}
		\underset{0\leq j \leq n}{\sup}\;d_{j,n} =O\Big(\frac{1}{n}\Big)~~\text{ and }~~N_{T_n}=O(nh)\label{supdjN_T_n},
		\end{equation}
		where $N_{T_n}$ and $d_{j,n}$ are defined as above. In addition, if $\underset{n \to \infty}{\lim} nh = \infty$ then the regular sequence of designs verifies Assumption $(E)$.
		\label{regularverifiesE}
	\end{lemma}
 \subsection{Evaluation of the bias}
In order to derive the asymptotic expression of the bias term of the projection estimator, we shall first give the asymptotic approximation of the   weights  ${m_{x,h}}_{|T_n}$ (defined by \eqref{defintionmxh})  in the following lemma.
\begin{lemma}\label{mxhlemma}
Suppose that  Assumptions $(A),(B)$ and $(C)$ are satisfied. Then for any $x \in ]0,1[$,
\begin {align*}
&~~~~~~~~~~~~~~~~~~~~~~~~~~~~~~~~~~~~ m_{x,h}(t_{i,n}) = \\
& \begin{cases}
	\frac{1}{2} \varphi_{x,h} (t_{i,n}) (t_{i+1,n}-t_{i-1,n})+O\big(\alpha_{n,h}+\beta_{n,h}\big)~~~~~if~~ i \notin \{1,n\} \text{ and }\\
	~~~~~~~~~~~~~~~~~~~~~~~~~~~~~~~~~~~~~~~~~~~~~~~~~~ [t_{i-1,n},t_{i+1,n}]\cap [x-h,x+h] \ne \emptyset ,\\[0.1cm]
	O\big(N_{T_n}\alpha_{n,h}+n\beta_{n,h}\big)~~~~~~~~~~~~~~~~~~~~~~~~~~~~~~~~\text{if}~~ i \in \{1,n\},\\
	O\big(\beta_{n,h}\big)~~~~~~~~~~~~~~~~~~~~~~~~~~~~~~~~~~~~~~~~~~~~~~~otherwise,
\end{cases}
\end{align*}
where,
\begin{align*}
\alpha_{n,h} & = \underset{0 \leq i \leq n}{\sup}~ \underset{t_{i,n} \leq s,t\leq t_{i+1,n}}{\sup}~d_{i,n}~|\alpha(s)\varphi_{x,h}(s)-\alpha(t)\varphi_{x,h}(t)|=O\Big(\frac{1}{h^2}\underset{0 \leq j \leq n}{\sup} d_j^2 \Big),\\
\beta_{n,h}& = \underset{0 \leq t \leq 1}{\sup}~\frac{1}{\alpha(t)}||R^{(0,2)}(.,t)||\frac{\sqrt C}{ \sqrt{h}}~ \underset{0 \leq j \leq n}{\sup} d_j^2=O\Big(\frac{1}{ \sqrt{h}}~ \underset{0 \leq j \leq n}{\sup} d_j^2\Big),
\end{align*}
and $C$ is a positive constant defined in Proposition \ref{propoflimit} below. 
\end{lemma}
\begin{remark}
	This Lemma shows that the weights of the projection estimator are asymptoticly equivalent  to those of some well known linear estimators of the regression function $g$. For instance,
	\begin{itemize}
		\item  Priestly and Chao (1972) \cite{Priestly and Chao} used the following weights: \[W_{x,h}(t_i)=(t_{i+1,n}-t_{i,n})\varphi_{x,h}(t_i)~~~\text{for}~i=1,\cdots,n.\]
		\item  Gasser and Müller (1979) \cite{Gasser and Müller 1979} used the following weights: \[W_{x,h}(t_i)=\int_{s_{i-1,n}}^{s_{i,n}}\varphi_{x,h}(s)~ds~~~~\text{for}~i=1,\cdots,n,\]
		where, $s_0=0$, $s_n=1$ and $s_{i,n}=(t_{i+1,n}+t_{i,n})/2$ for $i=1,\cdots,n-1$.
		\item Cheng and Lin (1981) \cite{Cheng and Lin} replaced $s_{i,n}$ by $t_{i,n}$, in the weights of the Gasser and Müller estimator.
	\end{itemize} 
\end{remark}
Using the asymptotic approximation of the weights given in Lemma \ref{mxhlemma}, we can obtain  the asymptotic expression of the bias of  the projection estimator as shows the following proposition.
\begin{proposition}
	Suppose that Assumptions $(A)-(D)$ are satisfied. If $T_n \cap ]x-h,x+h[ \ne \emptyset$ and $nh \geq 1$,  then for any $x \in ]0,1[$, \[\mathbb{E}(\hat{g}^{pro}_n (x) )  - g(x)=  \frac{1}{2} h^2 g''(x)B + o (h^2) +  O\Big(\frac{N_{T_n}}{h^3} \underset{0\leq j \leq 1}{\sup}\; d_{j,n}^3+ N_{T_n}\alpha_{n,h}+ n\beta_{n,h}\Big), \]
	where $\alpha_{n,h}$ and $\beta_{n,h}$ are given in Lemma \ref{mxhlemma} and $B = \int_{-1}^{1}t^2 K(t)  ~dt$.
	\label{biastheorem}
\end{proposition}
\begin{remark}\label{biasforregulardesign} Under the assumption of Lemma \ref{lemmaregulardesign} we have,
\[\mathbb{E}~(\hat{g}^{pro}_n (x))-g(x)=\frac{1}{2} h^2 g''(x)B + o (h^2) + O\Big(\frac{1}{nh}\Big).\]
\end{remark}
In the case of a  Wiener error process, a direct computation of the bias term of the projection estimator \eqref{winnerestimator}, with $\beta=0$, shows that the order term $O\Big(\frac{1}{nh}\Big)$ can be improved. The result is given by the following proposition.
\begin{proposition} \label{biasofwinner}
	Consider Model \eqref{model} with a Wiener error process of autocovariance function $R(s,t)=\min(s,t)$. Let $(T_n)_{ n \geq 1}$ be a regular sequence of designs generated by a density function $f$ (cf. Definition \ref{regulardesignexample}) and let $K$ be a kernel satisfying Assumption $(D)$. If $T_n \cap ]x-h,x+h[ \ne \emptyset$ and $nh \geq 1$  then,
	\begin{equation*}
	\mathbb{E}~(\hat{g}^{pro}_n (x))-g(x)=
	\frac{1}{2} h^2 g''(x)B+ o (h^2) + O\Big (\frac{1}{n^2h}\Big),
	\end{equation*}
where $B$ is given in Proposition \ref{biastheorem} above.
\end{proposition}
\subsection{Evaluation of the variance}
It is shown in Lemma \ref{finFlemma} of the Appendix  that $f_{x,h}$ defined by \eqref{f_xh} belongs to  the RKHS($R$) equipped with its norm $||~||$ and, \begin{equation}
||f_{x,h}||^2 = \int_0^1 \int_0^1  \varphi_{x,h}(s)R(s,t)\varphi_{x,h}(t) ds~dt\overset{\Delta}{=}\sigma_{x,h}^2.\label{sigmadef}
\end{equation}
In addition if $P_{|T_n}f_{x,h}$ is the projection of $f_{x,h}$ on the subspace of $\mathcal{F}$ spanned by $\{R(.,t),t \in T_n \}$ then it is shown by (F2) in the supplementary facts of the Appendix that,
\begin{equation}
||P_{|T_n}f_{x,h}||^2=m\Var\;\hat{g}_n^{pro}(x).\end{equation}
The following proposition controls the residual variance $\frac{\sigma_{x,h}^2}{m}-\Var\;\hat{g}_n^{pro}(x)$.
\begin{proposition}\label{lemmasigma-var}
Suppose that Assumptions $(A)$ and $(B)$ are satisfied. Moreover, assume that $\frac{1}{h} \underset{1 \leq i \leq n}{\sup}~d_i \leq 1$ and let,
\[K_\infty=\underset{t \in [-1,1]}{\sup} |K(t)|,~~ R_1=\underset{t,s \in [0,1]}{\sup} |R^{(1,1)}(s-,t+)|~~\text{and}~~ R_2=\underset{t,s \in [0,1]}{\sup} |R^{(0,2)}(s,t+)|.\]Then we have, for any $x \in ]0,1[$,
\begin{equation*} 
0 \leq \frac{\sigma_{x,h}^2}{m}-\Var\;\hat{g}_n^{pro}(x)  \leq \frac{C}{mh}~\underset{0\leq j\leq n}{\sup}d_{j,n}^2,
\end{equation*}
\noindent where
$C=\begin{cases}
&K^2_\infty(\frac{4}{3}\alpha_1+R_1+\frac{4}{3}R_2) ~~~~~~~~~~~ \text{ if } (x-h) \text{ and } (x+h) \in T_n,\\
&K^2_\infty(\frac{8}{3}\alpha_1+\frac{5}{3}R_1 + \frac{8}{3}R_2 )~~~~~~~~~\text{ otherwise}.
\end{cases}$\label{propoflimit}
\end{proposition}
If moreover $\{T_n, n \geq 1\}$  satisfies Assumption $(E)$ then  Proposition \ref{propoflimit} gives, \[\underset{n,m \to \infty}{\lim} \Big( \Var\;\hat{g}_n^{pro}(x) -\frac{\sigma_{x,h}^2}{m} \Big) = 0. \]
 The next proposition gives the rate of convergence of this residual variance.
\begin{proposition} \label{variancetheorem}
Suppose that  Assumptions $(A),(B)$ and $(C)$ are satisfied. Moreover, assume that  $(T_n)_{n \geq 1}$ is a sequence of designs verifying Assumption $(E)$. Then for any $x \in ]0,1[$ and for any positive integer $m$, 
\begin{equation}
\underset{n \to \infty}{\underline{\lim}} ~\frac{mN_{T_n}^2}{h} \Big(\frac{\sigma_{x,h}^2}{m}-\Var\;\hat{g}_n^{pro}(x)\Big) \geq \frac{1}{12} \alpha(x) \bigg\{ \int_{-1}^{1} K^{2/3}(t)dt \bigg\}^3,
\label{varianceformula}
\end{equation}
where $\sigma_{x,h}^2$ is given by \eqref{sigmadef}.
\end{proposition}
Using Propositions \ref{propoflimit} and  \ref{variancetheorem} we can obtain the optimal convergence rate $1/(mn^2h)$ of the residual variance. The result is given by the following proposition.
\begin{proposition}\label{optimalrateproposition}
	Suppose that all the assumptions of Lemma \ref{lemmaregulardesign}, Propositions \ref{lemmasigma-var} and  \ref{variancetheorem} are satisfied. Then there exist some positive constants $C$ and $C'$ such that for any $x \in ]0,1[$ and for any positive integer $m$,
	\begin{align}
	\underset{n \to \infty}{\overline \lim }~ mn^2h\Big( \frac{\sigma_{x,h}^2}{m} - \Var(\hat{g}^{pro}_n (x)) \Big) & \leq C, \label{firstineq}
	\end{align}
	and,
	\begin{align}
	\underset{n \to \infty}{\underline \lim }~ mn^2h\Big( \frac{\sigma_{x,h}^2}{m} - \Var~\hat{g}^{pro}_n(x) \Big) & \geq C'\label{secondineq}.
	\end{align}
\end{proposition}
Under the stronger assumption $(D)$ on the kernel $K$ and using a regular sequence of designs (see Definition \ref{regulardesignexample}), we obtain the asymptotic expression of the variance as shown by the following proposition. 
\begin{proposition}\label{exactvarianceformula} Suppose that Assumptions $(A)-(D)$ are satisfied. Moreover assume that $(T_n)_{n \geq 1}$ is a regular sequence of designs generated by a density function $f$ (see Definition \ref{regulardesignexample}). If $\underset{n \to \infty}{\lim} h =0$ and $\underset{n \to \infty}{\lim} nh =\infty$ then for any $x \in ]0,1[$,
	\begin{equation}
	\Var(\hat{g}^{pro}_n (x)) = \frac{\sigma_{x,h}^2}{m}-  \frac{1}{12mn^2}\int_{x-h}^{x+h}\frac{\alpha(t)}{f^2(t)} \varphi_{x,h}^{2}(t)dt + O\bigg(\frac{1}{mn^3h^2}\bigg), \label{varexpress}
	\end{equation}
	where $\sigma_{x,h}^2$ is given by \eqref{sigmadef}.
\end{proposition}
The following lemma (proved in Benhenni and Rachdi (2007) \cite{Benhenni and Rachdi})  gives the expression  of the main term of the asymptotic variance $\sigma_{x,h}^2/m$ in terms of $h$.
\begin{lemma}\label{devesigmalemma}
Suppose that Assumptions $(A)$, $(B)$ and $(D)$ are satisfied. If $\underset{n \to \infty}{\lim} h =0$ then, for any $x \in ]0,1[$, $\sigma_{x,h}^2$ (as given by \eqref{sigmadef}) has the following asymptotic expression
\begin{equation}\label{devesigma}
\sigma_{x,h}^2 = \big( R(x,x)-\frac{1}{2}\alpha(x)C_K h \big) + o(h),
\end{equation}
where $C_K = \int_{-1}^1 \int_{-1}^1 |u-v|K(u)K(v) du dv$.
\end{lemma}
\subsection{IMSE and optimal bandwidth}
Proposition \ref{exactvarianceformula} and Remark \ref{biasforregulardesign} allow to derive the asymptotic expression of the Mean Squared Error ($\MSE$) and the Integrated Mean  Squared Error ($\IMSE$) of the projection estimator \eqref{estimator} given, without proof, in the next theorem.
\begin{theorem}\label{MISEGEN}
	If all the assumptions of Propositions \ref{biastheorem} and \ref{exactvarianceformula} are satisfied and if $(T_n)_{n \geq 1}$ is a regular sequence of designs generated by some density function (see Definition \ref{regulardesignexample}) then for any $x \in ]0,1[$,
	\begin{align*}
	\MSE (\hat{g}^{pro}_n (x)) & = \frac{1}{m}\Big( R(x,x)-\frac{1}{2}\alpha(x)C_K h \Big)+\frac{1}{4} h^4 (g''(x))^2 B^2+o\Big(h^4+\frac{h}{m}\Big) \\
	&~~~~~~+O\Big(\frac{1}{mn^2h}+\frac{h}{n}+\frac{1}{n^2h^2}\Big),
	\end{align*}
	\begin{align*}
	\IMSE(\hat{g}_n^{pro})
	& = \frac{1}{m}\int_0^1 R(x,x)w(x)\;dx-\frac{C_K h}{2m}\int_0^1\alpha(x)w(x)\;dx \\
	&~~~+ \frac{B^2}{4} h^4 \int_0^1 [g''(x)]^2 w(x)\;dx +o\Big(h^4+\frac{h}{m}\Big)+O\Big(\frac{1}{mn^2h}+\frac{h}{n}+\frac{1}{n^2h^2}\Big),
	\end{align*}
	where $w$ is a positive density function, $B$ is given in Proposition \ref{biastheorem} and $ C_K$ is given in  Lemma \ref{devesigmalemma}.
\end{theorem}
\begin{remark}\label{remark4}
We note here that the term $\frac{1}{12mn^2}\int_{x-h}^{x+h}\frac{\alpha(t)}{f^2(t)} \varphi_{x,h}^{2}(t)dt$  appearing in the asymptotic variance, does not appear in the asymptotic $\MSE$ and $\IMSE$,  because it is negligible comparing to the squared bias, precisely due to the term $O\Big(\frac{1}{nh}\Big)$.
	
However in the case of a Wiener error process, we have proven (see Proposition \ref{biasofwinner}) that the previous term can be replaced by $O\Big(\frac{1}{n^2h}\Big)$ when using exact weights of the projection estimator (and not their asymptotic expression). Therefor, when $\varepsilon$ is a Wiener process, the asymptotic expressions of the $\MSE$ and $\IMSE$ of the projection estimator \eqref{winnerestimator} (with $ \beta=0$) are given by the following theorem.
\end{remark}
\begin{theorem}\label{thmise}
	Consider Model \eqref{model} with a Wiener error process and suppose that the kernel $K$ verifies Assumption $(D)$. Moreover, assume that $(T_n)_{n \geq 1}$ is a regular sequence of designs generated by a function $f$ (see Definition \ref{regulardesignexample}). If $\underset{n \to \infty}{\lim} h =0$ and $\underset{n \to \infty}{\lim} nh =\infty$ then for any $x \in ]0,1[$,
	\begin{align*}
	&\MSE(\hat{g}_n^{pro}(x)) = \frac{1}{m}\Big( R(x,x)-\frac{1}{2}\alpha(x)C_K h \Big) - \frac{1}{mn^2h}\frac{\alpha(x)}{f^2(x)}\int_{-1}^1  K^2(t)\;dt \\
	&~~~~+ \frac{1}{4} h^4 [g''(x)]^2 B^2   + o\Big(\frac{h}{m}+h^4\Big)+O\Big(\frac{h}{n²}+\frac{1}{mn^3h^2}+\frac{1}{mn^2}+\frac{1}{n^4h^2}\Big),
	\end{align*}
	and,
	\begin{align*}
	& \IMSE(\hat{g}_n^{pro})= \frac{1}{m}\int_0^1 R(x,x)w(x)\;dx-\frac{C_K h}{2m}\int_0^1\alpha(x)w(x)\;dx \\
	&-  \frac{A}{12mn^2h}\int_0^1\frac{\alpha(x)}{f^2(x)}w(x)\;dx
	+ \frac{B^2}{4} h^4 \int_0^1 [g''(x)]^2 w(x)\;dx   + o\Big(\frac{h}{m}+h^4\Big)\\
	&+O\Big(\frac{h}{n²}+\frac{1}{mn^3h^2}+\frac{1}{mn^2}+\frac{1}{n^4h^2}\Big),
	\end{align*}
	where $A=\int_{-1}^1  K^2(t)\;dt$, $w$, $B$ and $ C_K$ are given in Theorem \ref{MISEGEN}.
\end{theorem}
The asymptotic optimal bandwidth is obtained by minimizing  the asymptotic IMSE and is given in the following corollary.
\begin{corollary}[Optimal bandwidth]\label{hoptcor}
	Suppose that the assumptions of Theorem \ref{MISEGEN} are satisfied. Moreover assume that $\frac{n}{m} = O(1)$ as $n,m \to \infty$. Denote
	by IMSE(h) the IMSE of the projection estimator when the bandwidth h is
	used.  Then the bandwidth,
	\begin{equation} \label{hoptimal}
	h^* = \Bigg( \frac{C_K\int_0^1 \alpha(x)w(x)\;dx}{2B\int_0^1 [g''(x)]^2 w(x)\;dx}\Bigg)^{1/3}m^{-1/3},
	\end{equation}
	 is optimal in the sense that, \[\underset{n,m \to \infty}{\overline{\lim}} \frac{\IMSE(h^*)}{\IMSE(h_{n,m})}   \leq 1,\]
	for any sequence of bandwidths $h_{n,m}$ verifying:
	\[\underset{n,m \to \infty}{\lim} h_{n,m} = 0 ~~~\text{and}~~\underset{n,m \to \infty}{\overline\lim} mh_{n,m}^3<+\infty.\]
\end{corollary}
 \subsection{Asymptotic normality}
  The next theorem presents the asymptotic normality of the projection estimator \eqref{estimator} for any error process $\varepsilon$.
\begin{theorem}\label{asymptotic normality theorem}
	Suppose that the assumptions of  Theorem \ref{MISEGEN} are satisfied. Moreover assume that $\frac{n}{m} = O(1)$ as $n,m \to \infty$, that $\underset{n, m \to \infty}{\lim} n h^2= \infty $  and  that $\underset{n, m \to \infty}{\lim} \sqrt m h^2=0 $. Then for any $x \in ]0,1[$,
	\[ \sqrt m \Big( {\hat{g}^{pro}_n}(x) - g(x) \Big) \overset{\mathscr D}{\longrightarrow} Z~~~\text{with}~Z \sim \mathcal N (0,{R(x,x)})~~~\text{as}~n,m \to \infty, \]
	where $\mathscr D$ denotes the convergence in distribution and $\mathcal N$ is the normal distribution.
\end{theorem}
\section{Comparison with the Gasser and Müller's estimator}\label{comparisonGM}
In this section, we shall perform a theoretical comparison between the projection estimator given in \eqref{estimator} and the classical estimator proposed by Gasser and Müller (1979) \cite{Gasser and Müller 1979} that we recall in the definition below. 
\begin{definition}
	The Gasser and Müller's estimator of the regression function $g$ based on the observations $(t_i,Y_j(t_i))_{\underset{1 \leq j \leq m}{1 \leq i \leq n}}$ is given for any $x \in [0,1]$ by,
	\begin{equation}
	\hat{g}_n^{\text{GM}}(x) = \sum_{i=1}^n\overline{Y}(t_i) \int_{s_{i-1}}^{s_i} \varphi_{x,h}(s)\;ds\; , \label{Müllerestimatorformula}
	\end{equation}
	where $\overline{Y}, \varphi_{x,h}$ and $h$ are given in Definition \ref{definitionestimator}. The  midpoints $(s_i)_{1 \leq i \leq n}$ are such that: $s_0=0, s_n=1$  and for $i=1,\ldots,n-1$, $s_i=(t_i+t_{i+1})/2$.
\end{definition}
In order to compare this estimator to the projection estimator with respect to the IMSE, we recall in the next theorem the asymptotic expression of the IMSE of the Gasser and Müller's estimator (for the proof see Benelmadani et al. (2019) \cite{Benelmadani and Benhenni and Louhichi and Su} and Benhenni and Rachdi (2007) \cite{Benhenni and Rachdi} for further detailed results).
\begin{theorem}\label{MSEofMüller}
	Suppose that Assumptions $(A)$, $(B)$ and $(D)$ are satisfied. Moreover assume that $(T_n)_{n \geq 1}$ is a regular sequence of designs generated by a density function $f$ (see Definition \ref{regulardesignexample}). If $\underset{n \to \infty}{\lim} h =0$ and $\underset{n \to \infty}{\lim} nh =\infty$ then for any $x \in ]0,1[$,
	\begin{align*}
	\MSE (\hat{g}^{GM}_n (x))  &= \frac{1}{m}\Big( R(x,x)-\frac{1}{2}\alpha(x)C_K h \Big)+\frac{1}{4} h^4 (g''(x))^2 B^2 +o\Big(h^4+\frac{h}{m}\Big)\\
	&+ O\Big(\frac{h}{n^2}+\frac{1}{n^4h^2}+\frac{1}{mn^3h^2}+\frac{1}{mn^2}\Big),
	\end{align*}
	and,
	\begin{align*}
	\IMSE(\hat{g}_n^{GM}) &= \frac{1}{m}\int_0^1 R(x,x)w(x)\;dx-\frac{C_K h}{2m}\int_0^1\alpha(x)w(x)\;dx+ \frac{B^2}{4} h^4 \int_0^1 [g''(x)]^2 w(x)\;dx\\
	&~~ +o\Big(h^4+\frac{h}{m}\Big)+ O\Big(\frac{h}{n^2}+\frac{1}{n^4h^2}+\frac{1}{mn^3h^2}+\frac{1}{mn^2}\Big),
	\end{align*}
	where $B$ and $C_K$ are given in Propositions \ref{biastheorem} and \ref{exactvarianceformula} and $w$ is a continuous positive density.
\end{theorem}
The following theorem gives an asymptotic comparison in term of the variance of the projection estimator \eqref{estimator} and  the Gasser and Müller's estimator \eqref{Müllerestimatorformula}.
\begin{theorem}\label{compavariancegen}
	Suppose that Assumptions $(A), (B)$ and $(D)$ are satisfied. Moreover assume that $(T_n)_{n \geq 1}$ is a regular sequence of designs generated by a density function $f$ (see Definition \ref{regulardesignexample}). If $\underset{n \to \infty}{\lim} h =0$ and $\underset{n \to \infty}{\lim} nh =\infty$ then for any $x \in ]0,1[$,
	\begin{align*}
	\underset{n, m \to \infty}{\lim }~~ mn^2h ~\Big(\Var\;\hat{g}_n^{GM}(x)-\Var\;\hat{g}_n^{pro}(x) \Big)  = \frac{1}{12} \frac{\alpha(x)}{f^2(x)}> 0.
	\end{align*}
\end{theorem}
For a comparison of the bias of these estimators, we mention that the Gasser and Müller's estimator converges to zero slightly faster than the bias of the projection estimator, i.e., the term $O(\frac{1}{nh})$ in the bias of the projection estimator (see Remark \ref{biasforregulardesign}) is replaced by $O(\frac{1}{n^2h})$ in the bias of the Gasser and Müller's estimator (see Benelmadani at al. (2019) \cite{Benelmadani and Benhenni and Louhichi and Su}). However, for the Wiener error process both estimators have the same bias convergence rates, thus we can compare the asymptotic IMSE of both estimators in the following theorem.
\begin{theorem}\label{compavariance}
	Consider Model \eqref{model} where $\varepsilon$ is a Wiener error process. Suppose that the assumptions of Theorem \ref{thmise} are satisfied. Moreover, assume that  $\underset{n \to \infty}{\lim} nh^2 = 0$ and that $\frac{m}{n}=O(1)$ then,
	\begin{align*}
	\underset{n, m \to \infty}{ \lim }~~ mn^2h ~(\IMSE\;(\hat{g}_n^{GM})-\IMSE\;(\hat{g}_n^{pro}) )  = \frac{\sigma^2}{12}\int_{0}^{1} \frac{w(x)}{f^2(x)}~dx > 0.
	\end{align*}
\end{theorem}
\begin{remark} \text{}
	Theorems \ref{compavariancegen} and \ref{compavariance} show that, the projection estimator has an asymptotically smaller variance than the Gasser and Müller's estimator for any error process, it also has an asymptotically smaller IMSE when $\varepsilon$ is a Wiener error process.  However the Gasser  and Müller's estimator doesn't require the knowledge of the autocovariance function whereas the projection estimator does.
\end{remark}
\section{Simulation study}\label{simulationprojection}
\begin{figure}[H]
	\begin{minipage}[c]{.3\linewidth}
		\begin{center}
			\includegraphics[scale=0.45]{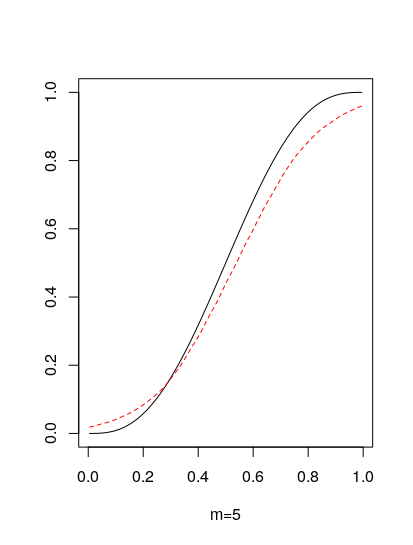}
		\end{center}
	\end{minipage} \hfill
	\begin{minipage}[c]{.3\linewidth}
		\begin{center}
			\includegraphics[scale=0.45]{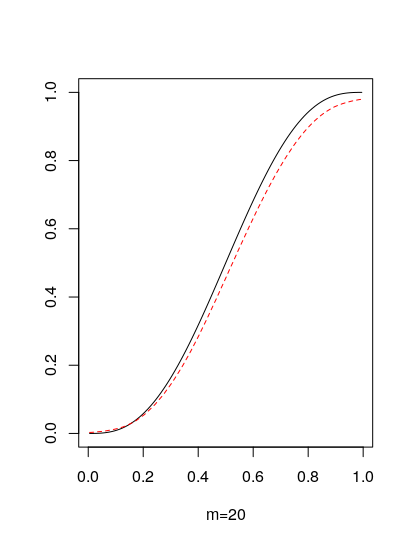}
		\end{center}
	\end{minipage} \hfill
	\begin{minipage}[c]{.3\linewidth}
		\begin{center}
			\includegraphics[scale=0.45]{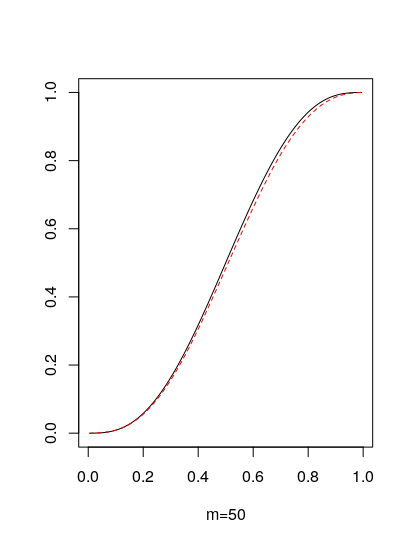}
		\end{center}
	\end{minipage}
	\caption{The regression function of model (M1)  is in plain line and the projection estimator is in dashed line.}
	\label{graphsdeprojection1}
\end{figure}

\begin{figure}[H]
	\begin{minipage}[c]{.3\linewidth}
		\begin{center}
			\includegraphics[scale=0.55]{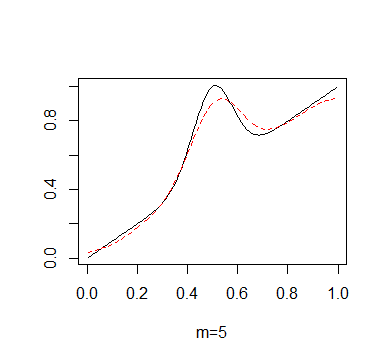}
		\end{center}
	\end{minipage} \hfill
	\begin{minipage}[c]{.3\linewidth}
		\begin{center}
			\includegraphics[scale=0.55]{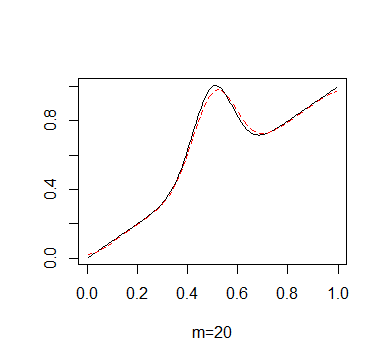}
		\end{center}
	\end{minipage} \hfill
	\begin{minipage}[c]{.3\linewidth}
		\begin{center}
			\includegraphics[scale=0.55]{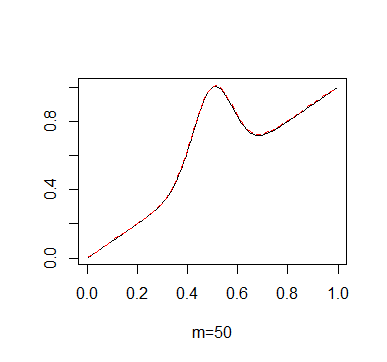}
		\end{center}
	\end{minipage}
	\caption{The regression function of model (M2)  is in plain line and the projection estimator is in dashed line.}
	\label{graphsdeprojection}
\end{figure}

In this section, we  investigate the performance of the proposed estimator \eqref{estimator} using  finite values of experimental units $m$ and  sampling points $n$. The following growth curves are considered: 
\begin{align*}
\text{(M1)}~~~~~~~~g(x)&=10x^3-15x^4+6x^5~~~\text{for}~~0<x<1.\\
\text{(M2)}~~~~~~~~g(x)&=x+0.5~e^{-80(x-0.5)^2}~~~\text{for}~~0<x<1.
\end{align*}
This growth curves were used by Hart and Wherly (1986) \cite{Hart and Wherly} and Benhenni and Rachdi (2007) \cite{Benhenni and Rachdi}  due to its similarity in shape to that of the logistic function, which is frequently found in growth curve analysis as noted by Hart and Wherly (1986) \cite{Hart and Wherly}.
The sampling points are taken to be:
\begin{equation}
t_i=(i-0.5)/n~~\text{for}~i=1,\cdots,n. \label{designinsimulation}
\end{equation}
The error process $\varepsilon$ is taken to be the Wiener error process with autocovariance function $R(s,t)=\sigma^2\min(s,t)$.
The Kernel used here is the quartic  kernel given by $K(u)=\frac{15}{16}(1-u^2)^2I_{[-1,1]}(u)$ and the bandwidth  is the optimal one with respect to the exact $\IMSE$, obtained using the Conjugated Gradiant Algorithm (CGA).
We consider the mean of all estimators obtained from 100 simulations. We take $\sigma^2=0.5$, simulations for other values of $\sigma^2$ gave similar results. The results are given in Figures \ref{graphsdeprojection1} and \ref{graphsdeprojection} for Models (M1) and (M2) respectively, for  a fixed number of observations $n=100$ and three different values of experimental units $m=5, 20, 50$.
\noindent We can see for Model (M1), from Figure \ref{graphsdeprojection1}, that the projection estimator gets closer to the regression function when $m$ gets bigger, which proves its good performance and consistency when $m$ increases. These results are confirmed for the growth curve  Model (M2) given in Figure \ref{graphsdeprojection}.
 
In this simulation study, we consider the comparison of the proposed estimator \eqref{estimator} to the Gasser and Müller \eqref{Müllerestimatorformula} (referred by GM estimator) with respect to the exact IMSE in a finite sample set. For this, we consider the cubic growth curve of model (M1).
We consider also the uniform design given by \eqref{designinsimulation} and the quartic kernel $K(u)=\frac{15}{16}(1-u^2)^2I_{[-1,1]}(u)$. For the error process, we shall consider both the Wiener of autocovariance function $R(s,t)=\min(s,t)$, and the Ornstein-Uhlenbeck process with autocovariance $R(s,t)=e^{-|s-t|}$.

The weight $w$, chosen here, is the uniform density on $[0,1]$, i.e., $w \equiv 1$ on $[0,1]$, we consider the optimal bandwidth with respect to the exact $\IMSE$ of the two estimators, i.e., $\underset{0 < h< 1}{\inf}~ \IMSE (h)$. The bandwidth $h$ is chosen through the algorithm CGA. The results are given in Tables \ref{wiener} and \ref{ornstein} for  $n=10$  and for different values of $m$. These tables present the integrated bias squared denoted by $Ibias^2$, integrated variance denoted by  $Ivar$ and the $\IMSE$ together with the optimal bandwidth  associated to each estimator.

First, we can see from these two tables that, the optimal bandwidth decreases when $m$ increases, as shown in Corollary \ref{hoptcor}. In addition, the optimal bandwidth of the projection estimator is slightly smaller than that of the GM estimator.

It is also seen that both the $Ivar$ and the $Ibias^2$, of the two estimators decrease when $m$ increases. In addition, the projection estimator has a smaller $Ibias^2$ and $Ivar$ than that of the GM estimator, which leads to a smaller $\IMSE$.

Another way to look at these results is as follows: for a fixed  number of experimental units  $m=10$ and when the error process is a Wiener process (similar results for the Ornstein-Uhlenbeck error process), the projection estimator would only need $n=10$ observations on each experimental unit to obtain the performance $\IMSE=4.53 \times 10^{-2}$ (see Table \ref{wiener}), whereas the GM estimator would need to have $n=18$ observations to obtain the same performance, and thus requires $80\%$ more samples in order to achieve the same performance.

The results of this simulation study show that, even for small number of observations, the projection estimator outperforms the GM estimator with respect to IMSE.

It should be noted here that, in order to solve the problem at the edges $[0,h]\cap[1-h,1]$, it was necessary to adjust the kernel as suggested by Hart and Wherly (1986) \cite{Hart and Wherly}.

\begin{table}[H]
	\caption{The integrated squared bias, integrated variance, IMSE and the optimal bandwidth for $n=10$ and different values of $m$ under the Wiener error process, for the GM and the projection estimators.} \label{wiener}
	\begin{center}
		\begin{tabular}{|c|c|c|c|c|c|}
			\hline
			$n=10$ &$m$ & $Ibias^2$ & $Ivar$ & $\IMSE$ & $h_{opt}$ \\
			\hlinewd{1.5pt}
			$GM$ &\multirow{2}{*}{10}& $1.508\times 10^{-3}$ & $4.507\times 10^{-2}$ & $4.658\times 10^{-2}$ &0.335\\
			\cline{1-1}\cline{3-6}
			$Pro$&& $1.304\times 10^{-3}$ &$4.399\times 10^{-2}$ &$4.530\times 10^{-2}$ &0.321 \\
			\hlinewd{1.5pt}

			$GM$&\multirow{2}{*}{50} &$2.662\times 10^{-4}$ &$9.494\times 10^{-3}$& $9.760\times 10^{-3}$ &0.198\\
			\cline{1-1}\cline{3-6}
			$Pro$& & $1.981\times 10^{-4}$ &$9.228\times 10^{-3}$ &$9.426\times 10^{-3}$ &0.187 \\
			\hlinewd{1.5pt}
			
			$GM$&\multirow{2}{*}{100} & $1.505\times 10^{-4}$ &$4.826\times 10^{-3}$ &$4.977\times 10^{-3}$ &0.154\\
			\cline{1-1}\cline{3-6}
			$Pro$& & $0.897\times 10^{-4}$ &$4.689\times 10^{-3}$ &$4.778\times 10^{-3}$ &0.142 \\
			\hline
		\end{tabular}
	\end{center}
	\label{tab:multicol}
\end{table}

\begin{table}[H]
	\caption{The integrated squared bias, integrated variance, IMSE and the optimal bandwidth for $n=10$ and different values of $m$ under the Ornstein-Uhlenbeck error process, for the GM and the projection estimators.} \label{ornstein}
	\begin{center}
		\begin{tabular}{|c|c|c|c|c|c|}
			\hline
			$n=10$ &$m$ & $Ibias^2$ & $Ivar$ & $\IMSE$ & $h_{opt}$ \\
			\hlinewd{1.5pt}
			$GM$ &\multirow{2}{*}{10}& $2.596\times 10^{-3}$ & $8.821\times 10^{-2}$ & $9.080\times 10^{-2}$ &0.387\\
			\cline{1-1}\cline{3-6}
			$Pro$&& $2.494\times 10^{-3}$ &$8.703\times 10^{-2}$ &$8.952\times 10^{-2}$ & 0.386 \\
			\hlinewd{1.5pt}

			$GM$&\multirow{2}{*}{50} &$4.481\times 10^{-4}$ &$1.848\times 10^{-2}$& $1.893\times 10^{-2}$ &0.236\\
			\cline{1-1}\cline{3-6}
			$Pro$& & $4.097\times 10^{-4}$ &$1.822\times 10^{-2}$ &$1.863\times 10^{-2}$ &0.237 \\\hlinewd{1.5pt}

			$GM$&\multirow{2}{*}{100} & $2.299\times 10^{-4}$ &$9.390\times 10^{-3}$ &$9.620\times 10^{-3}$ & 0.186 \\
			\cline{1-1}\cline{3-6}
			$Pro$& & $1.885\times 10^{-4}$ &$9.265\times 10^{-3}$ &$9.453\times 10^{-3}$ & 0.187 \\
			\hline
		\end{tabular}
	\end{center}
	\label{tab:multicol}
\end{table}
\section{Proofs}In this section, we shall omit the index $n$ in $t_{i,n}$ when there is no ambiguity.
\subsection{Proof of Proposition \ref{winnerestimatorproposition}.}
It is known that (see, for instance Su and Cambanis (1993) \cite{Su Cambanis} page 88) if $R(s,t) = \int_{0}^{\min(s,t)} u^\beta \;du$ then for any functions $u$ and $v$ and for any sampling design $T_n$ we have,
\begin{equation*}
{u_{|T_n}}'R_{|T_n}^{-1} v_{|T_n} = \frac{u(t_1)v(t_1)}{t_1^{\beta+1}}+\sum_{k=1}^{n-1} \frac{(u(t_{k+1})-u(t_k))(v(t_{k+1})-v(t_k))}{t_{k+1}^{\beta+1}-t_k^{\beta+1}}.
\end{equation*}
Replacing $u=f_{x,h}$ and $v=\overline{Y}$ we have,
\begin{equation}
\hat{g}_n^{pro} (x) = \frac{f_{x,h}(t_1)\overline{Y}(t_1)}{t_1^{\beta+1}}+\sum_{i=1}^{n-1} \frac{(f_{x,h}(t_{i+1})-f_{x,h}(t_i))(\overline{Y}(t_{i+1})-\overline{Y}(t_i))}{t_{i+1}^{\beta+1}-t_i^{\beta+1}}.
\end{equation}
Recall that $R(s,t)=\frac{1}{\beta+1} \min(s,t)^{\beta+1}$ and,
\begin{align*}
f_{x,h}(t_i) &=\int_0^1 R(s,t_i) \varphi_{x,h}(s)\;ds
=\frac{1}{\beta+1}\Big(\int_{0}^{t_i} s^{\beta+1}~ \varphi_{x,h}(s)\;ds +t_i^{\beta+1}~\int_{t_i}^{1}\varphi_{x,h}(s)\;ds\Big).
\end{align*}
Thus,
\begin{align*}
& f_{x,h}(t_{i+1})-f_{x,h}(t_{i})  = \frac{1}{\beta+1}\Big(\int_{0}^{t_{i+1}} s^{\beta+1} \varphi_{x,h}(s)~ds +  t_{i+1}^{\beta+1} \int_{t_{i+1}}^{1} \varphi_{x,h}(s)~ds \\
&-\int_{0}^{t_{i}} s^{\beta+1} \varphi_{x,h}(s)~ds- t_{i}^{\beta+1}  \int_{t_{i}}^{1}\varphi_{x,h}(s)~ds+t_{i+1}^{\beta+1}  \int_{t_{i}}^{1}\varphi_{x,h}(s)~ds-t_{i+1}^{\beta+1}  \int_{t_{i}}^{1}\varphi_{x,h}(s)~ds\Big)\\
& =\frac{1}{\beta+1}\Big( \int_{t_{i}}^{t_{i+1}} (s^{\beta+1}-t_{i+1}^{\beta+1}) \varphi_{x,h}(s)~ds  +(t_{i+1}^{\beta+1}-t_{i}^{\beta+1})\int_{t_{i}}^{1} \varphi_{x,h}(s)~ds\Big).
\end{align*}
Thus,
\begin{align*}
\hat{g}_n^{pro}(x) & = \frac{f_{x,h}(t_1)\overline{Y}(t_1)}{t_1^{\beta+1}}+\frac{1}{\beta+1}\Big(\sum_{i=1}^{n-1} (\overline{Y}(t_{i+1})-\overline{Y}(t_i))\int_{t_{i}}^{1} \varphi_{x,h}(s)~ds \\
&~~~+\sum_{i=1}^{n-1} \frac{\overline{Y}(t_{i+1})-\overline{Y}(t_i)}{t_{i+1}^{\beta+1}-t_i^{\beta+1}} \int_{t_{i}}^{t_{i+1}} (s^{\beta+1}-t_{i+1}^{\beta+1}) \varphi_{x,h}(s)~ds\Big)\\
& =  \frac{f_{x,h}(t_1)\overline{Y}(t_1)}{t_1^{\beta+1}}+\frac{1}{\beta+1}\Big(\sum_{i=2}^{n-1} \overline{Y}(t_i)\int_{t_{i-1}}^{t_i}\varphi_{x,h}(s)~ds-\overline{Y}(t_1)\int_{t_{1}}^{1} \varphi_{x,h}(s)~ds\\
&~~~ +\overline{Y}(t_{n})\int_{t_{n-1}}^{1} \varphi_{x,h}(s)~ds +\sum_{i=1}^{n-1} \frac{\overline{Y}(t_{i+1})-\overline{Y}(t_i)}{t_{i+1}^{\beta+1}-t_i^{\beta+1}} \int_{t_{i}}^{t_{i+1}} (s^{\beta+1}-t_{i+1}^{\beta+1}) \varphi_{x,h}(s)~ds\Big).
\end{align*}
Letting $t_0=\overline{Y}(t_{0})=0$ we have,
\begin{align*}
\frac{f_{x,h}(t_1)\overline{Y}(t_1)}{t_1^{\beta+1}}  &=\frac{1}{\beta+1}\Big( \frac{\overline{Y}(t_{1})}{t_1^{\beta+1}}\int_0^{t_1}s^{\beta+1}\varphi_{x,h}(s)~ds+\overline{Y}(t_{1})\int_{t_1}^1\varphi_{x,h}(s)~ds\Big)\\
& =\frac{1}{\beta+1}\Big( \frac{\overline{Y}(t_{1})-\overline{Y}(t_{0})}{t_1^{\beta+1}-t_0^{\beta+1}}\int_0^{t_1}(s^{\beta+1}-t_1^{\beta+1})\varphi_{x,h}(s)~ds\\
&~~~~+\overline{Y}(t_{1})\int_{0}^{t_{1}}\varphi_{x,h}(s)~ds+\overline{Y}(t_{1})\int_{t_1}^1\varphi_{x,h}(s)~ds\Big).
\end{align*}
Finally,
\begin{align*}
&\hat{g}_n^{pro}(x) = \frac{1}{\beta+1}\Big(\sum_{i=2}^{n-1} \overline{Y}(t_i)\int_{t_{i-1}}^{t_i}\varphi_{x,h}(s)~ds-\overline{Y}(t_1)\int_{t_{1}}^{1} \varphi_{x,h}(s)~ds +\overline{Y}(t_{n})\int_{t_{n-1}}^{1} \varphi_{x,h}(s)~ds \\
& +\overline{Y}(t_{n})\int_{t_{n-1}}^{t_n} \varphi_{x,h}(s)~ds -\overline{Y}(t_{n})\int_{t_{n-1}}^{t_n} \varphi_{x,h}(s)~ds \\
&+\sum_{i=1}^{n-1} \frac{\overline{Y}(t_{i+1})-\overline{Y}(t_i)}{t_{i+1}^{\beta+1}-t_i^{\beta+1}} \int_{t_{i}}^{t_{i+1}} (s^{\beta+1}-t_{i+1}^{\beta+1}) \varphi_{x,h}(s)~ds\\
&+\frac{\overline{Y}(t_{1})-\overline{Y}(t_{0})}{t_1^{\beta+1}-t_0^{\beta+1}}\int_0^{t_1}(s^{\beta+1}-t_1^{\beta+1})\varphi_{x,h}(s)~ds+\overline{Y}(t_{1})\int_{0}^{t_{1}}\varphi_{x,h}(s)~ds+\overline{Y}(t_{1})\int_{t_1}^1\varphi_{x,h}(s)~ds\Big)\\
& = \frac{1}{\beta+1}\Big(\sum_{i=1}^{n+1} \overline{Y}(t_i)\int_{t_{i-1}}^{t_i}\varphi_{x,h}(s)~ds+\sum_{i=0}^{n-1} \frac{\overline{Y}(t_{i+1})-\overline{Y}(t_i)}{t_{i+1}^{\beta+1}-t_i^{\beta+1}} \int_{t_{i}}^{t_{i+1}} (s^{\beta+1}-t_{i+1}^{\beta+1}) \varphi_{x,h}(s)~ds\Big),
\end{align*}
where $t_{n+1}=1$ and $\overline{Y}(t_{n+1})=\overline{Y}(t_{n})$. This concludes the proof of Proposition \ref{winnerestimatorproposition}. $\Box$
\subsection{Proof of Proposition \ref{ornsteinestimator}.}
It is known (see Anderson (1960) \cite{TW Anderson}  page 210) that for every functions $u$ and $v$ and for every design $T_n$ we have,
\begin{align*}
u_{|T_n}' R_{|T_n}^{-1}v_{|T_n} & = \frac{u(t_1)v(t_1)}{1-e^{-2(t_2-t_1)}} +\frac{u(t_n)v(t_n)}{1-e^{-2(t_n-t_{n-1})}} + \sum_{i=2}^{n-1} \frac{u(t_i)v(t_i) (1-e^{-2(t_{i+1}-t_{i-1})})}{(1-e^{-2(t_{i+1}-t_i)})(1-e^{-2(t_i-t_{i-1})})}\\
& ~~~~ - \sum_{i=1}^{n-1} \frac{u(t_i)v(t_{i+1}) +u(t_{i+1})v(t_i) }{1-e^{-2(t_{i+1}-t_i)}} e^{-(t_{i+1}-t_i)}.
\end{align*}
Taking $u = f_{x,h}$ and $v = \overline{Y}$  we get,
\begin{align}
\hat{g}_n^{pro}(x)& = \frac{f_{x,h}(t_1)\overline{Y}(t_1)}{1-e^{-2(t_2-t_1)}} +\frac{f_{x,h}(t_n)\overline{Y}(t_n)}{1-e^{-2(t_n-t_{n-1})}} + \sum_{i=2}^{n-1} \frac{f_{x,h}(t_i)\overline{Y}(t_i) (1-e^{-2(t_{i+1}-t_{i-1})})}{(1-e^{-2(t_{i+1}-t_i)})(1-e^{-2(t_i-t_{i-1})})}\nonumber\\
& ~~~~ - \sum_{i=1}^{n-1} \frac{f_{x,h}(t_i)\overline{Y}(t_{i+1}) +f_{x,h}(t_{i+1})\overline{Y}(t_i) }{1-e^{-2(t_{i+1}-t_i)}} e^{-(t_{i+1}-t_i)}\nonumber\\
& \overset{\Delta}{=} \frac{f_{x,h}(t_1)\overline{Y}(t_1)}{1-e^{-2(t_2-t_1)}} +\frac{f_{x,h}(t_n)\overline{Y}(t_n)}{1-e^{-2(t_n-t_{n-1})}} +A. \label{gnpro}
\end{align}
Note that,
\[1-e^{-2(t_{i+1}-t_{i-1})}=(1-e^{-2(t_{i+1}-t_{i})})+(1-e^{-2(t_{i}-t_{i-1})})-(1-e^{-2(t_{i}-t_{i-1})})(1-e^{-2(t_{i+1}-t_{i})}).\]
Thus,
\begin{align}
 A  &= \sum_{i=2}^{n-1} \frac{f_{x,h}(t_i)\overline{Y}(t_i)}{1-e^{-2(t_i-t_{i-1})}}+\sum_{i=2}^{n-1} \frac{f_{x,h}(t_i)\overline{Y}(t_i) }{1-e^{-2(t_{i+1}-t_i)}}-\sum_{i=2}^{n-1} f_{x,h}(t_i)\overline{Y}(t_i) \nonumber\\
& ~~~- \sum_{i=2}^{n} \frac{f_{x,h}(t_{i-1})\overline{Y}(t_{i})}{1-e^{-2(t_{i}-t_{i-1})}} e^{-(t_{i}-t_{i-1})}-\sum_{i=1}^{n-1} \frac{f_{x,h}(t_{i+1})\overline{Y}(t_i) }{1-e^{-2(t_{i+1}-t_i)}} e^{-(t_{i+1}-t_i)} \nonumber\\
& = \sum_{i=2}^{n-1}\frac{\overline{Y}(t_i)}{1-e^{-2(t_i-t_{i-1})}}\Big(f_{x,h}(t_i)-f_{x,h}(t_{i-1})e^{-(t_{i}-t_{i-1})} \Big)-\frac{f_{x,h}(t_{n-1})\overline{Y}(t_{n})}{1-e^{-2(t_{n}-t_{n-1})}} e^{-(t_{n}-t_{n-1})}\nonumber\\
& ~~~+ \sum_{i=2}^{n-1}\frac{\overline{Y}(t_i)}{1-e^{-2(t_{i+1}-t_{i})}}\Big(f_{x,h}(t_i)-f_{x,h}(t_{i+1})e^{-(t_{i+1}-t_{i})} \Big)-\frac{f_{x,h}(t_{2})\overline{Y}(t_1)}{1-e^{-2(t_{2}-t_{1})}} e^{-(t_{2}-t_{1})}\nonumber\\
&~~~~-\sum_{i=2}^{n-1} f_{x,h}(t_i)\overline{Y}(t_i) \label{A+B}
\end{align}
Simple calculations yield,
\begin{align}
&f_{x,h}(t_i)- f_{x,h}(t_{i-1})e^{-(t_i-t_{i-1})}=\nonumber\\
&  e^{-t_i}\int_{t_{i-1}}^{t_i}  e^{s}\varphi_{x,h}(s)\;ds - e^{t_i}\int_{t_{i-1}}^{t_i}  e^{-s}\varphi_{x,h}(s)\;ds+e^{t_i}(1-e^{-2(t_i-t_{i-1})})\int_{t_{i-1}}^{1}  e^{-s}\varphi_{x,h}(s)\;ds. \label{dif1ornstein}
\end{align}
In the same way we have,
\begin{align}
&f_{x,h}(t_i)-f_{x,h}(t_{i+1})e^{-(t_{i+1}-t_{i})}=\nonumber\\
&  e^{t_i}\int_{t_{i}}^{t_{i+1}}  e^{-s}\varphi_{x,h}(s)\;ds - e^{-t_i}\int_{t_{i}}^{t_{i+1}}  e^{s}\varphi_{x,h}(s)\;ds+e^{-t_i}(1-e^{-2(t_{i+1}-t_{i})})\int_{0}^{t_{i+1}}  e^{s}\varphi_{x,h}(s)\;ds.
\label{dif2ornstein}
\end{align}
It is easy to verify that,
\begin{align}
\sum_{i=2}^{n-1} f_{x,h}(t_i)\overline{Y}(t_i)=\sum_{i=2}^{n-1} \overline{Y}(t_i)e^{-t_i}\int_{0}^{t_i}  e^{s}\varphi_{x,h}(s)\;ds+\sum_{i=2}^{n-1} \overline{Y}(t_i)e^{t_i}\int_{t_{i}}^{1}  e^{-s}\varphi_{x,h}(s)\;ds.\label{mas}
\end{align}
We obtain using Equations \eqref{A+B}, \eqref{dif1ornstein}, \eqref{dif2ornstein} and \eqref{mas},
\begin{align*}
& A  = \sum_{i=2}^{n-1}\overline{Y}(t_i)e^{t_i}\int_{t_{i-1}}^{1}  e^{-s}\varphi_{x,h}(s)\;ds+ \sum_{i=2}^{n-1}\frac{\overline{Y}(t_i)e^{-t_i}}{1-e^{-2(t_i-t_{i-1})}}\int_{t_{i-1}}^{t_i}  e^{s}\varphi_{x,h}(s)\;ds\\
&~~~~- \sum_{i=2}^{n-1}\frac{\overline{Y}(t_i)e^{t_i}}{1-e^{-2(t_i-t_{i-1})}}\int_{t_{i-1}}^{t_i}  e^{-s}\varphi_{x,h}(s)\;ds+ \sum_{i=2}^{n-1}\overline{Y}(t_i)e^{-t_i}\int_{0}^{t_{i+1}}  e^{s}\varphi_{x,h}(s)\;ds\\
&~~~~+ \sum_{i=2}^{n-1}\frac{\overline{Y}(t_i)e^{t_i}}{1-e^{-2(t_{i+1}-t_{i})}}\int_{t_{i}}^{t_{i+1}}  e^{-s}\varphi_{x,h}(s)\;ds- \sum_{i=2}^{n-1}\frac{\overline{Y}(t_i)e^{-t_i}}{1-e^{-2(t_{i+1}-t_{i})}}\int_{t_{i}}^{t_{i+1}}  e^{s}\varphi_{x,h}(s)\;ds\\
&~~~~ -\frac{f_{x,h}(t_{2})\overline{Y}(t_1)}{1-e^{-2(t_{2}-t_{1})}} e^{-(t_{2}-t_{1})}-\frac{f_{x,h}(t_{n-1})\overline{Y}(t_{n})}{1-e^{-2(t_{n}-t_{n-1})}} e^{-(t_{n}-t_{n-1})}\\
&~~~~-\sum_{i=2}^{n-1} \overline{Y}(t_i)e^{-t_i}\int_{0}^{t_i}  e^{s}\varphi_{x,h}(s)\;ds -\sum_{i=2}^{n-1} \overline{Y}(t_i)e^{t_i}\int_{t_{i}}^{1}  e^{-s}\varphi_{x,h}(s)\;ds.
\end{align*}
Replacing this expression of A in \eqref{gnpro} gives,
\begin{align}
&\hat{g}_n^{pro}(x) = \sum_{i=2}^{n-1} \overline{Y}(t_i)\int_{t_{i-1}}^{t_{i+1}}  e^{|t_i-s|}\varphi_{x,h}(s)\;ds+\sum_{i=2}^{n-2}\frac{\overline{Y}(t_{i+1})e^{-t_{i+1}}-\overline{Y}(t_{i})e^{-t_{i}}}{1-e^{-2(t_{i+1}-t_{i})}}\int_{t_{i}}^{t_{i+1}}  e^{s}\varphi_{x,h}(s)\;ds\nonumber\\
&~~~-\sum_{i=2}^{n-2}\frac{\overline{Y}(t_{i+1})e^{t_{i+1}}-\overline{Y}(t_{i})e^{t_{i}}}{1-e^{-2(t_{i+1}-t_{i})}}\int_{t_{i}}^{t_{i+1}}  e^{-s}\varphi_{x,h}(s)\;ds+\frac{\overline{Y}(t_2)e^{-t_2}}{1-e^{-2(t_2-t_{1})}}\int_{t_{1}}^{t_2}  e^{s}\varphi_{x,h}(s)\;ds\nonumber\\
&~~~-\frac{\overline{Y}(t_{n-1})e^{-t_{n-1}}}{1-e^{-2(t_n-t_{n-1})}}\int_{t_{n-1}}^{t_n}  e^{s}\varphi_{x,h}(s)\;ds-\frac{\overline{Y}(t_2)e^{t_2}}{1-e^{-2(t_2-t_{1})}}\int_{t_{1}}^{t_2}  e^{-s}\varphi_{x,h}(s)\;ds\nonumber\\
&~~~~+\frac{\overline{Y}(t_{n-1})e^{t_{n-1}}}{1-e^{-2(t_n-t_{n-1})}}\int_{t_{n-1}}^{t_n}  e^{-s}\varphi_{x,h}(s)\;ds+\frac{f_{x,h}(t_1)\overline{Y}(t_1)}{1-e^{-2(t_2-t_1)}} -\frac{f_{x,h}(t_{2})\overline{Y}(t_1)}{1-e^{-2(t_{2}-t_{1})}} e^{-(t_{2}-t_{1})}\nonumber\\
&~~~+\frac{f_{x,h}(t_n)\overline{Y}(t_n)}{1-e^{-2(t_n-t_{n-1})}}-\frac{f_{x,h}(t_{n-1})\overline{Y}(t_{n})}{1-e^{-2(t_{n}-t_{n-1})}} e^{-(t_{n}-t_{n-1})}.\label{gnproOU}
\end{align}
Note that Equation \eqref{dif2ornstein} yields,
\begin{align}
&\frac{\overline{Y}(t_1)}{1-e^{-2(t_2-t_1)}}\big( f_{x,h}(t_1) -f_{x,h}(t_{2})e^{-(t_{2}-t_{1})}\big)= \frac{\overline{Y}(t_1)e^{t_1}}{1-e^{-2(t_2-t_1)}}\int_{t_{1}}^{t_{2}}  e^{-s}\varphi_{x,h}(s)\;ds\nonumber \\
& ~~~~~~-\frac{\overline{Y}(t_1)e^{-t_1}}{1-e^{-2(t_2-t_1)}}\int_{t_1}^{t_{2}}  e^{s}\varphi_{x,h}(s)\;ds + \overline{Y}(t_1)e^{-t_1}\int_{t_{1}}^{t_{2}}  e^{s}\varphi_{x,h}(s)\;ds. \label{born1}
\end{align}
Similarly, Equation \eqref{dif1ornstein} yields,
\begin{align}
&\frac{\overline{Y}(t_n)}{1-e^{-2(t_n-t_{n-1})}}\big(f_{x,h}(t_n)-f_{x,h}(t_{n-1})e^{-(t_{n}-t_{n-1})}\big)=\frac{\overline{Y}(t_n)e^{-t_n}}{1-e^{-2(t_n-t_{n-1})}}\int_{t_{n-1}}^{t_n}  e^{s}\varphi_{x,h}(s)\;ds\nonumber \\
&~~ - \frac{\overline{Y}(t_n)e^{t_n}}{1-e^{-2(t_n-t_{n-1})}}\int_{t_{n-1}}^{t_n}  e^{-s}\varphi_{x,h}(s)\;ds+\overline{Y}(t_n)e^{t_n}\int_{t_{n-1}}^{1}  e^{-s}\varphi_{x,h}(s)\;ds.\label{bornn}
\end{align}
We obtain using \eqref{born1} and \eqref{bornn} in  \eqref{gnproOU},
\begin{align*}
&\hat{g}_n^{pro}(x) = \sum_{i=2}^{n-1} \overline{Y}(t_i)\int_{t_{i-1}}^{t_{i+1}}  e^{|t_i-s|}\varphi_{x,h}(s)\;ds+\sum_{i=2}^{n-2}\frac{\overline{Y}(t_{i+1})e^{-t_{i+1}}-\overline{Y}(t_{i})e^{-t_{i}}}{1-e^{-2(t_{i+1}-t_{i})}}\int_{t_{i}}^{t_{i+1}}  e^{s}\varphi_{x,h}(s)\;ds\\
&~~~-\sum_{i=2}^{n-2}\frac{\overline{Y}(t_{i+1})e^{t_{i+1}}-\overline{Y}(t_{i})e^{t_{i}}}{1-e^{-2(t_{i+1}-t_{i})}}\int_{t_{i}}^{t_{i+1}}  e^{-s}\varphi_{x,h}(s)\;ds+\frac{\overline{Y}(t_2)e^{-t_2}}{1-e^{-2(t_2-t_{1})}}\int_{t_{1}}^{t_2}  e^{s}\varphi_{x,h}(s)\;ds\\
&~~~-\frac{\overline{Y}(t_{n-1})e^{-t_{n-1}}}{1-e^{-2(t_n-t_{n-1})}}\int_{t_{n-1}}^{t_n}  e^{s}\varphi_{x,h}(s)\;ds-\frac{\overline{Y}(t_2)e^{t_2}}{1-e^{-2(t_2-t_{1})}}\int_{t_{1}}^{t_2}  e^{-s}\varphi_{x,h}(s)\;ds\\
&~~~~+\frac{\overline{Y}(t_{n-1})e^{t_{n-1}}}{1-e^{-2(t_n-t_{n-1})}}\int_{t_{n-1}}^{t_n}  e^{-s}\varphi_{x,h}(s)\;ds+\frac{\overline{Y}(t_1)e^{t_1}}{1-e^{-2(t_2-t_1)}}\int_{t_{1}}^{t_{2}}  e^{-s}\varphi_{x,h}(s)\;ds\\
&~~~~-\frac{\overline{Y}(t_1)e^{-t_1}}{1-e^{-2(t_2-t_1)}}\int_{t_1}^{t_{2}}  e^{s}\varphi_{x,h}(s)\;ds + \overline{Y}(t_1)e^{-t_1}\int_{t_{1}}^{t_{2}}  e^{s}\varphi_{x,h}(s)\;ds\\
&~~~~+\frac{\overline{Y}(t_n)e^{-t_n}}{1-e^{-2(t_n-t_{n-1})}}\int_{t_{n-1}}^{t_n}  e^{s}\varphi_{x,h}(s)\;ds - \frac{\overline{Y}(t_n)e^{t_n}}{1-e^{-2(t_n-t_{n-1})}}\int_{t_{n-1}}^{t_n}  e^{-s}\varphi_{x,h}(s)\;ds\\
&~~~~+\overline{Y}(t_n)e^{t_n}\int_{t_{n-1}}^{1}  e^{-s}\varphi_{x,h}(s)\;ds\end{align*}\begin{align*}
& = \sum_{i=2}^{n-1}\overline{Y}(t_i) \int_{t_{i-1}}^{t_{i+1}} e^{|s-t_i|}\varphi_{x,h}(s)\;ds +\overline{Y}(t_1)\int_{0}^{t_2} e^{s-t_1}\varphi_{x,h}(s)\;ds+\overline{Y}(t_{n})\int_{t_{n-1}}^{1} e^{t_n-s}\varphi_{x,h}(s)\;ds\\
& ~~~~  -\sum_{i=1}^{n-1} \frac{e^{t_{i+1}}\overline{Y}(t_{i+1})-e^{t_i}\overline{Y}(t_i)}{1-e^{-2(t_{i+1}-t_i)}} \int_{t_i}^{t_{i+1}} e^{-s}\varphi_{x,h}(s)\;ds\\
& ~~~~+\sum_{i=1}^{n-1} \frac{e^{-t_{i+1}}\overline{Y}(t_{i+1})-e^{-t_i}\overline{Y}(t_i)}{1-e^{-2(t_{i+1}-t_i)}} \int_{t_i}^{t_{i+1}} e^{s}\varphi_{x,h}(s)\;ds .
\end{align*}
This concludes the proof of Proposition \ref{ornsteinestimator}. $\Box$
\subsection{Proof of Lemma \ref{alphapos}.}
Let $(u,v) \in [-1,1]^2$. We first consider the triangle $\{ -1<u<v<1\}$ which is further split into smaller triangles:
\[D_1 = \{ 0<u<v<1\},~ D_2 = \{ -1<u<0<v<1\}~~\text{and}~~D_3 = \{-1<u<v<0\}.\]
Let $b \in ]0,1[$. For $(u,v) \in D_1$, using Assumption $(A)$, Taylor expansion of $R$ around $(x,x)$ gives,
\begin{align*}
 R(x+b u,x+b v)&  = R(x,x+b v) +buR^{(1,0)}(x,x+bv) +\frac{1}{2}b^2u^2R^{(2,0)}(\varepsilon_x,x+bv)  \\
& = R(x,x)+ bvR^{(0,1)}(x,\eta_x) +buR^{(1,0)}(x,x+bv)+\frac{1}{2}b^2u^2R^{(2,0)}(\varepsilon_x,x+bv),
\end{align*}
where $x<\varepsilon_x <x+bu < x+bv$  and $ x< \eta_x <  x+bv$. Thus,
\[R(x+b u,x+b v) = R(x,x)+ bvR^{(0,1)}(x,x^+) +buR^{(0,1)}(x,x^-) + o(b). \]
Now, for $(u,v) \in D_2$  we obtain in the same way,
\begin{align*}
& R(x+b u,x+b v)  = R(x,x+b v) +buR^{(1,0)}(x,x+bv)+\frac{1}{2}b^2u^2R^{(2,0)}(\varepsilon_x,x+bv)   \\
& = R(x,x)+ bvR^{(0,1)}(x,\eta_x) +buR^{(1,0)}(x,x+bv)+\frac{1}{2}b^2u^2R^{(2,0)}(\varepsilon_x,x+bv) ,
\end{align*}
where $x +bu<\varepsilon_x <x< x+bv$  and $ x< \eta_x <  x+bv$. Thus,
\[R(x+b u,x+b v) = R(x,x)+ bvR^{(0,1)}(x,x^+) +buR^{(0,1)}(x,x^-) + o(b). \]
Finally, for $(u,v) \in D_3$  we get,
\begin{align*}
& R(x+b u,x+b v)  = R(x+bu,x) +bvR^{(0,1)}(x+bu,x) +\frac{1}{2}b^2v^2R^{(0,2)}(x+bu,\eta_x)  \\
& = R(x,x)+ubR^{(1,0)}(\varepsilon_x,x)+ bvR^{(0,1)}(x+bu,\eta_x) +\frac{1}{2}b^2v^2R^{(0,2)}(x+bu,\eta_x) ,
\end{align*}
where $x+hu < x+bv< \eta_x <  x$  and $x+bu<\varepsilon_x < x$. Thus,
\[R(x+b u,x+b v) = R(x,x)+ bvR^{(0,1)}(x,x^+) +buR^{(0,1)}(x,x^-) + o(b). \]
\noindent Hence for $v>u$ we have,
\begin{align*}
R(x+b u,x+b v) & =  R(x,x) +\frac{b}{2} \big(R^{(0,1)}(x,x^+) + R^{(0,1)}(x,x^-) \big) (u+v) \\
&~~~+ \frac{b}{2} \big(R^{(0,1)}(x,x^+) - R^{(0,1)}(x,x^-) \big) (v-u) +o(b).
\end{align*}
Similarly, we obtain for the triangular $\{1 > u > v > -1\}$,
\begin{align*}
R(x+b u,x+b v) & = R(x,x) +\frac{b}{2} \big(R^{(0,1)}(x,x^+) + R^{(0,1)}(x,x^-) \big) (u+v) \\
&~~~+ \frac{b}{2} \big(R^{(0,1)}(x,x^+) - R^{(0,1)}(x,x^-) \big) (u-v).
\end{align*}
Thus, for $(u,v) \in [-1,1]^2$  we have,
\begin{align}
R(x+b u,x+b v) & = R(x,x) +\frac{b}{2} \big(R^{(0,1)}(x,x^+) + R^{(0,1)}(x,x^-) \big) (u+v) \nonumber\\
&~~~+ \frac{b}{2} \big(R^{(0,1)}(x,x^+) - R^{(0,1)}(x,x^-) \big) |u-v|.
\label{formuleRuh}
\end{align}
Consider now a function $g$, bounded and integrable on $[-1,1]$. The  Dominated Convergence Theorem yields that $R(.,t)\times g$ is an integrable function for every $t \in [-1,1]$.
Using \eqref{formuleRuh} and putting,
\[\gamma(x)=\frac{1}{2}\big(R^{(0,1)}(x,x^+) + R^{(0,1)}(x,x^-) \big),\] we obtain,
\begin{equation}\label{dominated convergence}
\begin{split}
&\iint_{[-1,1]^2}  R(x+bu,x+bv)g(u)g(v) dudv = R(x,x) \left( \int_{-1}^{1} g(u)du \right)^2 \\
&+ 2\gamma (x)b \int_{-1}^1 g(u)du \int_{-1}^1 vg(v)dv -\frac{b}{2}\alpha (x)\iint_{[-1,1]^2} g(u)g(v) |u-v| dudv + o(b) .
\end{split}
\end{equation}
The left side of \eqref{dominated convergence} is non-negative since the autocovariance function $R$ is a non-negative definite function. Taking $g(u) = u  \mathrm{1}_{[-1,1]}(u)$ we obtain,
\[ \int_{-1}^1 g(u)du = 0 ~~\text{and} \iint_{[-1,1]^2} uv |u-v| dudv = - \frac{8}{15 }.\]
Thus,
\[\frac{4}{15 } \alpha(x) + o(b) \geq 0. \]
Taking $b$ small enough concludes the proof of Lemma \ref{alphapos}. $\Box$
\subsection{Proof of Lemma \ref{mxhlemma}.}
The great lines of this proof are based on the work of Sacks and Ylvisaker (1966) \cite{Sacks and Ylvisaker} (c.f.  Lemma 3.2 there). Let $x,h \in ]0,1[$ and put $g_n=P_{T_n}f_{x,h}$, it is shown by \eqref{mvar=|||PT_n|} in the Appendix  that,\[g_n(t_i) = \sum_{j=1}^{n} m_{x,h}(t_j)R(t_j,t_i)~~~\text{for all}~i=1,\cdots,n.\]
On the one hand, Assumption $(A)$ yields that $g_n$ is twice differentiable on $[0,1]$ except on $T_n$, but it has left and right derivatives.  Thus, for every $i = 1,\dots,n$ we have,
\[g_n'(t_i^-)= \sum_{j=1}^{n}m_{x,h}(t_j)R^{(0,1)}(t_j,t_i^-)~~\text{and}~~g_n'(t_i^+)= \sum_{j=1}^{n}m_{x,h}(t_j)R^{(0,1)}(t_j,t_i^+). \]
Since for $j\ne i$, $R^{(0,1)}(t_j,t_i^-)=R^{(0,1)}(t_j,t_i^+)$ then Assumption $(B)$ yields,
\begin{equation}
g_n'(t_i^-)-g_n'(t_i^+) = \alpha(t_i)m_{x,h}(t_i).
\label{malpha}
\end{equation}
On the other hand, Assumption $(A)$ yields that $f_{x,h}$ (as defined by \eqref{f_xh}) is twice differentiable on $]0,1[$, thus
for $i=1,\dots,n-1$, Taylor expansion of $f_{x,h}-g_n$ around $t_i$ gives,
\begin{align*}
f_{x,h}(t_{i+1})-g_n(t_{i+1})& = (f_{x,h}(t_i)-g_n(t_i))+d_i(f'_{x,h}(t_i)-g_n'(t_i^+))+\frac{1}{2}d_i^2(f_{x,h}''(\sigma_i)-g_n''(\sigma_i)),
\end{align*}
where $d_i=t_{i+1}-t_i $ and $\sigma_i \in ]t_i,t_{i+1}[$. Recall that, for all $i= 1,\dots,n$,  $f_{x,h}(t_i)=g_n(t_i)$ (see the Appendix,  Equation \eqref{appendixPtn=fxh}). Thus,
\begin{equation}
f'_{x,h}(t_i)-g_n'(t_i^+)=-\frac{1}{2}d_i(f_{x,h}''(\sigma_i)-g_n''(\sigma_i)),
\label{2diff'}
\end{equation}
 Similarly, for $i=2,\dots,n$, we have,
\begin{equation}
f'_{x,h}(t_i)-g_n'(t_i^-)=\frac{1}{2}d_{i-1}(f_{x,h}''(\theta_i)-g_n''(\theta_i)),
\label{2diff''}
\end{equation}
for some $\theta_i \in ]t_{i-1},t_i[$.   We obtain subtracting \eqref{2diff''} from \eqref{2diff'} and using \eqref{malpha} for $i=2,\dots,n-1$,
\begin{align}\label{alphamxh}
\alpha(t_i)m_{x,h}(t_i)& = -\frac{1}{2}d_i (f_{x,h}''(\sigma_i)-g_n''(\sigma_i))-\frac{1}{2}d_{i-1} (f_{x,h}''(\theta_i)-g_n''(\theta_i)).
\end{align}
We shall now control the last expression. On the one hand we have,
\begin{equation}
f_{x,h}'(t) = \int_{0}^{t} R^{(0,1)}(s,t^+)\varphi_{x,h}(s)\;ds+\int_{t}^{1} R^{(0,1)}(s,t^-)\varphi_{x,h}(s)\;ds,~~~~~~
\label{f't}
\end{equation}
and,
\begin{align}
f_{x,h}''(t) & = (R^{(0,1)}(t,t^+)-R^{(0,1)}(t,t^-)) \varphi_{x,h}(t)+ \int_{0}^{1}R^{(0,2)}(s,t^+)\varphi_{x,h}(s)\;ds\nonumber\\
&-\alpha(t) \varphi_{x,h}(t)+ \int_{0}^{1}R^{(0,2)}(s,t^+)\varphi_{x,h}(s)\;ds.
\label{f''t}
\end{align}
On the other hand we know, using (F3) in the Appendix, that  every function in the RKHS($R$), noted by $\mathcal F (\varepsilon)$, is continuous,  hence Assumption $(C)$ implies that  $R^{(0,2)}(\cdot,t^+) $ is a continuous function on $[0,1]$ for every fixed $t \in [0,1]$. Thus, \[R^{(0,2)}(t,t^+)=\lim_{s \downarrow t } R^{(0,2)}(s,t^+)=\lim_{s \downarrow t } R^{(0,2)}(s,t^-)= R^{(0,2)}(t,t^-),\] from which we get that $R^{(0,2)}(t,t)$ exists. Hence for $i = 1,\dots,n$ we have,
\begin{align}
g_n''(t_i^-)=g_n''(t_i^+)=\sum_{j=1}^{n} m_{x,h}(t_j)R^{(0,2)}(t_j,t_i).\label{gsecond}
\end{align}
In addition, it is shown by (F4) in the Appendix  that for every $t \in [0,1]$,
\begin{align}
f_{x,h}'' (t) - g_n''(t) & =  -\alpha(t) \varphi_{x,h}(t) + \langle R^{(0,2)}(\cdot,t),f_{x,h}-g_n\rangle, \label{f''-g''=inner}
\end{align}
where $\langle \cdot,\cdot \rangle$ is the inner product on $\mathcal F (\varepsilon) $.
Injecting \eqref{f''-g''=inner} in \eqref{alphamxh} we obtain,
\begin{align*}
\alpha(t_i)m_{x,h}(t_i)& = \frac{1}{2}d_i\alpha(\sigma_i)\varphi_{x,h}(\sigma_i)+\frac{1}{2}d_{i-1}\alpha(\theta_i)\varphi_{x,h}(\theta_i)-\frac{1}{2}d_i\langle R^{(0,2)}(\cdot,\sigma_i),f_{x,h}-g_n\rangle\\
&~~~~-\frac{1}{2}d_{i-1}\langle R^{(0,2)}(\cdot,\theta_i),f_{x,h}-g_n\rangle.
\end{align*}
Using Assumption $(B)$ we obtain for $i=2,\dots,n-1$,
\begin{align}
m_{x,h}(t_i) &=\frac{1}{2}(d_i+d_{i-1})\varphi_{x,h}(t_i)~ + \frac{1}{2\alpha(t_i)}d_i\big(\alpha(\sigma_i)\varphi_{x,h}(\sigma_i)-\alpha(t_i)\varphi_{x,h}(t_i)\big)\nonumber \\
& + \frac{1}{2\alpha(t_i)}d_{i-1}\big(\alpha(\theta_i)\varphi_{x,h}(\theta_i)-\alpha(t_i)\varphi_{x,h}(t_i)\big)-\frac{1}{2\alpha(t_i)}d_i\langle R^{(0,2)}(\cdot,\sigma_i),f_{x,h}-g_n\rangle \nonumber\\
& -\frac{1}{2\alpha(t_i)}d_{i-1}\langle R^{(0,2)}(\cdot,\theta_i),f_{x,h}-g_n\rangle\nonumber\\
&\overset{\Delta}{=} \frac{1}{2}(d_i+d_{i-1})\varphi_{x,h}(t_i)~ + A_{i}^{(1)}+ A_{i}^{(2)}- A_{i}^{(3)}-A_{i}^{(4)}, \label{mti}
\end{align}
Using the Cauchy-Schwartz inequality, Assumption $(C)$ and Equation \eqref{lamajoration} (in the proof of Proposition \ref{propoflimit} below) we obtain,
\begin{equation}
 |A_{i}^{(3)}+A_{i}^{(4)}| \leq ~ \underset{0\leq t \leq 1 }{\sup}~\frac{1}{2\alpha(t)}||R^{(0,2)}(.,t)||\frac{\sqrt C}{ \sqrt{h}}~ \underset{0 \leq j \leq n}{\sup} d_j^2 \overset{\Delta}{=} \beta_{n,h}, \label{N}
\end{equation}
 where $C$ is a positive constant defined in Proposition \ref{propoflimit} below.\\
 Recall that $\varphi_{x,h}$ is of support $[x-h,x+h]$, thus for $t_i $  such that $[t_{i-1},t_{i+1}]\cap ]x-h,x+h[ = \emptyset$, $\varphi_{x,h}(t)=0$ so that $A_{i}^{(1)}=0$ and $A_{i}^{(2)}=0$.
For $t_i$ such that $[t_{i-1},t_{i+1}]\cap ]x-h,x+h[ \ne \emptyset$, let,
\begin{equation}
\alpha_{n,h} =  \underset{0 \leq i \leq n}{\sup}~~ \underset{t_i \leq s,t\leq t_{i+1}}{\sup} \frac{1}{2\alpha(t)}d_i|\alpha(s)\varphi_{x,h}(s)-\alpha(t)\varphi_{x,h}(t)|.
\label{M}
\end{equation}
\noindent We obtain using \eqref{N} and \eqref{M} together with \eqref{mti} for $i = 2,\dots,n-1$,
\[m_{x,h}(t_i) =
\begin{cases}
 \frac{1}{2} \varphi_{x,h} (t_i) (t_{i+1}-t_{i-1})~+~O\big(\alpha_{n,h}+\beta_{n,h}\big)~~~if~ [t_{i-1},t_{i+1}]\cap ]x-h,x+h[ \ne \emptyset \\[0.1cm]
O\big(\beta_{n,h} \big)~~~~~~~~~~~~~~~~~~~~~~~~~~~~~~~~~~~~~~~~~~~~otherwise.
\end{cases}\]
After having obtained $m_{x,h}(t_i)$ for $i = 2,\dots,n-1$, we are now  able to obtain $m_{x,h}(t_1)$ and $m_{x,h}(t_n)$. We have for $i = 1,\dots,n$,
\begin{equation}
R(t_1,t_i)m_{x,h}(t_1)+R(t_n,t_i)m_{x,h}(t_n)=f_{x,h}(t_i)-\sum_{j=2}^{n-1} m_{x,h}(t_j)R(t_j,t_i).
\label{m1mn}
\end{equation}
Recall that $N_{T_n}=\Card~ I_{x,h} = \Card~\{i=1,\cdots,n~:~[t_{i-1},t_{i+1}]\cap ]x-h,x+h[ \ne \emptyset \}$  and that $t_{x,i}$ are the points of $T_n$ for which $i \in I_{x,h}$. We have,
\begin{align*}
\sum_{j=2}^{n-1} m_{x,h}(t_j)R(t_j,t_i) & = \sum_{j=1}^{N_{T_n}} m_{x,h}(t_{x,j})R(t_{x,j},t_i)+\sum_{j=2}^{n-1}1_{\{ j \notin I_{x,h}\}} m_{x,h}(t_j)R(t_j,t_i).
\end{align*}
On the one hand, we have using \eqref{mti} (where $A_{x,j}$ stands for $A_j$ with $t_j$ replaced by $t_{x,j}$),
{\small
\begin{align}
&\sum_{j=2}^{n-1} m_{x,h}(t_j)R(t_j,t_i)= \frac{1}{2}\sum_{j=1}^{N_{T_n}} (d_{x,j}+d_{x,j-1})\varphi_{x,h}(t_{x,j})R(t_{x,j},t_i)\nonumber\\
&~~~~~~ +\sum_{j=1}^{N_{T_n}} ({A_{x,j}^1} + {A_{x,j}^2}- {A_{x,j}^3}- {A_{x,j}^4})R(t_{x,j},t_i) -\sum_{j=2}^{n-1}1_{\{ j \notin I_{x,h}\}} ({A_{j}^3	}+{A_{j}^4})R(t_j,t_i)\nonumber\\
& =\frac{1}{2}\sum_{j=1}^{N_{T_n}} (d_{x,j}+d_{x,j-1})\varphi_{x,h}(t_{x,j})R(t_{x,j},t_i) +\sum_{j=1}^{N_{T_n}}({A_{x,j}^1} + {A_{x,j}^2})R(t_{x,j},t_i)-\sum_{j=1}^{n} ({A_{j}^3}+{A_{j}^4})R(t_j,t_i). \label{sum2n}
\end{align}}
On the other hand,
\begin{align}
&f_{x,h}(t_i)  = \int_0^1 R(s,t_i) \varphi_{x,h}(s)\;ds = \int_{x-h}^{x+h} R(s,t_i) \varphi_{x,h}(s)\;ds \nonumber=  \frac{1}{2} \sum_{j=1}^{N_{T_n}} \int_{t_{x,j-1}}^{t_{x,j+1}} R(s,t_i)\varphi_{x,h}(s)\;ds\nonumber\\
& = \frac{1}{2} \sum_{j=1}^{N_{T_n}} (d_{x,j}+d_{x,j-1})R(t_{x,j},t_i)\varphi_{x,h}(t_j)+\frac{1}{2} \sum_{j=1}^{N_{T_n}}\int_{t_{x,j-1}}^{t_{x,j+1}} (R(s,t_i)\varphi_{x,h}(s)- R(t_{x,j},t_i)\varphi_{x,h}(t_{x,j}))\;ds.
\label{fxhsomme}
\end{align}
Inserting \eqref{sum2n} and \eqref{fxhsomme} in \eqref{m1mn} we obtain for $i = 1,\dots,n$,
\begin{align*}
&R(t_1,t_i)m_{x,h}(t_1)+R(t_n,t_i)m_{x,h}(t_n) =\frac{1}{2} \sum_{j=1}^{N_{T_n}}\int_{t_{x,j-1}}^{t_{x,j+1}} (R(s,t_i)\varphi_{x,h}(s)- R(t_{x,j},t_i)\varphi_{x,h}(t_{x,j}))\;ds \\
& ~~~-\sum_{j=1}^{N_{T_n}}({A_{x,j}^1} + {A_{x,j}^2})R(t_{x,j},t_i)+\sum_{j=1}^{n} ({A_{j}^3}+{A_{j}^4})R(t_j,t_i) \overset{\Delta}= \Phi_{x,h}(t_i) .
\end{align*}
We then obtain the following linear system,
\begin{align}
\begin{cases}
& R(t_1,t_1)m_{x,h}(t_1) + R(t_n,t_1)m_{x,h}(t_1)  = \Phi_{x,h}(t_1).\\
& \label{system}\\
& R(t_1,t_n)m_{x,h}(t_1) + R(t_n,t_n)m_{x,h}(t_n)  = \Phi_{x,h}(t_n).
\end{cases}
\end{align}
Solving \eqref{system} for $m_{x,h}(t_1)$ and $m_{x,h}(t_n)$ we obtain,
\begin{align}
m_{x,h}(t_1) & = \frac{R(t_n,t_n)\Phi_{x,h}(t_1)-R(t_1,t_n)\Phi_{x,h}(t_n)}{R(t_1,t_1)R(t_n,t_n)-R(t_1,t_n)^2}. \\
m_{x,h}(t_n) & = \frac{R(t_1,t_1)\Phi_{x,h}(t_n)-R(t_1,t_n)\Phi_{x,h}(t_1)}{R(t_1,t_1)R(t_n,t_n)-R(t_1,t_n)^2} .
\end{align}
Finally, simple calculations yield,
\begin{equation*}
m_{x,h}(t_1)=O(N_{T_n}\alpha_{n,h}+n\beta_{n,h}\big) ~~~\text{and}~~~~m_{x,h}(t_n)=O(N_{T_n}\alpha_{n,h}+n\beta_{n,h}\big).
\end{equation*}
This completes the proof of Lemma \ref{mxhlemma}. $\Box$
\subsection{Proof of Proposition \ref{biastheorem}.}
Recall that $N_{T_n}=\Card~ I_{x,h} = \Card~\{i=1,\cdots,n~:~[t_{i-1},t_{i+1}]\cap ]x-h,x+h[ \ne \emptyset \}$ and denote by $t_{x,i}$ the points of $T_n$ for which $i \in I_{x,h}$, that is $T_n \cap ]x-h,x+h[ = \{t_{x,2},\cdots,t_{x,N_{T_n}-1}\}$. Since $\mathbb E (\overline{Y} (t_i)) = g(t_i)$ then,
\begin{align*}
	\mathbb{E}(\hat{g}^{pro}_n (x) ) &= \sum_{j=1}^n m_{x,h} (t_j) g (t_j)\\
	& =\sum_{i=1}^{N_{T_n}} m_{x,h} (t_{x,i}) g (t_{x,i})+ \sum_{j=2}^{n-1} 1_{\{i \notin I_{x,h}\}} m_{x,h} (t_j) g (t_j)+ m_{x,h} (t_1) g (t_1)+ m_{x,h} (t_n) g (t_n).
\end{align*}
Using the asymptotic approximation of ${m_{x,h}}_{|T_n}$ given in Lemma \ref{mxhlemma} we obtain,
\begin{align}
E(\hat{g}^{pro}_n (x) )& = \frac{1}{2}\sum_{i=1}^{N_{T_n}} (t_{x,i+1}-t_{x,i-1}) \varphi_{x,h}(t_{x,i}) g (t_{x,i})+ O\big(N_{T_n}\alpha_{n,h}+n \beta_{n,h}\big),\label{Egdansproof}
\end{align}
For $x \in [0,1]$ let, \[I_h(x) = \int_{x-h}^{x+h} \varphi_{x,h}(t)g(t)~dt=\frac{1}{2}\sum_{i=1}^{N_{T_n}}\int_{t_{x,i-1}}^{t_{x,i+1}} \varphi_{x,h}(t)g(t)~dt,\]
and write,
\begin{align}
\mathbb{E}(\hat{g}^{pro}_n (x) ) & =  \mathbb{E}(\hat{g}^{pro}_n (x) )- I_h(x)+I_h(x) = \Delta_{x,h} +I_h(x) +O\big(N_{T_n}\alpha_{n,h}+n \beta_{n,h}\big), \label{thissana}
\end{align}
where, \[\Delta_{x,h} =\frac{1}{2}\sum_{i=1}^{N_{T_n}}\int_{t_{x,i-1}}^{t_{x,i+1}} \Big(  \varphi_{x,h}(t_{x,i}) g (t_{x,i})- \varphi_{x,h}(t)g(t) \Big)~dt.\]
We first control $\Delta_{x,h}$.
We have,
\begin{align*}
\Delta_{x,h}
&= \frac{1}{2}\sum_{i=1}^{N_{T_n}}\varphi_{x,h}(t_{x,i})\int_{t_{x,i-1}}^{t_{x,i+1}}( g (t_{x,i})- g(t))~dt+\frac{1}{2}\sum_{i=1}^{N_{T_n}}\int_{t_{x,i-1}}^{t_{x,i+1}}g(t)(\varphi_{x,h}(t_{x,i})- \varphi_{x,h}(t))~dt.
\end{align*}
Since $\varphi_{x,h}$ is in $C^1$ and $g$ is  in $C^2$ then Taylor expansions of $\varphi_{x,h}$ and $g$ give,
\begin{align*}
& g(t)  = g(t_{x,i}) + (t-t_{x,i})g'(t_{x,i})+\frac{1}{2} (t-t_{x,i})^2g''(\theta_{x,i}),
\end{align*}
and,
\begin{align*}
& \varphi_{x,h}(t) = \varphi_{x,h}(t_{x,i}) + (t-t_{x,i}) \varphi_{x,h}'(\eta_{x,i}),
\end{align*}
for some $\theta_{x,i}$ and $\eta_{x,i}$ between $t$ and $t_{x,i}$.  Thus,
\begin{align*}
&\Delta_{x,h} =   - \frac{1}{2}\sum_{i=1}^{N_{T_n}}\varphi_{x,h}(t_{x,i})g'(t_{x,i})\int_{t_{x,i-1}}^{t_{x,i+1}}(t-t_{x,i}) ~dt- \frac{1}{4}\sum_{i=1}^{N_{T_n}}\varphi_{x,h}(t_{x,i})\int_{t_{x,i-1}}^{t_{x,i+1}} g''(\theta_{x,i})(t-t_{x,i})^2 ~dt\\
& -\frac{1}{2}\sum_{i=1}^{N_{T_n}}g(t_{x,i})\int_{t_{x,i-1}}^{t_{x,i+1}}\varphi_{x,h}'(\eta_{x,i}) (t-t_{x,i}) ~dt-\frac{1}{2}\sum_{i=1}^{N_{T_n}}g'(t_{x,i})\int_{t_{x,i-1}}^{t_{x,i+1}}\varphi_{x,h}'(\eta_{x,i}) (t-t_{x,i})^2 ~dt\\
&-\frac{1}{4}\sum_{i=1}^{N_{T_n}}\int_{t_{x,i-1}}^{t_{x,i+1}} g''(\theta_{x,i})\varphi_{x,h}'(\eta_{x,i}) (t-t_{x,i})^3 ~dt.
\end{align*}
Recall that $g'$ and $g''$ are both bounded and that,
\begin{equation}
\underset{0 \leq t \leq 1}{\sup}~|\varphi_{x,h}(t)| < \frac{c}{h}~~\text{and}~~\underset{0 \leq t \leq 1}{\sup}~|\varphi_{x,h}'(t)| < \frac{c'}{h^2}, \label{phiboundries}
\end{equation}
for appropriate positive constants $c$ and $c'$. Using this we obtain,
\begin{align*}
&\frac{1}{4}\sum_{i=1}^{N_{T_n}}\varphi_{x,h}(t_{x,i})\int_{t_{x,i-1}}^{t_{x,i+1}} g''(\theta_{x,i})(t-t_{x,i})^2 ~dt = O\bigg(\frac{N_{T_n}}{h} \underset{0\leq j \leq 1}{\sup}\; d_{j,n}^3\bigg)\\
&\frac{1}{2}\sum_{i=1}^{N_{T_n}}g'(t_{x,i})\int_{t_{x,i-1}}^{t_{x,i+1}} \varphi_{x,h}'(\eta_{x,i})(t-t_{x,i})^2 ~dt~= O\bigg(\frac{N_{T_n}}{h^2} \underset{0\leq j \leq 1}{\sup}\; d_{j,n}^3\bigg)\\
&\frac{1}{4}\sum_{i=1}^{N_{T_n}}\int_{t_{x,i-1}}^{t_{x,i+1}}g''(\theta_{x,i})\varphi_{x,h}'(\eta_{x,i}) (t-t_{x,i})^2 ~dt = O\bigg(\frac{N_{T_n}}{h^2}  \underset{0\leq j \leq 1}{\sup}\; d_{j,n}^3\bigg).
\end{align*}
Thus,
\begin{align*}
\Delta_{x,h}
& =  - \frac{1}{2}\sum_{i=1}^{N_{T_n}}\varphi_{x,h}(t_{x,i})g'(t_{x,i})\int_{t_{x,i-1}}^{t_{x,i+1}}(t-t_{x,i}) dt-\frac{1}{2}\sum_{i=1}^{N_{T_n}}g(t_{x,i})\varphi_{x,h}'(t_{x,i})\int_{t_{x,i-1}}^{t_{x,i+1}} (t-t_{x,i}) ~dt\\
& ~~-\frac{1}{2}\sum_{i=1}^{N_{T_n}}g(t_{x,i})\int_{t_{x,i-1}}^{t_{x,i+1}} (t-t_{x,i})\Big(\varphi_{x,h}'(\eta_{x,i})-\varphi_{x,h}'(t_{x,i}) \Big) ~dt+ O\bigg(\frac{N_{T_n}}{h^2}  \underset{0\leq j \leq 1}{\sup}\; d_{j,n}^3\bigg).
\end{align*}
Since $\varphi_{x,h}'$ is Lipschitz then,
\begin{align*}
\sum_{i=1}^{N_{T_n}}g(t_{x,i})\int_{t_{x,i-1}}^{t_{x,i+1}} (t-t_{x,i})\Big(\varphi_{x,h}'(\eta_{x,i})-\varphi_{x,h}'(t_{x,i}) \Big) ~dt = O\bigg(\frac{N_{T_n}}{h^3}  \underset{0\leq j \leq 1}{\sup}\; d_{j,n}^3\bigg).
\end{align*}
Thus,
\begin{align*}
\Delta_{x,h}
& =  - \frac{1}{2}\sum_{i=1}^{N_{T_n}}\varphi_{x,h}(t_{x,i})g'(t_{x,i})\int_{t_{x,i-1}}^{t_{x,i+1}}(t-t_{x,i}) dt-\frac{1}{2}\sum_{i=1}^{N_{T_n}}g(t_{x,i})\varphi_{x,h}'(t_{x,i})\int_{t_{x,i-1}}^{t_{x,i+1}} (t-t_{x,i}) ~dt\\
&~~+O\bigg(\frac{N_{T_n}}{h^3}  \underset{0\leq j \leq 1}{\sup}\; d_{j,n}^3\bigg).
\end{align*}
Basic integration gives,
\begin{align*}
\Delta_{x,h} & =  - \frac{1}{4}\sum_{i=1}^{N_{T_n}}\varphi_{x,h}(t_{x,i})g'(t_{x,i})(d_{x,i}^2-d_{x,i-1}^2) -\frac{1}{4}\sum_{i=1}^{N_{T_n}}g(t_{x,i})\varphi_{x,h}'(t_{x,i})(d_{x,i}^2-d_{x,i-1}^2)\\
&~~~~ + O\bigg(\frac{N_{T_n}}{h^3}  \underset{0\leq j \leq 1}{\sup}\; d_{j,n}^3\bigg).
\end{align*}
We shall show that,
\begin{align*}
A \overset{\Delta}{=} \sum_{i=1}^{N_{T_n}}\varphi_{x,h}(t_{x,i})g'(t_{x,i})(d_{x,i}^2-d_{x,i-1}^2) & =  O\bigg(\frac{N_{T_n}}{h^2}  \underset{0\leq j \leq 1}{\sup}\; d_{j,n}^3\bigg), \\
B \overset{\Delta}{=} \sum_{i=1}^{N_{T_n}}g(t_{x,i})\varphi_{x,h}'(t_{x,i})(d_{x,i}^2-d_{x,i-1}^2) & =  O\bigg(\frac{N_{T_n}}{h^3}  \underset{0\leq j \leq 1}{\sup}\; d_{j,n}^3\bigg).
\end{align*}
Starting with the term $A$. Recall that, since $\varphi$ is of support $[x-h,x+h]$ and $t_{x,1}, t_{x,N_{T_n}-1} \notin ]x-h,x+h[$, then $\varphi_{x,h}(t_{x,N_{T_n}})=\varphi_{x,h}(t_{x,1})=0$ thus,
\begin{align*}
A & = \sum_{i=2}^{N_{T_n}-1}\varphi_{x,h}(t_{x,i})g'(t_{x,i})d_{x,i}^2-\sum_{i=1}^{N_{T_n}-2}\varphi_{x,h}(t_{x,i+1})g'(t_{x,i+1})d_{x,i}^2\\
& = \sum_{i=2}^{N_{T_n}-2}\big( \varphi_{x,h}(t_{x,i})g'(t_{x,i})-\varphi_{x,h}(t_{x,i+1})g'(t_{x,i+1})\big)d_{x,i}^2 +\bigg(\varphi_{x,h}(t_{x,N_{T_n}-1})g'(t_{x,N_{T_n}-1})d_{x,N_{T_n}-1}^2\\
&~~~~~~~~-\varphi_{x,h}(t_{x,2})g'(t_{x,2})d_{x,1}^2\bigg)\\
& \overset{\Delta}{=} A_1+A_2.
\end{align*}
On the one hand, Taylor expansions of $\varphi_{x,h}$ around $t_{x,N_{T_n}}$ and $t_{x,1}$ yield,
\begin{align*}
& \varphi_{x,h}(t_{x,N_{T_n}-1})= (t_{x,N_{T_n}-1}-t_{x,N_{T_n}})\varphi_{x,h}'(\gamma_{x,N_{T_n}}),\\
& \varphi_{x,h}(t_{x,2}) =(t_{x,2}-t_{x,1})\varphi_{x,h}'(\gamma_{x,1}),
\end{align*}
for some $\gamma_{x,N_{T_n}} \in ]t_{x,N_{T_n}-1},t_{x,N_{T_n}}[$ and some $\gamma_{x,1} \in ]t_{x,1},t_{x,2}[$.
Using \eqref{phiboundries} and the fact that $g'$ is bounded we obtain, \[A_2=O\bigg(\frac{1}{h^2}  \underset{0\leq j \leq 1}{\sup}\; d_{j,n}^3\bigg).\]
On the other hand we have,
\begin{align*}
A_1 &= \sum_{i=2}^{N_{T_n}-2}\big( \varphi_{x,h}(t_{x,i})g'(t_{x,i})-\varphi_{x,h}(t_{x,i+1})g'(t_{x,i+1})\big)d_{x,i}^2 \\
& =  \sum_{i=2}^{N_{T_n}-2} \varphi_{x,h}(t_{x,i})\big( g'(t_{x,i})-g'(t_{x,i+1})\big)d_{x,i}^2 + \sum_{i=2}^{N_{T_n}-2}g'(t_{x,i+1}) \big( \varphi_{x,h}(t_{x,i})-\varphi_{x,h}(t_{x,i+1})\big)d_{x,i}^2.
\end{align*}
Since $\varphi_{x,h}$ is in $C^1$ and $g$ is in $C^2$ then using \eqref{phiboundries}, we obtain, \[ 	A_1	 = O\Big(\frac{N_{T_n}}{h^2}  \underset{0\leq j \leq 1}{\sup}\; d_{j,n}^3\Big).\]
In a similar way and from Assumption $(D)$, we obtain,  \[ B = O\Big(\frac{N_{T_n}}{h^3}  \underset{0\leq j \leq 1}{\sup}\; d_{j,n}^3\Big). \] Hence,
\begin{align*}
\Delta_{x,h}  & = O\Big(\frac{N_{T_n}}{h^3} \underset{0\leq j \leq 1}{\sup}\; d_{j,n}^3\Big).
\end{align*}
Thus using \eqref{thissana},
\begin{align*}
\mathbb{E}(\hat{g}^{pro}_n (x) ) & =  I_h(x) +O(N_{T_n}\alpha_{n,h}+n \beta_{n,h})+O\Big(\frac{N_{T_n}}{h^3} \underset{0\leq j \leq 1}{\sup}\; d_{j,n}^3\Big).
\end{align*}
The control of $I_h(x)$ is classical and it can bee seen from Gasser and Müller (1984) \cite{Gasser and Müller } that,
\begin{equation}
I_h(x)  = g(x) +  \frac{1}{2} h^2 g''(x)\int_{-1}^{1}t^2 K(t)~dt + o (h^2).\label{chepa}
\end{equation} Finally,
\begin{align*}
\mathbb{E}(\hat{g}^{pro}_n (x) ) & - g(x)=  \frac{1}{2} h^2 g''(x)\int_{-1}^{1}t^2 K(t)~dt + o (h^2) +  O\Big(\frac{N_{T_n}}{h^3} \underset{0\leq j \leq 1}{\sup}\; d_{j,n}^3+ N_{T_n}\alpha_{n,h}+ n\beta_{n,h}\Big).
\end{align*}
 This concludes the proof of Proposition \ref{biastheorem}. $\Box$
\subsection{Proof of Proposition \ref{biasofwinner}.}
Let $t_0=0,~ t_{n+1}=1$ and set $\overline{Y}(t_0)=0$ and $\overline{Y}(t_{n+1})=\overline{Y}(t_n)$. Recall that,
\[\hat{g}^{pro}_n (x) = \sum_{i=1}^{n+1}\overline{Y} (t_i) \int_{t_{i-1}}^{t_i}\varphi_{x,h}(s)ds  + \sum_{i=0}^n \frac{\overline{Y} (t_{i+1})-\overline{Y} (t_i)}{t_{i+1}-t_i} \int_{t_i}^{t_{i+1}}(s-t_{i+1})\varphi_{x,h}(s)ds.\]
Since $\mathbb{E}~(\overline{Y} (t_i))=g(t_i)$ then,\[\mathbb{E}~(\hat{g}^{pro}_n (x) )= \sum_{i=1}^{n+1}g(t_i) \int_{t_{i-1}}^{t_i}\varphi_{x,h}(s)ds  + \sum_{i=0}^n \frac{g(t_{i+1})-g(t_i)}{t_{i+1}-t_i} \int_{t_i}^{t_{i+1}}(s-t_{i+1})\varphi_{x,h}(s)ds.\]
Recall that $N_{T_n}=\Card~I_{x,h}  =\{i=1,\cdots,n~:~ [t_{i-1,n},t_{i+1,n}]\cap]x-h,x+h[ \ne \emptyset\}$ and denote by $t_{x,i}$ the points of $T_n$ for which $i \in I_{x,h}$.
Using the support of $\varphi_{x,h}$ we obtain,
\begin{align*}
\mathbb{E}~(\hat{g}^{pro}_n (x) ) &  = \sum_{i=1}^{N_{T_n}} g (t_{x,i})\int_{t_{x,i-1}}^{t_{x,i}}\varphi_{x,h}(s)ds+ \sum_{i=1}^{N_{T_n}} \frac{g (t_{x,i+1})-g (t_{x,i})}{t_{x,i+1}-t_{x,i}} \int_{t_{x,i}}^{t_{x,i+1}}(s-t_{x,i+1})\varphi_{x,h}(s)ds.
\end{align*}
Let $d_{x,i}=t_{x,i+1}-t_{x,i}$. Since $g$ is in $C^2$ and $\varphi_{x,h}$ is in $C^1$ then Taylor expansions of $g$ around $t_{x,i}$ and of $\varphi_{x,h}$ around $t_{x,i+1}$ yield,
\begin{align*}
g(t_{x,i+1}) & = g(t_{x,i}) + d_{x,i}~g'(t_{x,i})+\frac{1}{2}d_{x,i}^2 ~g''(\theta_{x,i}),\\
\varphi_{x,h}(s) & = \varphi_{x,h}(t_{x,i+1})+(s-t_{x,i+1})\varphi_{x,h}'(s_i).
\end{align*}
for some $\theta_{x,i} \in ]t_{x,i},t_{x,i+1}[$ and some $s_i \in ]s,t_{x,i+1}[$. Recall that, using the support of $\varphi$,   $\varphi_{x,h}(t_{x,1})=\varphi_{x,h}(t_{x,N_{T_n}})=0$ thus,
\begin{align*}
&\mathbb{E}~(\hat{g}^{pro}_n (x) )   = \sum_{i=1}^{N_{T_n}}g (t_{x,i}) \int_{t_{x,i-1}}^{t_{x,i}}\varphi_{x,h}(s)ds~ +\sum_{i=1}^{N_{T_n}-2} g'(t_{x,i})\varphi_{x,h}(t_{x,i+1})\int_{t_{x,i}}^{t_{x,i+1}}(s-t_{x,i+1})ds \\
&+\sum_{i=1}^{N_{T_n}} g'(t_{x,i})\int_{t_{x,i}}^{t_{x,i+1}}(s-t_{x,i+1})^2\varphi_{x,h}'(s_i)ds+\frac{1}{2}\sum_{i=1}^{N_{T_n}}\varphi_{x,h}(t_{x,i+1})g''(\theta_{x,i})d_{x,i}~\int_{t_{x,i}}^{t_{x,i+1}}(s-t_{x,i+1})ds \\
&+\frac{1}{2}\sum_{i=1}^{N_{T_n}}g''(\theta_{x,i})d_{x,i}~\int_{t_{x,i}}^{t_{x,i+1}}(s-t_{x,i+1})^2\varphi_{x,h}'(s_i)ds .
\end{align*}
Recall that $g'$ and $g''$ are bounded, Lemma \ref{lemmaregulardesign} yields  $N_{T_n}=O(nh)$ and $d_{x,i}=O(\frac{1}{n})$ and using \eqref{phiboundries} we obtain,
\begin{align*}
& \sum_{i=1}^{N_{T_n}} g'(t_{x,i})\int_{t_{x,i}}^{t_{x,i+1}}(s-t_{x,i+1})^2\varphi_{x,h}'(s_i)ds = O\bigg(\frac{1}{n^2h}\bigg).\\
& \frac{1}{2}\sum_{i=1}^{N_{T_n}}\varphi_{x,h}(t_{x,i+1})g''(\theta_{x,i})d_{x,i}\int_{t_{x,i}}^{t_{x,i+1}}(s-t_{x,i+1})ds=O\bigg(\frac{1}{n^2}\bigg).\\
&\frac{1}{2}\sum_{i=1}^{N_{T_n}}g''(\theta_{x,i})d_{x,i}\int_{t_{x,i}}^{t_{x,i+1}}(s-t_{x,i+1})^2\varphi_{x,h}'(s_i)ds=O\bigg(\frac{1}{n^3h}\bigg).
\end{align*}
It follows that by simple integration,
\begin{align*}
\mathbb{E}~(\hat{g}^{pro}_n (x) ) &  = \sum_{i=1}^{N_{T_n}}g (t_{x,i}) \int_{t_{x,i-1}}^{t_{x,i}}\varphi_{x,h}(s)ds -\frac{1}{2}\sum_{i=1}^{N_{T_n}-2} g'(t_{x,i})\varphi_{x,h}(t_{x,i+1})d_{x,i}^2 +O\bigg(\frac{1}{n^2h}\bigg)\\
&  = \sum_{i=1}^{N_{T_n}} \int_{t_{x,i-1}}^{t_{x,i}}\varphi_{x,h}(s)g(s)~ds+\sum_{i=1}^{N_{T_n}} \int_{t_{x,i-1}}^{t_{x,i}}\varphi_{x,h}(s)(g (t_{x,i})-g(s))~ds\\
&~~~~~~-\frac{1}{2}\sum_{i=1}^{N_{T_n}-2} g'(t_{x,i})\varphi_{x,h}(t_{x,i+1})d_{x,i}^2 +O\bigg(\frac{1}{n^2h}\bigg).
\end{align*}
On the one hand, we have, \[\sum_{i=1}^{N_{T_n}} \int_{t_{x,i-1}}^{t_{x,i}}\varphi_{x,h}(s)g(s)~ds = \int_{x-h}^{x+h}\varphi_{x,h}(s)g(s)~ds. \]
On the other hand, Taylor expansion of $g$ and $\varphi_{x,h}$ arround $t_{x,i}$ yield,
\begin{align*}
g(t_{x,i}) & = g(s) + (t_{x,i}-s)g'(t_{x,i})-\frac{1}{2}(t_{x,i}-s)^2g''(s_i'),\\
\varphi_{x,h}(s) & = \varphi_{x,h}(t_{x,i})+(s-t_{x,i})\varphi_{x,h}'(s_i'').
\end{align*}
for some $s_i'$ and $s_i''$ in $]s,t_{x,i}[$. Thus,
\begin{align*}
&\mathbb{E}~(\hat{g}^{pro}_n (x) )  = \int_{x-h}^{x+h}\varphi_{x,h}(s)g(s)~ds+\sum_{i=2}^{N_{T_n}-1}g'(t_{x,i})\varphi_{x,h}(t_{x,i})\int_{t_{x,i-1}}^{t_{x,i}}(t_{x,i}-s)~ds\\
&~~-\sum_{i=1}^{N_{T_n}}g'(t_{x,i}) \int_{t_{x,i-1}}^{t_{x,i}} (t_{x,i}-s)^2\varphi_{x,h}'(s_i'')~ds-\frac{1}{2}\sum_{i=1}^{N_{T_n}}\varphi_{x,h}(t_{x,i}) \int_{t_{x,i-1}}^{t_{x,i}}g''(s_i')(t_{x,i}-s)^2  ~ds\\
&~~+\frac{1}{2} \sum_{i=1}^{N_{T_n}} \int_{t_{x,i-1}}^{t_{x,i}}g''(s_i') \varphi_{x,h}'(s_i'') (t_{x,i}-s)^3~ds-\frac{1}{2}\sum_{i=1}^{N_{T_n}-2} g'(t_{x,i})\varphi_{x,h}(t_{x,i+1})d_{x,i}^2 \\&~~+O\bigg(\frac{1}{n^2h}\bigg).
\end{align*}
Using the boundedness of  $g'$ and $g''$ in addition to Lemma \ref{lemmaregulardesign}  and Equation  \eqref{phiboundries}, we obtain,
\begin{align*}
& \sum_{i=1}^{N_{T_n}}g'(t_{x,i}) \int_{t_{x,i-1}}^{t_{x,i}} (t_{x,i}-s)^2\varphi_{x,h}'(s_i'')~ds = O\bigg(\frac{1}{n^2h}\bigg).\\
& \frac{1}{2}\sum_{i=1}^{N_{T_n}}\varphi_{x,h}(t_{x,i}) \int_{t_{x,i-1}}^{t_{x,i}}g''(s_i')(t_{x,i}-s)^2  ~ds = O\bigg(\frac{1}{n^2}\bigg).\\
& \frac{1}{2} \sum_{i=1}^{N_{T_n}} \int_{t_{x,i-1}}^{t_{x,i}}g''(s_i') \varphi_{x,h}'(s_i'') (t_{x,i}-s)^3~ds= O\bigg(\frac{1}{n^3h}\bigg).
\end{align*}
Thus,
\begin{align*}
\mathbb{E}~(\hat{g}^{pro}_n (x) )  &= \int_{x-h}^{x+h}\varphi_{x,h}(s)g(s)~ds+\frac{1}{2}\sum_{i=2}^{N_{T_n}-2}g'(t_{x,i})\varphi_{x,h}(t_{x,i})d_{x,i-1}^2\\
&~~~~-\frac{1}{2}\sum_{i=1}^{N_{T_n}-2} g'(t_{x,i})\varphi_{x,h}(t_{x,i+1})d_{x,i}^2+O\bigg(\frac{1}{n^2h}\bigg)\\
&  = \int_{x-h}^{x+h}\varphi_{x,h}(s)g(s)~ds+\frac{1}{2}\sum_{i=1}^{N_{T_n}-2}\Big(g'(t_{x,i+1})-g'(t_{x,i})\Big)\varphi_{x,h}(t_{x,i+1})d_{x,i}^2+O\bigg(\frac{1}{n^2h}\bigg).
\end{align*}
Since $g'$ is Lipschitz, then we have,
\begin{align}
\mathbb{E}~(\hat{g}^{pro}_n (x) ) &= \int_{x-h}^{x+h} \varphi_{x,h}(s)\;g (s)\;ds+ O\bigg(\frac{1}{n^2h}\bigg).\label{espewinner}
\end{align}
Finally, from \eqref{chepa} we obtain,
\begin{align*}
\mathbb{E}~(\hat{g}^{pro}_n (x) ) -g(x) & = \frac{1}{2} h^2 g''(x)\int_{-1}^{1} t^2K(t) dt + o (h^2)+ O\bigg(\frac{1}{n^2h}\bigg).
\end{align*}
This concludes the proof of Proposition \ref{biasofwinner}. $\Box$
\subsection{Proof of Proposition \ref{propoflimit}.}
The great lines of this proof are based on  Sacks and Ylvisaker (1966) \cite{Sacks and Ylvisaker}. From the definition of the orthogonal projection (see the Appendix) and using the Pythagore theorem we obtain,
\begin{equation}
m\Big(\frac{\sigma_{x,h}^2}{m}-\Var g_n^{pro}(x)\Big)=||f_{x,h}||^2 -||P_{|T_n}f_{x,h}||^2=||f_{x,h} - P_{|T_n}f_{x,h}||^2, \label{devaranorm}
\end{equation}
where $P_{|T_n}f_{x,h}$ is the orthogonal projection of $f_{x,h}$ on the subspace of $\mathcal F(\varepsilon)$ spanned by $\{R(\cdot,t_i), t_i \in T_n\}$, denoted here by $V_{T_n}$.
 We shall then prove that,
\begin{align}
||f_{x,h} - P_{|T_n}f_{x,h}||^2  \leq \frac{C}{h}\;\underset{0\leq j\leq n }{\sup}\;d_{j,n}^2. \label{lamajoration}
\end{align}
Recall that $N_{T_n}=\Card~I_{x,h} = \Card~I_{x,h} =\{i=1,\cdots,n~:~ [t_{i-1,n},t_{i+1,n}]\cap]x-h,x+h[ \ne \emptyset\}$ and denote by $t_{x,i}$ the points of $T_n$ for which $i \in I_{x,h}$. Let $g_n:=g_{n,x} =\sum_{i=1}^{n} \gamma_{x,i} R(\cdot,t_{x,i})$ with $\gamma_{x,i}=0$ for every $i \notin I_{x,h}$. It is clear that  $g_n\in V_{T_n}$ and thus from the definition of the orthogonal projection we have,
\[||f_{x,h} - P_{|T_n}f_{x,h}||^2  \leq ||f_{x,h} -g_n||^2.\]
Now using (F1) in the Appendix and the support of $\varphi_{x,h}$  we obtain,
\begin{align}
||f_{x,h} -g_n||^2& = \int_{0}^{1}(f_{x,h}(t) -g_n(t))\varphi_{x,h}(t)\;dt - \sum_{i=1}^{n}(f_{x,h}(t_i) -g_n(t_i))\gamma_{x,i}\nonumber\\
& = \int_{x-h}^{x+h}(f_{x,h}(t) -g_n(t))\varphi_{x,h}(t)\;dt - \sum_{i=1}^{N_{T_n}}(f_{x,h}(t_{x,i}) -g_n(t_{x,i}))\gamma_{x,i}\label{casgen}
\end{align}
In what follows, we distinguish between three cases according to the location of $t_{x,1}$ and $t_{x,N_{T_n}}$ in the interval $[x-h,x+h]$.\\
\textbf{First case.} Suppose first that $t_{x,1}=x-h$ and $t_{x,N_{T_n}}=x+h$ and take, \begin{equation}\gamma_{x,i}=\int_{t_{x,i}}^{t_{x,{i+1}}} \varphi_{x,h}(t)\;dt~~ \text{for}~~ i =1,\dots ,N_{T_n}-1\label{gammafprappendix}.
\end{equation}
we have in this case,
\begin{align}
||f_{x,h} -g_n||^2
& = \sum_{i=1}^{N_{T_n}} \int_{t_{x,i}}^{t_{x,{i+1}}} \bigg( (f_{x,h}(t) -g_n(t))-(f_{x,h}(t_{x,i}) -g_n(t_{x,i})) \bigg)\varphi_{x,h}(t)\;dt.\label{okhhhh}
\end{align}
Assumption $(A)$ yields that $f_{x,h}$ is twice differentiable on $[0,1]$, while $g_n$ is twice differentiable everywhere except on $T_n$, but it has left and right derivatives.
Taylor expansion of $f_{x,h} -g_n$  around $t_{x,i}$ for $i=1,\cdots,N_{T_n}-1$ and $t \in ]t_{x,i},t_{x,i+1}[$ gives,
\begin{align}
f_{x,h}(t) -g_n(t) &= (f_{x,h}(t_{x,i}) -g_n(t_{x,i}))+(t-t_{x,i})(f_{x,h}'(t_{x,i}) -g_n'(t_{x,i}^+)) \nonumber \\
&~~~~+\frac{1}{2}(t-t_{x,i})^2(f_{x,h}''(\theta_{x,t}) -g_n''(\theta_{x,t}^+)),
\label{expansion}
\end{align}
for some $\theta_{x,t} \in ]t_{x,i},t[$. On the one hand, we have,
\begin{equation}
g_n'(t_{x,i}^+) =\sum_{j=1}^{N_{T_n}-1} R^{(0,1)}(t_{x,j},t_{x,i}^+)\gamma_{x,j}.
\label{g'}
\end{equation}
On the other hand, using \eqref{f't} we obtain,
\begin{align}
& f_{x,h}'(t_{x,i})  = \int_{x-h}^{x+h} R^{(0,1)}(s,t_{x,i}^+)\varphi_{x,h}(s)\;ds \nonumber = \sum_{j=1}^{N_{T_n}-1} \int_{t_{x,j}}^{t_{x,{j+1}}} R^{(0,1)}(s,t_{x,i}^+)\varphi_{x,h}(s)\;ds\nonumber \\
& = \underset{j \ne i}{\sum_{j=1}^{N_{T_n}-1}} \int_{t_{x,j}}^{t_{x,{j+1}}} R^{(0,1)}(s,t_{x,i}^+)\varphi_{x,h}(s)\;ds\;+\;\int_{t_{x,i}}^{t_{x,{i+1}}} R^{(0,1)}(s,t_{x,i}^+)\varphi_{x,h}(s)\;ds. \label{f'i}
\end{align}
When $j \ne i$ we have,
{\small
\begin{align}
& \int_{t_{x,j}}^{t_{x,{j+1}}} R^{(0,1)}(s,t_{x,i}^+)\varphi_{x,h}(s)\;ds  =R^{(0,1)}(t_{x,j},t_{x,i})\gamma_{x,j}+\int_{t_{x,j}}^{t_{x,{j+1}}}(s-t_{x,j}) R^{(1,1)}(\delta_{s,j},t_{x,i})\varphi_{x,h}(s)\;ds, \label{inej}
\end{align}}
for some $\delta_{s,j} \in ]t_{x,j},s[$, while for $j=i$ we have,
\begin{align}
&\int_{t_{x,i}}^{t_{x,{i+1}}} R^{(0,1)}(s,t_{x,i}^+)\varphi_{x,h}(s)\;ds = \int_{t_{x,i}}^{t_{x,{i+1}}} R^{(0,1)}(s,t_{x,i}^-)\varphi_{x,h}(s)\;ds \nonumber \\
&~~~~~~ = R^{(0,1)}(t_{x,i},t_{x,i}^-)\gamma_{x,i}+\int_{t_{x,i}}^{t_{x,{i+1}}}(s-t_{x,i}) R^{(1,1)}(\delta_{s,i}^+,t_{x,i}^-)\varphi_{x,h}(s)\;ds. \label{i=j}
\end{align}
Collecting \eqref{g'}, \eqref{f'i}, \eqref{inej} and \eqref{i=j} we obtain,
\begin{align*}
f_{x,h}'&(t_{x,i}) -g_n'(t_{x,i}^+)=\underset{j \ne i}{\sum_{j=1}^{N_{T_n}-1}} R^{(0,1)}(t_{x,j},t_{x,i})\gamma_{x,j}+\underset{j \ne i}{\sum_{j=1}^{N_{T_n}-1}}\int_{t_{x,j}}^{t_{x,{j+1}}}(s-t_{x,j}) R^{(1,1)}(\delta_{s,j},t_{x,i})\varphi_{x,h}(s)\;ds\\
&~~+ R^{(0,1)}(t_{x,i},t_{x,i}^-)\gamma_{x,i}+\int_{t_{x,i}}^{t_{x,{i+1}}} R^{(1,1)}(\delta_{s,i}^+,t_{x,i}^-)\varphi_{x,h}(s)\;ds -\sum_{j=1}^{N_{T_n}-1} R^{(0,1)}(t_{x,j},t_{x,i}^+)\gamma_{x,j}\\
& = \alpha(t_{x,i})\gamma_{x,i}+\sum_{j=1}^{N_{T_n}-1}\int_{t_{x,j}}^{t_{x,{j+1}}} (s-t_{x,j})R^{(1,1)}(\delta_{s,j}^+,t_{x,i}^-)\varphi_{x,h}(s)\;ds.
\end{align*}
It is easy to see that,
\begin{align}
|f_{x,h}'(t_{x,i}) -g_n'(t_{x,i}^+)| & \leq \alpha_1 \gamma_{x,i} + \frac{K_{\infty}}{h} R_1 \sum_{j=1}^{N_{T_n}-1}\int_{t_{x,j}}^{t_{x,{j+1}}} (s-t_{x,j})\;ds \nonumber\\
& \leq \frac{K_{\infty}}{h} \alpha_1 d_{x,i}  +  \frac{K_{\infty}}{2h} R_1 \sum_{j=1}^{N_{T_n}-1} d_{x,j}^2. \label{dif'}
\end{align}
 We deduce from \eqref{f''t}  that for all $\theta_{x,t} \in ]t_{x,i},t_{x,{i+1}}[$ we have,
\begin{equation*}
|f_{x,h}''(\theta_{x,t})|\leq \frac{K_{\infty}}{h} \alpha_1 + \frac{K_{\infty}}{h} R_2 \times 2h = \frac{K_{\infty}}{h} \alpha_1 + 2K_{\infty} R_2.
\end{equation*}
In addition, for $\theta_{x,t} \in ]t_{x,i},t_{x,{i+1}}[$ we have,
\begin{align*}
|g_n''(\theta_{x,t}^+)| & = \bigg|\sum_{j=1}^{N_{T_n}-1} R^{(0,2)}(t_{x,j},\theta_{x,t}^+)\gamma_{x,j} \bigg| \leq \frac{K_{\infty}}{h} R_2 \sum_{j=1}^{N_{T_n}-1} d_{x,j}= \frac{K_{\infty}}{h} R_2 \times 2h = 2K_{\infty} R_2,
\end{align*} Thus,
\begin{equation}
|f_{x,h}''(\theta_{x,t})-g_n''(\theta_{x,t}^+)| \leq \frac{K_{\infty}}{h} \alpha_1 + 4K_{\infty} R_2 \label{dif''}.
\end{equation}
Equations \eqref{expansion}, \eqref{dif'} and \eqref{dif''} yield that for  $i=1,\cdots,N_{T_n}-1,$
{\small
\begin{align}
&~~~~\bigg| \int_{t_{x,i}}^{t_{x,{i+1}}} \big[ (f_{x,h}(t) -g_n(t))-(f_{x,h}(t_{x,i}) -g_n(t_{x,i})) \big]\varphi_{x,h}(t)\;dt \bigg|\nonumber\\
& \leq \int_{t_{x,i}}^{t_{x,{i+1}}} (t-t_{x,i}) |f_{x,h}'(t_{x,i}) -g_n'(t_{x,i}^+)||\varphi_{x,h}(t)|\;dt\nonumber\\
& \;\;\;\;\;\;\;\;\; + \frac{1}{2} \int_{t_{x,i}}^{t_{x,{i+1}}} (t-t_{x,i}) ^2 |f_{x,h}''(\theta_{x,t})-g_n''(\theta_{x,t}^+)||\varphi_{x,h}(t)|\;dt\nonumber\\
& \leq \bigg( \frac{K_{\infty}}{h} \alpha_1 d_{x,i}  +  \frac{K_{\infty}}{2h} R_1 \sum_{j=1}^{N_{T_n}-1} d_{x,j}^2 \bigg) \int_{t_{x,i}}^{t_{x,{i+1}}} (t-t_{x,i})|\varphi_{x,h}(t)|\;dt \nonumber\\
& \;\;\;\;\;\;\;\;\;+ \frac{1}{2} \bigg( \frac{K_{\infty}}{h} \alpha_1 + 4K_{\infty} R_2\bigg)\int_{t_{x,i}}^{t_{x,{i+1}}} (t-t_{x,i}) ^2|\varphi_{x,h}(t)|\;dt \nonumber\\
& \leq \bigg( \frac{K_{\infty}}{h} \alpha_1 d_{x,i}  +  \frac{K_{\infty}}{2h} R_1 \sum_{j=1}^{N_{T_n}-1} d_{x,j}^2 \bigg)\frac{K_{\infty}}{2h} d_{x,i}^2 +  \frac{1}{2} \bigg( \frac{K_{\infty}}{h} \alpha_1 + 4K_{\infty} R_2\bigg)\frac{K_{\infty}}{3h} d_{x,i}^3 \nonumber\\
& \leq \frac{K_{\infty}^2}{4h^2}R_1 d_{x,i}^2 \sum_{j=1}^{N_{T_n}-1} d_{x,j}^2 +\frac{2K_{\infty}^2}{3h}\Big(\frac{\alpha_1}{h}+R_2\Big) d_{x,i}^3.\label{cases}
\end{align}}

Injecting this inequality in \eqref{okhhhh} yields,
\begin{align*}
&||f_{x,h} - P_{|T_n}f_{x,h}||^2  \leq   \frac{K_{\infty}^2}{4h^2}R_1 \bigg(\sum_{i=1}^{N_{T_n}-1} d_{x,i}^2\bigg)^2 +\frac{2K_{\infty}^2}{3h}\Big(\frac{\alpha_1}{h}+R_2\Big)  \sum_{i=1}^{N_{T_n}-1} d_{x,i}^3\\
& \leq   \frac{K_{\infty}^2}{4h^2}R_1 ~\underset{1 \leq i\leq n}{\sup}\;d_{i,n}^2\bigg(\sum_{i=1}^{N_{T_n}-1} d_{x,i}\bigg)^2 +\frac{2K_{\infty}^2}{3h}\big(\frac{\alpha_1}{h}+R_2\big)~\underset{1 \leq i\leq n}{\sup}\;d_{i,n}^2  \sum_{i=1}^{N_{T_n}-1} d_{x,i}.
\end{align*}
Since $\sum_{i=1}^{N_{T_n}-1} d_{x,i}=2h$ then,
\[||f_{x,h} - P_{|T_n}f_{x,h}||^2  \leq  \bigg(\frac{4}{3h}\alpha_1+ R_1 +\frac{4}{3}R_2  \bigg)K_{\infty}^2\underset{1 \leq i\leq n}{\sup}\;d_{i,n}^2\]
Finally, since $h<1$ then,
\begin{align*}
&||f_{x,h} - P_{|T_n}f_{x,h}||^2  \leq \bigg(\frac{4}{3}\alpha_1+ R_1 +\frac{4}{3}R_2  \bigg)K_{\infty}^2\frac{1}{h} \underset{1 \leq i\leq n}{\sup}\;d_{i,n}^2. 
\end{align*}
Proposition \ref{propoflimit} is then proved for the first case.\\[0.2cm]
\textbf{Second case.}  Consider now the case where  $t_{x,1} < x-h$ and $t_{x,{N_{T_n}}} > x+h$. For $i=2,\dots,N_{T_n}-2$ set,
{\small
\begin{equation}
\gamma_{x,i}=\int_{t_{x,i}}^{t_{x,{i+1}}} \varphi_{x,h}(t)\;dt,~
\gamma_{x,1}=\int_{x-h}^{t_{x,2}} \varphi_{x,h}(t)\;dt,~\gamma_{x,N_{T_n}-1}=\int_{t_{x,N_{T_n}-1}}^{x+h} \varphi_{x,h}(t)\;dt~\text{and}~\gamma_{x,N_{T_n}}=0.
\end{equation}}
Using this we obtain,
\begin{align}
||f_{x,h} - g_n||^2 
& = \int_{x-h}^{t_{x,2}}\Big( (f_{x,h}(t) - g_n (t))-(f_{x,h}(t_{x,1}) - g_n (t_{x,1}))\Big)\varphi_{x,h}(t) \;dt \nonumber \\
&~ + \sum_{i=2}^{N_{T_n}} \int_{t_{x,i}}^{t_{x,{i+1}}}\Big( (f_{x,h}(t) - g_n (t))-(f_{x,h}(t_{x,i}) - g_n (t_{x,i}))\Big)\varphi_{x,h}(t) \;dt\nonumber \\
&~+\int_{t_{x,N_{T_n}}}^{x+h}\Big( (f_{x,h}(t) - g_n (t))-(f_{x,h}(t_{x,N_{T_n}}) - g_n (t_{x,N_{T_n}}))\Big)\varphi_{x,h}(t)\;dt. \label{diff}
\end{align}
We first control the first term of \eqref{diff}. Let,\begin{align*}
A_{x,h}^{(1)} & = \int_{x-h}^{t_{x,2}}\Big( (f_{x,h}(t) - g_n (t))-(f_{x,h}(t_{x,1}) - g_n (t_{x,1}))\Big)\varphi_{x,h}(t) \;dt.
\end{align*}
For $t \in ]x-h,t_{x,2}[$ we have,
\begin{align}
f_{x,h}(t) -g_n(t) &= (f_{x,h}(t_{x,1}) -g_n(t_{x,1}))+(t-t_{x,1})(f_{x,h}'(t_{x,1}) -g_n'(t_{x,1}^+)) \nonumber \\
& ~~~~+\frac{1}{2}(t-t_{x,1})^2(f_{x,h}''(\theta_{x,1}) -g_n''(\theta_{x,1}^+)),\label{expansion1}
\end{align}
for some $\theta_{x,1} \in ]x-h,t[$. Equation \eqref{f't} yields,
\begin{align}
&f_{x,h}'(t_{x,1})  = \int_{x-h}^{x+h} R^{(0,1)}(s,t_{x,1}^+)\varphi_{x,h}(s)\;ds \nonumber = \sum_{j=1}^{N_{T_n}-1} \int_{t_{x,j}}^{t_{x,{j+1}}} R^{(0,1)}(s,t_{x,1}^+)\varphi_{x,h}(s)\;ds\nonumber \\
& = \int_{x-h}^{t_{x,2}} R^{(0,1)}(s,t_{x,1}^-)\varphi_{x,h}(s)\;ds\;+{\sum_{j=2}^{N_{T_n}-1}} \int_{t_{x,j}}^{t_{x,{j+1}}} R^{(0,1)}(s,t_{x,1}^+)\varphi_{x,h}(s)\;ds\nonumber\\
& = R^{(0,1)}(t_{x,1},t_{x,1}^-)\gamma_{x,1}+\int_{x-h}^{t_{x,2}} (s-t_{x,1})R^{(1,1)}(\delta_{s,1}^+,t_{x,1}^-)\varphi_{x,h}(s)\;ds \nonumber\\
&\;\;\;\;+{\sum_{j=2}^{N_{T_n}-1}} R^{(0,1)}(t_{x,j},t_{x,1})\gamma_{x,j} \;+ \;{\sum_{j=2}^{N_{T_n}-1}} \int_{t_{x,j}}^{t_{x,{j+1}}} (s-t_{x,j})R^{(1,1)}(\delta_{s,j},t_{x,1}^+)\varphi_{x,h}(s)\;ds. \label{f't1}
\end{align}
Recall that,
\begin{equation}
g_n'(t_{x,1}^+) =R^{(0,1)}(t_{x,1},t_{x,1}^+)\gamma_{x,1}+ \sum_{j=2}^{N_{T_n}-1} R^{(0,1)}(t_{x,j},t_{x,1})\gamma_{x,j}. \label{g't1}
\end{equation}
Equations \eqref{f't1} and \eqref{g't1} give,
\begin{align*}
f_{x,h}'(t_{x,1})-g_n'(t_{x,1}^+) & = \alpha(t_{x,1})\gamma_{x,1} + {\sum_{j=2}^{N_{T_n}-1}} \int_{t_{x,j}}^{t_{x,{j+1}}} (s-t_{x,j})R^{(1,1)}(\delta_{s,j},t_{x,1}^+)\varphi_{x,h}(s)\;ds\\
&\;\;\; +\int_{x-h}^{t_{x,2}}(s-t_{x,1})R^{(1,1)}(\delta_{s,1}^+  ,t_{x,1}^-)\varphi_{x,h}(s)\;ds.
\end{align*}
Note that $t_{x,2}-(x-h)\leq \underset{1 \leq i \leq n}{\sup}~d_{i,n}$.  We obtain,
\begin{align}
|f_{x,h}'(t_{x,1})-g_n'(t_{x,1}^-) | &\leq \frac{K_{\infty}}{h}\alpha_1\underset{1 \leq i \leq n}{\sup}~d_{i,n} + \frac{K_{\infty}}{2h} R_1 {\sum_{j=2}^{N_{T_n}-1}} d_{x,j}^2+\frac{K_{\infty}}{2h} R_1 \underset{1 \leq i \leq n}{\sup}~d_{i,n}^2 \nonumber\\
& \leq \frac{K_{\infty}}{h}\alpha_1\underset{1 \leq i \leq n}{\sup}~d_{i,n} + K_{\infty}R_1\underset{1 \leq i \leq n}{\sup}~d_{i,n}+\frac{K_{\infty}}{2h} R_1 \underset{1 \leq i \leq n}{\sup}~d_{i,n}^2\nonumber\\
&\leq K_{\infty}\Big(\frac{\alpha_1}{h}+\frac{3}{2}R_1\Big)\underset{1 \leq i \leq n}{\sup}~d_{i,n} \label{dif'1}
\end{align}
By \eqref{dif''} we have,
\begin{equation}
|f_{x,h}''(\theta_{x,t})-g_n''(\theta_{x,t}^-)| \leq \frac{K_{\infty}}{h} \alpha_1 + 4K_{\infty} R_2 \label{dif''1}.
\end{equation}
Equations \eqref{expansion1}, \eqref{dif'1} and \eqref{dif''1} yield,
\begin{align}
|A_{x,h}^{(1)}| &  \leq ~~|f_{x,h}'(t_{x,1}) -g_n'(t_{x,1}^+)|\int_{x-h}^{t_{x,2}} (t-t_{x,1}) |\varphi_{x,h}(t)|\;dt \nonumber\\
&~~+ \frac{1}{2} \int_{x-h}^{t_{x,2}} (t-t_{x,1}) ^2 |f_{x,h}''(\theta_{x,1})-g_n''(\theta_{x,1}^+)||\varphi_{x,h}(t)|\;dt \nonumber \\
& \leq \bigg(K_{\infty}\Big(\frac{\alpha_1}{h}+\frac{3}{2}R_1\Big)\underset{1 \leq i \leq n}{\sup}~d_{i,n} \bigg) \frac{K_{\infty}}{2h}\underset{1 \leq i \leq n}{\sup}~d_{i,n}^2  +\Big (\frac{K_{\infty}}{h} \alpha_1 + 4K_{\infty} R_2\Big) \frac{K_{\infty}}{6h} \underset{1 \leq i \leq n}{\sup}~d_{i,n}^3 \nonumber \\
& \leq  \bigg(\frac{2}{3}\alpha_1+\frac{3}{4}R_1+\frac{2}{3}R_2\bigg)\frac{K_{\infty}^2}{h^2} \underset{1 \leq i \leq n}{\sup}~d_{i,n}^3. \label{A1}
\end{align}
\noindent Similarly we obtain,
\begin{align}
A_{x,h}^{(2)} & \overset{\Delta}{=}\int_{t_{x,N_{T_n}}}^{x+h}\bigg( (f_{x,h}(t) - g_n (t))-(f_{x,h}(t_{x,N_{T_n}}) - g_n (t_{x,N_{T_n}}))\bigg)\varphi_{x,h}(t)\;dt\nonumber\\
|A_{x,h}^{(2)}| & \leq  \bigg(\frac{2}{3}\alpha_1+\frac{3}{4}R_1+\frac{2}{3}R_2\bigg)\frac{K_{\infty}^2}{h^2} \underset{1 \leq i \leq n}{\sup}~d_{i,n}^3.\label{A2}
\end{align}
Thus,
\begin{align*}
|A_{x,h}^{(1)}+A_{x,h}^{(2)}| & \leq  \bigg(\frac{4}{3}\alpha_1+\frac{3}{2}R_1+\frac{4}{3}R_2\bigg)\frac{K_{\infty}^2}{h^2} \underset{1 \leq i \leq n}{\sup}~d_{i,n}^3.
\end{align*}
For $i = 2,\dots,N_{T_n}-2$, similar calculations as those leading to \eqref{cases} give,
\begin{align*}
&\bigg| \int_{t_{x,i}}^{t_{x,{i+1}}} \big( (f_{x,h}(t) -g_n(t))-(f_{x,h}(t_{x,i}) -g_n(t_{x,i})) \big)\varphi_{x,h}(t)\;dt \bigg|\\
& \leq \frac{K_{\infty}^2}{4h^2}R_1 d_{x,i}^2 \sum_{j=1}^{N_{T_n}} d_{x,j}^2 +\frac{2K_{\infty}^2}{3h}(\frac{\alpha_1}{h}+R_2) d_{x,i}^3.
\end{align*}
Thus,
\begin{align}
& \bigg| \sum_{i=2}^{N_{T_n}-2} \int_{t_{x,i}}^{t_{x,{i+1}}} \bigg( (f_{x,h}(t) -g_n(t))-(f_{x,h}(t_{x,i}) -g_n(t_{x,i})) \bigg)\varphi_{x,h}(t)\;dt \bigg|\nonumber \\
& \leq \bigg(\frac{4}{3}\alpha_1+ R_1 +\frac{4}{3}R_2  \bigg)\frac{K_{\infty}^2}{h} \underset{1 \leq i \leq n}{\sup}~d_{i,n}^2.
\label{jsup2}
\end{align}
Then, Equations \eqref{A1}, \eqref{A2} and \eqref{jsup2} yield,
\begin{align*}
|| f_{x,h} - P_{T_n}f_{x,h}  ||^2 & \leq  \bigg(\frac{4}{3}\alpha_1+\frac{3}{2}R_1+\frac{4}{3}R_2\bigg)\frac{K_{\infty}^2}{h} \underset{1 \leq i \leq n}{\sup}~d_{i,n}^2+\bigg(\frac{4}{3}\alpha_1+ R_1 +\frac{4}{3}R_2  \bigg)\frac{K_{\infty}^2}{h^2} \underset{1 \leq i \leq n}{\sup}~d_{i,n}^3\\
& =  \bigg(\frac{8}{3}\alpha_1+\frac{5}{2}R_1+\frac{8}{3}R_2\bigg)\frac{K_{\infty}^2}{h} \underset{1 \leq i \leq n}{\sup}~d_{i,n}^2
\end{align*}
\textbf{Third case.} Suppose now that $t_{x,1} = x-h$ and $t_{x,{N_{T_n}}} > x+h$ (respectively $t_{x,1} < x-h$ and $t_{x,{N_{T_n}}} = x+h$). Let $T_{n-1} = T_n -\{x-h\}$ (respectively $T_{n-1} = T_n -\{x+h\}$). Since $P_{T_{n-1}}f_{x,h} \in V_{T_n}$ we obtain,  \[||f_{x,h} - P_{T_n}f_{x,h}||^2 \leq ||f_{x,h} - P_{T_{n-1}}f_{x,h}||^2,\]
we can then apply the result of the second case to the right side of the previous inequality.
The proof of Proposition \ref{propoflimit} is complete. $\Box$
\subsection{Proof of Proposition \ref{variancetheorem}.}
The great lines of this proof are based on the work of Sacks and Ylvisaker (1966) \cite{Sacks and Ylvisaker}. Keeping  Equation \eqref{devaranorm} in mind we deduce  that Equation \eqref{varianceformula} is equivalent to,
\begin{equation}
\underset{n \to \infty}{\underline \lim} \frac{N_{T_n}^2}{h} ||f_{x,h} - P_{|T_n}f_{x,h}||^2 \geq \frac{1}{12} \alpha(x) \bigg\{ \int_{-1}^{1} K^{2/3}(t)dt \bigg\}^3.
\label{liminf}
\end{equation}We shall take the same notation as in the previous proof. Let $g_n=P_{|T_n}f_{x,h}$,  it is shown by Equation \eqref{mvar=|||PT_n|} in the Appendix  that:  \[ g_n(t_i)=f_{x,h}(t_i)=\sum_{j=1}^n R(t_j,t_i)m_{x,h}(t_j),~~\text{for}~i=1,\cdots,n.\] We have from  (F1) in the Appendix  that,
\begin{align}
&||f_{x,h}-g_n||^2 =\int_{0}^{1} (f_{x,h}(t) -g_n(t))\varphi_{x,h} (t)\;dt- \sum_{i=1}^{n} m_{x,h}(t_i) ( f_{x,h}(t_{i}) -g_n(t_{i}))\nonumber\\
& =\int_{x-h}^{x+h} (f_{x,h}(t)-g_n(t))\varphi_{x,h}(t)dt.\nonumber
\end{align}
 Suppose first that  $t_{x,1}=x-h$ and $t_{x,{N_{T_n}}}=x+h$, then the last equalities give,
\begin{align}
&||f_{x,h}-g_n||^2  = \sum_{i=1}^{N_{T_n}-1}\int_{t_{x,i}}^{t_{x,{i+1}}}(f_{x,h}(t)-g_n(t))\varphi_{x,h}(t)dt. \label{norm'}
\end{align}
Under Assumptions $(A)$ and $(B)$, the function $f_{x,h}$ is twice differentiable at every $t \in [0,1]$ and $g_n$ is twice differentiable at every $t \in [0,1]$ except on $T_n$, however, it has left and right derivatives.
We expand $(f_{x,h}-g_n)$ in a Taylor series around $t_{x,i}$ for $t \in ]t_{x,i},t_{x,{i+1}}[$ up to order 2 we obtain,
\begin{align*}
f_{x,h}(t)-g_n(t)& = (f_{x,h}(t_{x,i})-g_n(t_{x,i}))+(t-t_{x,i})(f'_{x,h}(t_{x,i})-g_n'(t_{x,i}^+)) \nonumber \\
& \;\;\;\;+\frac{1}{2}(t-t_{x,i})Â²(f_{x,h}''(\sigma_{x,t})-g_n''(\sigma_{x,t}^+)),
\end{align*}
for some $\sigma_{x,t} \in ]t_{x,i},t[$.  Since $g_n(t_{x,i})=f_{x,h}(t_{x,i})$ then,
\begin{align}
f_{x,h}(t)-g_n(t)& = (t-t_{x,i})(f'_{x,h}(t_{x,i})-g_n'(t_{x,i}^+))+\frac{1}{2}(t-t_{x,i})Â²(f_{x,h}''(\sigma_{x,t})-g_n''(\sigma_{x,t}^+)), \label{expansion2}
\end{align}
On the one hand, we have for $i \in 1,\dots,N_{T_n}-1$,
\begin{align*}
f_{x,h}(t_{x,{i+1}})-g_n(t_{x,{i+1}})& = d_{x,i}(f'_{x,h}(t_{x,i})-g_n'(t_{x,i}^+))+\frac{1}{2}d_{x,i}^2(f_{x,h}''(\sigma_{x,i})-g_n''(\sigma_{x,i}^+)).
\end{align*}
for some $\sigma_{x,i} \in ]t_{x,i},t_{x,{i+1}}[$. Thus,
\begin{equation}
f'_{x,h}(t_{x,i})-g_n'(t_{x,i}^+)=-\frac{1}{2}d_{x,i}(f_{x,h}''(\sigma_{x,i})-g_n''(\sigma_{x,i}^+)).\label{diff'}
\end{equation}
On the other  hand, it is shown by (F4) in the Appendix  that,
\begin{equation}
 f_{x,h}''(t)-g_n''(t^+) = -\alpha(t) \varphi_{x,h}(t)+\langle R^{(0,2)}(\cdot,t^+),f_{x,h}-g_n\rangle. \label{diff''}
\end{equation}Injecting \eqref{diff'} and \eqref{diff''} in \eqref{expansion2}  gives,
\begin{align*}
& f_{x,h}(t)-g_n(t) =-\frac{1}{2} (t-t_{x,i})d_{x,i}(f_{x,h}''(\sigma_{x,i})-g_n''(\sigma_{x,i}^+))+\frac{1}{2}(t-t_{x_i})^2 (f_{x,h}''(\sigma_{x,i})-g_n''(\sigma_{x,i}^+))\\
& =\frac{1}{2} d_{x,i}(t-t_{x,i})\alpha(\sigma_{x,i}) \varphi_{x,h}(\sigma_{x,i})-\frac{1}{2}(t-t_{x,i})^²\alpha(\sigma_{x,t}) \varphi_{x,h}(\sigma_{x,t})\\
&~~-\frac{1}{2}d_{x,i}(t-t_{x,i}) \langle R^{(0,2)}( \cdot,\sigma_{x,i}^+),f_{x,h}-g_n\rangle+\frac{1}{2}(t-t_{x,i})^2 \langle R^{(0,2)}(\cdot ,\sigma_{x,t}^+),f_{x,h}-g_n  \rangle .
\end{align*}
\noindent Thus,
\begin{align}
&\int_{t_{x,i}}^{t_{x,{i+1}}}(f_{x,h}(t)-g_n(t))\varphi_{x,h}(t)\;dt=\nonumber\\
&\frac{1}{2}d_{x,i}\alpha(\sigma_{x,i}) \varphi_{x,h}(\sigma_{x,i})\int_{t_{x,i}}^{t_{x,{i+1}}}(t-t_{x,i})\varphi_{x,h}(t)\;dt-\frac{1}{2}\int_{t_{x,i}}^{t_{x,{i+1}}}(t-t_{x,i})Â²\alpha(\sigma_{x,t}) \varphi_{x,h}(\sigma_{x,t})\varphi_{x,h}(t)\;dt\nonumber\\
& ~~-\frac{1}{2}d_{x,i}\langle R^{(0,2)}(\cdot,\sigma_{x,i}^+),f_{x,h}-g_n\rangle \int_{t_{x,i}}^{t_{x,{i+1}}}(t-t_{x,i})\varphi_{x,h}(t)\;dt\nonumber\\
&~~+\frac{1}{2}\int_{t_{x,i}}^{t_{x,{i+1}}}(t-t_{x,i})^2\langle R^{(0,2)}(\cdot,\sigma_{x,t}^+),f_{x,h}-g_n\rangle \varphi_{x,h}(t)\;dt\nonumber\\
& = \frac{1}{4}d_{x,i}^3\alpha(\sigma_{x,i}) \varphi_{x,h}^2(\sigma_{x,i})-\frac{1}{6}d_{x,i}^3\alpha(\sigma_{x,i}) \varphi_{x,h}^2(\sigma_{x,i}) \nonumber \\
&~~+\frac{1}{2}d_{x,i}\alpha(\sigma_{x_i}) \varphi_{x,h}(\sigma_{x,i})\int_{t_{x,i}}^{t_{x,{i+1}}}(t-t_{x,i})[\varphi_{x,h}(t)-\varphi_{x,h}(\sigma_{x,i})]\;dt \nonumber \\
&~~-\frac{1}{2}\alpha(\sigma_{x,i}) \varphi_{x,h}(\sigma_{x,i})\int_{t_{x,i}}^{t_{x,{i+1}}}(t-t_{x,i})^2[\varphi_{x,h}(t)-\varphi_{x,h}(\sigma_{x,i})]\;dt \nonumber\\
&~~-\frac{1}{2}\int_{t_{x,i}}^{t_{x,{i+1}}}(t-t_{x,i})^2[\alpha(\sigma_{x,t})\varphi_{x,h}(\sigma_{x,t})-\alpha(\sigma_{x,i})\varphi_{x,h}(\sigma_{x,i})]\varphi_{x,h}(t)\;dt \nonumber\\
&~~ -\frac{1}{2}d_{x,i}\langle R^{(0,2)}(\cdot,\sigma_{x,i}^+),f_{x,h}-g_n\rangle \int_{t_{x,i}}^{t_{x,{i+1}}}(t-t_{x,i})\varphi_{x,h}(t)\;dt \nonumber\\
&~~+\frac{1}{2}\int_{t_{x,i}}^{t_{x,{i+1}}}(t-t_{x,i})^2\langle R^{(0,2)}(\cdot,\sigma_{x,t}^+),f_{x,h}-g_n\rangle\varphi_{x,h}(t)\;dt \nonumber\\[0.3cm]
& = \frac{1}{12}d_{x,i}^3\alpha(\sigma_{x,i}) \varphi_{x,h}^2(\sigma_{x,i})+A_{x,i}^{(1)}-A_{x,i}^{(2)}-A_{x,i}^{(3)}-A_{x,i}^{(4)}+A_{x,i}^{(5)}, \label{touse}
\end{align}

\begin{align*}
\text{where,} ~~~~A_{x,i}^{(1)}&=\frac{1}{2}d_{x,i}\alpha(\sigma_{x,i}) \varphi_{x,h}(\sigma_{x,i})\int_{t_{x,i}}^{t_{x,{i+1}}}(t-t_{x,i})[\varphi_{x,h}(t)-\varphi_{x,h}(\sigma_{x,i})]\;dt.\\
A_{x,i}^{(2)} & =\frac{1}{2}\alpha(\sigma_{x,i}) \varphi_{x,h}(\sigma_{x,i})\int_{t_{x,i}}^{t_{x,{i+1}}}(t-t_{x,i})^2[\varphi_{x,h}(t)-\varphi_{x,h}(\sigma_{x,i})]\;dt.\\
A_{x,i}^{(3)} & = \frac{1}{2}\int_{t_{x,i}}^{t_{x,{i+1}}}(t-t_{x,i})^2[\alpha(\sigma_{x,t})\varphi_{x,h}(\sigma_{x,t})-\alpha(\sigma_{x,i})\varphi_{x,h}(\sigma_{x,i})]\varphi_{x,h}(t)\;dt.\\
A_{x,i}^{(4)} & = \frac{1}{2}d_{x,i}\langle R^{(0,2)}(\cdot,\sigma_{x,i}^+),f_{x,h}-g_n\rangle\int_{t_{x,i}}^{t_{x,{i+1}}}(t-t_{x,i})\varphi_{x,h}(t)\;dt.\\
A_{x,i}^{(5)} & = \frac{1}{2}\int_{t_{x,i}}^{t_{x,{i+1}}}(t-t_{x,i})^2\langle R^{(0,2)}(\cdot,\sigma_{x,t}^+),f_{x,h}-g_n\rangle\varphi_{x,h}(t)\;dt.
\end{align*}
We shall now control these quantities. Let,
\begin{align*}
B_{x,i}^{(1)} & = \underset{t_{x,i}< s,t< t_{x,{i+1}}}{\sup} |\varphi_{x,h}(t)-\varphi_{x,h}(s)|~\text{and}~
B_{x,i}^{(2)} = \underset{t_{x,i}< s,t< t_{x,{i+1}}}{\sup} |\alpha(t)\varphi_{x,h}(t)-\alpha(s)\varphi_{x,h}(s)|.
\end{align*}
Since $\alpha$ and $\varphi_{x,h}$ are Lipschitz 
then,
\begin{equation}
\underset{0\leq i \leq n}{\sup}\;B_{x,i}^{(1)} = O\bigg(\frac{1}{h^2} \underset{0\leq j \leq n}{\sup}~d_{j,n}\bigg) ~~~\text{and}~~ \underset{0\leq i \leq n}{\sup}\;B_{x,i}^{(2)} = O\bigg(\frac{1}{h^2} \underset{0\leq j \leq n}{\sup}~d_{j,n}\bigg).\label{B}
\end{equation}
\noindent Elementary calculations show that,
\begin{equation}
|A_{x,i}^{(1)}| \leq \frac{a_1}{h}B_{x,i}^{(1)} d_{x,i}^3 ,\;\;\;\;\;\;|A_{x,i}^{(2)}| \leq \frac{a_2}{h}B_{x,i}^{(1)} d_{x,i}^3 \;\;\;\;\text{and}\;\; |A_{x,i}^{(3)}| \leq \frac{a_3}{h}B_{x,i}^{(2)} d_{x,i}^3,
\label{majorationA}
\end{equation}
for appropriate constants $a_1,a_2$ and $a_3$. We obtain from the Cauchy-Schwartz inequality, Assumption $(C)$ and Proposition \ref{propoflimit} that,
\begin{align}
|A_{x,i}^{(4)}|+|A_{x,i}^{(5)}| & \leq \frac{a_4}{h} d_{x,i}^3||f_{x,h}-g_n|| \leq  \frac{1 }{h}d_{x,i}^3 \underbrace{a_4 \sqrt{\frac{C}{h}}}_{a_h} \underset{0 \leq j \leq n}{\sup}d_{j,n},\label{A45}
\end{align}
for an appropriate constant $a_4$ ($C$ is defined in Proposition \ref{propoflimit}).
\noindent Thus,
\begin{align}
&~~~~~~~~~~~~~~~~~\int_{t_{x,i}}^{t_{x,{i+1}}}(f_{x,h}(t)-g_n(t))\varphi_{x,h}(t)\;dt \nonumber\\
& = \frac{1}{12}d_{x,i}^3\alpha(\sigma_{x,i}) \varphi_{x,h}^2(\sigma_{x,i})+A_{x,i}^{(1)}-A_{x,i}^{(2)}-A_{x,i}^{(3)}-A_{x,i}^{(4)}+A_{x,i}^{(5)} \nonumber \\[0.2cm]
&\geq \frac{1}{12}d_{x,i}^3\alpha(\sigma_{x,i}) \varphi_{x,h}^2(\sigma_{x,i})
-d_{x,i}^3(\frac{a_1}{h}B_{x,i}^{(1)}+\frac{a_2}{h}B_{x,i}^{(2)}+\frac{a_h}{h}\underset{0 \leq j \leq n}{\sup}d_{j,n}). \label{ineq}
\end{align}
Let,
\[\rho_{h,N_{T_n}}=\underset{0 \leq i \leq N_{T_n}}{\sup}(\frac{a_1}{h}B_{x,i}^{(1)}+\frac{a_2}{h}B_{x,i}^{(2)}+\frac{a_h}{h} \underset{0 \leq j \leq n}{\sup} d_{j,n}).\]
\noindent Equation \eqref{B} implies that for an appropriate constant $c$ {and $c'$} we have,
 \[
  |\rho_{h,N_{T_n}}| \leq \Big(\frac{c}{h^3}\;\underset{0 \leq j \leq n}{\sup}d_{j,n}+{\frac{c'}{h^{3/2}}\;\underset{0 \leq j \leq n}{\sup}d_{j,n}}\Big).\]
Using \eqref{ineq} and \eqref{norm'}  together with Equation \eqref{ineq} in \eqref{norm'} we obtain,
\begin{align}
&||f_{x,h}-g_n||^2  \geq \sum_{i=1}^{N_{T_n}-1} \Big(\frac{1}{12}\alpha(\sigma_{x,i}) \varphi_{x,h}^2(\sigma_{x,i})-\rho_{h,N_{T_n}}\Big) d_{x,i}^3\nonumber \\
&  \geq  \frac{1}{12} \sum_{i=1}^{N_{T_n}-1}\alpha(\sigma_{x,i}) \varphi_{x,h}^2(\sigma_{x,i})d_{x,i}^3 - \frac{cN_{T_n}}{h^3}\;\underset{0 \leq j \leq n}{\sup}d_{j,n}^{\;4}-{\frac{c'N_{T_n}}{h^{3/2}}\;\underset{0 \leq j \leq n}{\sup}d_{j,n}^4}. \label{beforeholder}
\end{align}
Then the H{\"o}lder's inequality gives,
{\small
\begin{align*}
||f_{x,h}-g_n||^2 &\geq  \frac{1}{12(N_{T_n}-1)^2}\bigg\{\sum_{j=1}^{N_{T_n}-1} [\alpha(\sigma_{x,i}) \varphi_{x,h}^2(\sigma_{x,i})]^{\frac{1}{3}} d_{x,i}\bigg\}^3- \frac{cN_{T_n}}{h^3}\;\underset{0 \leq j \leq n}{\sup}d_{j,n}^{\;4}-{\frac{c'N_{T_n}}{h^{3/2}}\;\underset{0 \leq j \leq n}{\sup}d_{j,n}^4}.
\end{align*}}
We shall now control the first term of the right side of this inequality. We have,
\begin{align*}
&~~~~~~~~~~~~~~~~~~~~~~~~~~~~~\bigg\{ \sum_{j=1}^{N_{T_n}-1} \Big(
\alpha(\sigma_{x,i}) \varphi_{x,h}^2(\sigma_{x,i})\Big)^{\frac{1}{3}} d_{x,i}\bigg\}^3\\
& = \bigg\{\int_{x-h}^{x+h}\Big(\alpha(t) \varphi_{x,h}^2(t)\Big)^{\frac{1}{3}}\;dt-\sum_{j=1}^{N_{T_n}-1}\int_{t_{x,i}}^{t_{x,i+1}} \Big((\alpha(t) \varphi_{x,h}^2(t))^{\frac{1}{3}}-(\alpha(\sigma_{x,i}) \varphi_{x,h}^2(\sigma_{x,i}))^{\frac{1}{3}}\Big)\bigg\}^3 \nonumber\\
& = \bigg\{\int_{x-h}^{x+h} \Big(\alpha(t) \varphi_{x,h}^2(t)\Big)^{\frac{1}{3}}\;dt\bigg\}^3 - \bigg\{ \sum_{j=1}^{N_{T_n}-1}\int_{t_{x,i}}^{t_{x,i+1}} \Big((\alpha(t) \varphi_{x,h}^2(t))^{\frac{1}{3}}-(\alpha(\sigma_{x,i}) \varphi_{x,h}^2(\sigma_{x,i}))^{\frac{1}{3}}\Big)\bigg\}^3 \\
&~~~~-3\bigg\{\int_{x-h}^{x+h} \Big(\alpha(t) \varphi_{x,h}^2(t)\Big)^{\frac{1}{3}}\;dt\bigg\}^2 \bigg\{ \sum_{j=1}^{N_{T_n}-1}\int_{t_{x,i}}^{t_{x,i+1}} \Big((\alpha(t) \varphi_{x,h}^2(t))^{\frac{1}{3}}-(\alpha(\sigma_{x,i}) \varphi_{x,h}^2(\sigma_{x,i}))^{\frac{1}{3}}\Big)\bigg\}\\
&~~~~+3\bigg\{\int_{x-h}^{x+h} \Big(\alpha(t) \varphi_{x,h}^2(t)\Big)^{\frac{1}{3}}\;dt\bigg\} \bigg\{ \sum_{j=1}^{N_{T_n}-1}\int_{t_{x,i}}^{t_{x,i+1}} \Big((\alpha(t) \varphi_{x,h}^2(t))^{\frac{1}{3}}-(\alpha(\sigma_{x,i}) \varphi_{x,h}^2(\sigma_{x,i}))^{\frac{1}{3}}\Big)\bigg\}^2\\
& \overset{\Delta}{=}\bigg\{\int_{x-h}^{x+h} \Big(\alpha(t) \varphi_{x,h}^2(t)\Big)^{\frac{1}{3}}\;dt\bigg\}^3 + B,
\end{align*}
We obtain using \eqref{phiboundries} and the fact that $\alpha$ is Lipschitz,
\begin{align*}
B &= O\bigg(\Big(\frac{N_ {T_n}}{h^{5/3}}\;\underset{0 \leq j \leq n}{\sup}d_{j,n}^2\Big)^{3}\bigg)+O\bigg(\Big(\frac{N_ {T_n}}{h^{5/3}}\;\underset{0 \leq j \leq n}{\sup}d_{j,n}^2\Big)\;h^{\;2/3}\bigg)+O\bigg(\Big(\frac{N_ {T_n}}{h^{5/3}}\;\underset{0 \leq j \leq n}{\sup}d_{j,n}^2\Big)^{\;2}\;h^{1/3}\bigg).
\end{align*}
Assumption $(E)$ implies that for an appropriate constant $c''$ we have,
\begin{align*}
|B| \leq \frac{c''N_ {T_n}}{h}\;\underset{0 \leq j \leq n}{\sup}d_{j,n}^2.
\end{align*}Using the Riemann integrability of $\alpha $ and $\varphi_{x,h}$ we get,
\begin{align*}
&||f_{x,h}-g_n||^2  \geq \frac{1}{12(N_{T_n}-1)^2}\bigg\{\int_{x-h}^{x+h} \Big(\alpha(t) \varphi_{x,h}^2(t)\Big)^{\frac{1}{3}}\;dt\bigg\}^3-\frac{c''}{N_ {T_n}h}\;\underset{0 \leq j \leq n}{\sup}d_{j,n}^2\\
&~~~- \frac{cN_{T_n}}{h^3}\;\underset{0 \leq j \leq n}{\sup}d_{j,n}^{\;4}-{\frac{c'N_{T_n}}{h^{3/2}}\;\underset{0 \leq j \leq n}{\sup}d_{j,n}^4} \\
& \geq \frac{1}{12N_{T_n}^2}\bigg\{\int_{x-h}^{x+h} \Big(\alpha(t) \varphi_{x,h}^2(t)\Big)^{\frac{1}{3}}\;dt\bigg\}^3-\frac{c''}{N_ {T_n}h}\;\underset{0 \leq j \leq n}{\sup}d_{j,n}^2- \frac{cN_{T_n}}{h^3}\;\underset{0 \leq j \leq n}{\sup}d_{j,n}^{\;4}-{\frac{c'N_{T_n}}{h^{3/2}}\;\underset{0 \leq j \leq n}{\sup}d_{j,n}^4} \\
& = \frac{1}{12h^2N_{T_n}^2}\bigg\{\int_{x-h}^{x+h} \Big(\alpha(t) K^2(\frac{x-t}{h})\Big)^{\frac{1}{3}}\;dt\bigg\}^3-\frac{c''}{N_ {T_n}h}\;\underset{0 \leq j \leq n}{\sup}d_{j,n}^2\\
&~~~- \frac{cN_{T_n}}{h^3}\;\underset{0 \leq j \leq n}{\sup}d_{j,n}^{\;4}-{\frac{c'N_{T_n}}{h^{3/2}}\;\underset{0 \leq j \leq n}{\sup}d_{j,n}^4} \\
& = \frac{h}{12N_{T_n}^2}\bigg\{\int_{-1}^{1} \Big(\alpha(x-th) K^2(t)\Big)^{\frac{1}{3}}\;dt\bigg\}^3-\frac{c''}{N_ {T_n}h}\;\underset{0 \leq j \leq n}{\sup}d_{j,n}^2\\
&~~~- \frac{cN_{T_n}}{h^3}\;\underset{0 \leq j \leq n}{\sup}d_{j,n}^{\;4}-{\frac{c'N_{T_n}}{h^{3/2}}\;\underset{0 \leq j \leq n}{\sup}d_{j,n}^4}.
\end{align*}
Assumption $(E)$ implies that, \[\underset{ n \to \infty}{\lim}\frac{1}{h²}N_ {T_n}\;\underset{0 \leq j \leq n}{\sup}d_{j,n}^2=0,~~~\underset{ n \to \infty}{\lim} \frac{1}{h^4}\;\underset{0 \leq j \leq n}{\sup}d_{j,n}^{\;4}N_{T_n}^3= 0~~~\text{and}~~\underset{ n \to \infty}{\lim} {\frac{c'}{h^{3/2}}\;\underset{0 \leq j \leq n}{\sup}d_{j,n}^4N_{T_n}}=0.\] Finally  the continuity of $\alpha$ yields,
\[\underset{n \to \infty}{ \underline \lim }~\frac{N_{T_n}^2}{h}||f_{x,h}-g_n||^2 \geq \frac{1}{12}\alpha(x) \bigg\{\int_{-1}^{1} K(t)^{\frac{2}{3}}\;dt\bigg\}^3.\]
Inequality \eqref{liminf} is then proved for a sequence of designs containing $x-h$ and $x+h$.  Consider now any sequence of designs $ \{T_n, n\geq 1\}$ satisfying Assumption $(E)$ we can adjoin the points $\{x-h,x+h\}$ to $T_n$ (if they aren't present). Hence we form a sequence $\{S_n, n \geq 1\}$ with $S_n \in D_{n+2}$ and satisfying $\eqref{liminf}$. We have,
\[||f_{x,h} - P_{|S_n}f_{x,h}||^2 \leq ||f_{x,h} - P_{|T_n}f_{x,h}||^2.\]
Then,
\begin{equation}
N_{S_n}^2||f_{x,h} - P_{|S_n}f_{x,h}||^2 \leq N_{S_n}^2||f_{x,h} - P_{|T_n}f_{x,h}||^2 . \label{useitnext}
\end{equation}
We know that $N_{S_n} \in \{N_{T_n}+1, N_{T_n}+2\}$, replacing $N_{S_n}$ in the right term  of \eqref{useitnext} by $(N_{T_n}+2)$ (or ($N_{T_n}+1$) ) gives,
\[ \frac{N_{S_n}^2}{h}||f_{x,h} - P_{|S_n}f_{x,h}||^2-\frac{(4+2N_{T_n})}{h}||f_{x,h} - P_{|T_n}f_{x,h}||^2  \leq \frac{N_{T_n}^2}{h}||f_{x,h} - P_{|T_n}f_{x,h}||^2. \]
Assumption $(E)$ and Equation \eqref{lamajoration} yield,
\[ \underset{n \to \infty}{ \underline \lim }~\frac{(4+2N_{T_n})}{h}||f_{x,h} - P_{|T_n}f_{x,h}||^2 =0.\]
Hence, for any sequence $\{T_n , n\geq 1\}$ we have,
\[\underset{n \to \infty}{ \underline \lim }\frac{N_{T_n}^2}{h} ||f_{x,h} - P_{|T_n}f_{x,h}||^2 \geq \frac{1}{12} \alpha(x) \bigg\{ \int_{-1}^{1} K^{2/3}(t)dt \bigg\}^3.\]
This completes the proof of Proposition \ref{variancetheorem}. $\Box$
\subsection{Proof of Proposition \ref{optimalrateproposition}.}
On the one hand, Proposition \ref{propoflimit} yields that there exists a constant $c > 0$ such that,
\begin{equation*}
0 \leq \frac{\sigma_{x,h}^2}{m}-\Var~\hat{g}^{pro}_n (x)  \leq \frac{c}{mh}~\underset{0\leq j\leq n}{\sup}d_{j,n}^2.
\end{equation*}
Lemma \ref{lemmaregulardesign} implies that there exists a  constant $c' > 0$ such that, \[\underset{0\leq j\leq n}{\sup}d_{j,n}^2 \leq \frac{c'}{n^2}.\]
Thus, for $n \geq 1$ we have,
\begin{equation*}
0 \leq \frac{\sigma_{x,h}^2}{m}-\Var~\hat{g}^{pro}_n (x)  \leq \frac{c'c}{mn^2h}.
\end{equation*}
Finally, taking $C=cc'$ we obtain,
\begin{equation*}
\underset{n \to \infty}{\overline \lim }~ mn^2h\Big( \frac{\sigma_{x,h}^2}{m}-\Var~\hat{g}^{pro}_n (x) \Big) \leq C.
\end{equation*}
Inequality \eqref{firstineq} is then proved. On the  other hand, Proposition \ref{variancetheorem} yields,
\begin{equation*}
\frac{mN_{T_n}^2}{h} \Big( \frac{\sigma_{x,h}^2}{m}-\Var~\hat{g}^{pro}_n (x) \Big) \geq \frac{1}{12} \alpha(x) \bigg\{ \int_{-1}^{1} K^{2/3}(t)dt \bigg\}^3.
\end{equation*}
Lemma \ref{lemmaregulardesign} implies that there exists a constant $c'' > 0$ such that,  \[N_{T_n} < c'' nh,\]
which implies that,
\begin{equation*}
c''mn^2h~\Big( \frac{\sigma_{x,h}^2}{m}-\Var~\hat{g}^{pro}_n (x) \Big) \geq \frac{1}{12} \alpha(x) \bigg\{ \int_{-1}^{1} K^{2/3}(t)dt \bigg\}^3.
\end{equation*}
Finally, taking $C' = \frac{1}{12c''} \alpha(x) \bigg\{ \int_{-1}^{1} K^{2/3}(t)dt \bigg\}^3$ we obtain,
\[\underset{n \to \infty}{\underline \lim }~ mn^2h\Big( \frac{\sigma_{x,h}^2}{m}-\Var~\hat{g}^{pro}_n (x) \Big) \geq C'.\]
This concludes the proof of Proposition \ref{optimalrateproposition}.~$\Box$
\subsection{Proof of Proposition \ref{exactvarianceformula}}
The first part of this proof is the same as that of Proposition \eqref{variancetheorem}. Recall that, \[m\Big( \frac{\sigma_{x,h}^2}{m}-\Var~\hat{g}^{pro}_n (x) \Big)=||f_{x,h}||^2 -||P_{|T_n}f_{x,h}||^2=||f_{x,h} - P_{|T_n}f_{x,h}||^2.\]
Using \eqref{norm'} and \eqref{touse} we obtain,
\begin{align}
&\Var~\hat{g}^{pro}_n (x)- \frac{\sigma_{x,h}^2}{m}  = -\frac{1}{m} ||f_{x,h}-P_{|T_n}f_{x,h}||^2 \nonumber \\
& = -\frac{1}{m} \sum_{i=1}^{N_{T_n}} \bigg( \frac{1}{12}d_{x,i}^3\alpha(\sigma_{x,i}) \varphi_{x,h}^2(\sigma_{x,i})+A_{x,i}^{(1)}-A_{x,i}^{(2)}-A_{x,i}^{(3)}-A_{x,i}^{(4)}+A_{x,i}^{(5)}\bigg),\label{ladifference}
\end{align}
for some $\sigma_{x,i} \in ]t_{x,i},t_{x,{i+1}}[$  and some $\sigma_{x,t} \in ]t_{x,i},t[$, where,
\begin{align*}
A_{x,i}^{(1)}&=\frac{1}{2}d_{x,i}\alpha(\sigma_{x,i}) \varphi_{x,h}(\sigma_{x,i})\int_{t_{x,i}}^{t_{x,{i+1}}}(t-t_{x,i})[\varphi_{x,h}(t)-\varphi_{x,h}(\sigma_{x,i})]\;dt.\\
A_{x,i}^{(2)} & =\frac{1}{2}\alpha(\sigma_{x,i}) \varphi_{x,h}(\sigma_{x,i})\int_{t_{x,i}}^{t_{x,{i+1}}}(t-t_{x,i})^2[\varphi_{x,h}(t)-\varphi_{x,h}(\sigma_{x,i})]\;dt.\\
A_{x,i}^{(3)} & = \frac{1}{2}\int_{t_{x,i}}^{t_{x,{i+1}}}(t-t_{x,i})^2[\alpha(\sigma_{x,t})\varphi_{x,h}(\sigma_{x,t})-\alpha(\sigma_{x,i})\varphi_{x,h}(\sigma_{x,i})]\varphi_{x,h}(t)\;dt.\\
A_{x,i}^{(4)} & = \frac{1}{2}d_{x,i}\langle R^{(0,2)}(\cdot,\sigma_{x,i}^+),f_{x,h}-g_n\rangle\int_{t_{x,i}}^{t_{x,{i+1}}}(t-t_{x,i})\varphi_{x,h}(t)\;dt.\\
A_{x,i}^{(5)} & = \frac{1}{2}\int_{t_{x,i}}^{t_{x,{i+1}}}(t-t_{x,i})^2\langle R^{(0,2)}(\cdot,\sigma_{x,t}^+),f_{x,h}-g_n\rangle\varphi_{x,h}(t)\;dt.
\end{align*}
From the definition of the regular sequence of  designs (see Definition \ref{regulardesignexample})  and the mean value theorem we have for $i=1,\cdots,N_{T_n}$,
\begin{align*}
d_{x,i} &= t_{x,i+1}-t_{x,i} = F^{-1}\Big(\frac{i+1}{n}\Big)-F^{-1}\Big(\frac{i}{n}\Big) = \frac{1}{nf(t_{x,i}^*)},
\end{align*}
where $t_{x,i}^* \in ]t_{x,i},t_{x,i+1}[$. Using this together with \eqref{ladifference} we obtain,
\begin{align*}
\Var~\hat{g}^{pro}_n (x) - \frac{\sigma_{x,h}^2}{m}
& = -\frac{1}{12mn^2} \sum_{i=1}^{N_{T_n}} d_{x,i} \frac{1}{f^2(t_{x,i}^*)}\alpha(\sigma_{x,i}) \varphi_{x,h}^2(\sigma_{x,i})\\
&~~~~~~-\frac{1}{m}\sum_{i=1}^{N_{T_n}}\bigg( A_{x,i}^{(1)}-A_{x,i}^{(2)}-A_{x,i}^{(3)}-A_{x,i}^{(4)}+A_{x,i}^{(5)}\bigg).
\end{align*}
Lemma \ref{lemmaregulardesign} yields that $N_{T_n} = O(nh)$.
Using \eqref{majorationA}, \eqref{A45} and \eqref{B} we obtain, \[A_{x,i}^{(1)} = O\Big(\frac{1}{n^4h^3}\Big),~~A_{x,i}^{(2)} = O\Big(\frac{1}{n^4h^3}\Big),~~A_{x,i}^{(3)} = O\Big(\frac{1}{n^4h^3}\Big)~~\text{and}~~A_{x,i}^{(4)}+A_{x,i}^{(5)} = O\Big(\frac{1}{n^4h^{3/2}}\Big). \]Finally,
\begin{align*}
\Var~\hat{g}^{pro}_n (x)- \frac{\sigma_{x,h}^2}{m}
& = -\frac{1}{12mn^2} \sum_{i=1}^{N_{T_n}} d_{x,i} \frac{1}{f^2(t_{x,i}^*)}\alpha(\sigma_{x,i}) \varphi_{x,h}^2(\sigma_{x,i})+O\bigg(\frac{1}{mn^3h^2}+\frac{1}{mn^3 \sqrt{h}}\bigg).
\end{align*}
Using a classical approximation of a sum by an integral (see for instance, Lemma 2 in  \cite{Benelmadani and Benhenni and Louhichi}) and the fact that $0 < h < 1$ we obtain,
\begin{align*}
\Var~\hat{g}^{pro}_n (x)- \frac{\sigma_{x,h}^2}{m}
& = -\frac{1}{12mn^2}\int_{x-h}^{x+h} \frac{\alpha(t) }{f^2(t)}\varphi_{x,h}^2(t)\;dt+O\bigg(\frac{1}{mn^3h^2}\bigg).
\end{align*}
This concludes the proof of Proposition \ref{exactvarianceformula}. $\Box$
\subsection{Proof of Theorem \ref{thmise}.}
First, note that since $\alpha$ and $f$ are Lipschitz functions then the asymptotic expression of the integral in \eqref{varexpress} is:
{\small
	\begin{align*}
	\frac{1}{mn^2}\int_{x-h}^{x+h}\frac{\alpha(t)}{f^2(t)} \varphi_{x,h}^{2}(t)dt  &= \frac{1}{mn^2h}\int_{-1}^1 \frac{\alpha(x-th)}{f^2(x-th)} K^2(t)\;dt\\
	& = \frac{1}{mn^2h}\Big(\frac{\alpha(x)}{f^2(x)}\int_{-1}^1  K^2(t)\;dt+\int_{-1}^1 \Big(\frac{\alpha(x-th)}{f^2(x-th)} -\frac{\alpha(x)}{f^2(x)}\Big) K^2(t)\;dt\Big)\\
	& = \frac{1}{mn^2h}\frac{\alpha(x)}{f^2(x)}\int_{-1}^1  K^2(t)\;dt + O\Big(\frac{1}{mn^2}\Big).
	\end{align*}}
This last equality together with Proposition \ref{exactvarianceformula} and Proposition \ref{biasofwinner} concludes the proof of  Theorem \ref{thmise}. $\Box$
\subsection{Proof of Corollary \ref{hoptcor}.}
Let $I_1=\int_{0}^{1}R(x,x)w(x)\;dx$ and put, \[\Psi (h,m)=-\frac{C_K h}{2m} \int_{0}^{1}\alpha(x)  w(x)\;dx + \frac{1}{4} h^4 B^2 \int_{0}^{1}[g''(x)]^2 w(x)\;dx.\]
We have from Theorem \ref{MISEGEN},
\begin{align*}
\IMSE(h)
& = \frac{I}{m}+\Psi (h,m) +o\Big(h^4+\frac{h}{m}\Big)+O\Big(\frac{1}{mn^2h}+\frac{h}{n}+\frac{1}{n^2h^2}\Big),
\end{align*}
Let $h^*$ be as defined by \eqref{hoptimal}. It is clear that $h^*=\underset{0 < h <1}{\argmin}~\Psi (h,m)$ so that $\Psi (h,m) \geq \Psi (h^*,m)$  for every $0< h < 1$. Let $h_{n,m}$ be as defined in Corollary \ref{hoptcor}. We have,
\begin{align*}
&\frac{\IMSE(h^*)}{\IMSE(h_{n,m})}
= \frac{ \frac{I_1}{m}  +\Psi(h^*,m)  +o \Big({h^*}^4+\frac{h^*}{m}\Big)+O\Big(\frac{1}{mn^2h^*}+\frac{h^*}{n}+\frac{1}{n^2{h^*}^2}\Big)}{\frac{I_1}{m}  +\Psi(h_{n,m},m) +o \Big(h_{n,m}^4+\frac{h_{n,m}}{m}\Big)+O\Big(\frac{1}{mn^2h_{n,m}}+\frac{h_{n,m}}{n}+\frac{1}{n^2h_{n,m}^2}\Big)}\\
&\\
&\leq \frac{ I_1  +m\Psi(h_{n,m},m)  +o \Big(m{h^*}^4+h^*\Big)+O\Big(\frac{1}{n^2h^*}+\frac{mh^*}{n}+\frac{m}{n^2h^2}\Big)}{I_1  +m\Psi(h_{n,m},m) +o \Big(mh_{n,m}^4+h_{n,m}\Big)+O\Big(\frac{1}{n^2h_{n,m}}+\frac{mh_{n,m}}{n}+\frac{m}{n^2h_{n,m}^2}\Big)}.\\
\end{align*}
We have, using the definition of $h^*,$ $mh_{n,m}^3=O(1)$,  $\underset{n,m \to \infty}{\lim}~h_{n,m}=0$  and using the assumption $\frac{m}{n}=O(1)$ as $n,m \to \infty$ we know that $m\Psi(h_{n,m},m)=O(h_{n,m})$. Thus,  \[ \underset{n,m \to \infty}{\overline {\lim}} \frac{\IMSE(h^*)}{\IMSE(h_{n,m})} \leq 1.\]This concludes the proof of Corollary \ref{hoptcor}. $\Box$
\subsection{Proof of Theorem \ref{asymptotic normality theorem}.}
Let $x \in ]0,1[$ be fixed. We have the following decomposition,
{\small
	\begin{equation} \label{diffasym}
	\sqrt{m} \Big(\hat{g}_{n,m}^{pro}(x) -g(x) \Big)=\sqrt{m} \Big(\hat{g}_{n,m}^{pro}(x) -\mathbb{E}\big( \hat{g}_{n,m}^{pro}(x)\big) \Big)+\sqrt{m}~\Big( \mathbb E (\hat{g}_{n,m}^{pro}(x))-g(x)\Big).
	\end{equation}}
Since $\underset{n, m \to \infty}{\lim}\sqrt{m}h =0$, $\frac{n}{m} = O(1)$ as $n,m \to \infty$ and $\underset{n,m \to \infty}{\lim} nh^2 =\infty$ then Remark \ref{biasforregulardesign} implies that,
\begin{equation}
\underset{n,m \to \infty}{\lim}\sqrt{m}~\Big( \mathbb E (\hat{g}_{n,m}^{pro}(x))-g(x)\Big) =0.
\end{equation}
Consider now the first term of the right side of \eqref{diffasym}. Since $\overline{Y}(t_{x,i})-\mathbb{E}(\overline{Y}(t_{x,i})) =\overline{\varepsilon}(t_{x,i})$, we have, as done by Fraiman and Iribarren (1991) \cite{Fraiman},
\begin{align}
& \sqrt{m} \Big(\hat{g}_{n,m}^{pro}(x) -\mathbb{E}\big( \hat{g}_{n,m}^{pro}(x)\big) \Big)  =  \frac{1}{\sqrt{m}}\sum_{j=1}^{m}  \sum_{i=1}^{n}  m_{x,h}(t_i) \varepsilon_j(t_i)\nonumber \\
& =  \frac{1}{\sqrt{m}}\sum_{j=1}^{m}  \sum_{i=1}^{n}  m_{x,h}(t_i) \big(\varepsilon_j(t_i)- \varepsilon_j(x))  + \bigg(\sum_{i=1}^{n}  m_{x,h}(t_i) \bigg)  \bigg(\frac{1}{\sqrt{m}} \sum_{j=1}^{m} \varepsilon_j(x)\bigg). \label{equationinasym}
\end{align}
\noindent We start by controlling the second term of this last equation. Using Lemma  \ref{mxhlemma} together with Lemma \ref{lemmaregulardesign} we obtain,
\begin{align*}
&~~~~~~~~~~~~~~~~~~~~~~~~~~~~~~~~~~~~ m_{x,h}(t_{i,n}) = \\
& \begin{cases}
\frac{1}{2} \varphi_{x,h} (t_{i,n}) (t_{i+1,n}-t_{i-1,n})+O\Big(\frac{1}{n^2h^2}+\frac{1}{n^2 \sqrt h}\Big)~~~~~if~~ i \notin \{1,n\} \text{ and }\\
~~~~~~~~~~~~~~~~~~~~~~~~~~~~~~~~~~~~~~~~~~~~~~~~~~ [t_{i-1,n},t_{i+1,n}]\cap [x-h,x+h] \ne \emptyset ,\\[0.1cm]
O\Big(\frac{1}{n^2h^2}+\frac{1}{n^2 \sqrt h}\Big)~~~~~~~~~~~~~~~~~~~~~~~~~~~~~~~~~~~~~~\text{if}~~ i \in \{1,n\},\\
O\Big(\frac{1}{n^2 \sqrt h}\Big)~~~~~~~~~~~~~~~~~~~~~~~~~~~~~~~~~~~~~~~~~~~~~~~otherwise.
\end{cases}
\end{align*}
Recall that $N_{T_n}=\Card~ I_{x,h} = \Card~\{i=1,\cdots,n~:~[t_{i-1},t_{i+1}]\cap ]x-h,x+h[ \ne \emptyset \}$ and denote by $t_{x,i}$ the points of $T_n$ for which $i \in I_{x,h}$, Lemma \ref{lemmaregulardesign} yields that $N_{T_n}=O(nh)$. Thus,
\begin{align*}
\sum_{i=1}^{n}  m_{x,h}(t_i) & = \frac{1}{2} \sum_{i=2}^{N_{T_n}-1} \varphi_{x,h} (t_{x,i}) (t_{x,i+1}-t_{x,i-1}) +O\Big(\frac{1}{nh}\Big).
\end{align*}
Since $\underset{n \to \infty}{\lim} nh = +\infty$, then using the Riemann integrability of $K$, we obtain,
\begin{align*}
\underset{n,m \to \infty}{\lim}\sum_{i=1}^{n}  m_{x,h}(t_i) =\frac{1}{2}\underset{n,m \to \infty}{\lim} \sum_{i=2}^{N_{T_n}-1} \varphi_{x,h} (t_{x,i}) (t_{x,i+1}-t_{x,i-1}) = \int_{-1}^{1} K(t)~dt = 1.
\end{align*}
The Central Limit Theorem for i.i.d. variables yields,
\[ \frac{1}{\sqrt{m}} \sum_{j=1}^{m} \varepsilon_j(x) \underset{m \to \infty}{\overset{\mathscr D}{\longrightarrow}} Z ~~~where~~ Z \sim \mathcal{N}(0,R(x,x)). \]
We shall prove now that the first term of Equation \eqref{equationinasym} tends to 0 in probability as $n,m$ tends to infinity.  Let,
\begin{align*}
A_{m,n}(x) & = \frac{1}{\sqrt{m}} \sum_{j=1}^{m} \sum_{i=1}^{n} m_{x,h}(t_i) \big (\varepsilon_j (t_i )-\varepsilon_j(x) \big)\overset{\Delta}{=}\frac{1}{\sqrt{m}} \sum_{j=1}^{m} T_{n,j}(x).
\end{align*}
From the Chebyshev inequality, it suffices to prove that $\underset{n,m \to \infty}{\lim} \mathbb{E}(A_{m,n}^2(x)) =0 $. We have for $j \ne l$, $\mathbb{E} (\varepsilon_j(x)\varepsilon_l(y))=0$ then $\mathbb{E} (T_{n,j}(x) T_{n,l}(x)) =0$. Hence,
\begin{align*}
\mathbb{E}(A_{m,n}^2(x)) & = \frac{1}{m} \sum_{j=1}^{m} \sum_{l=1}^{m} \mathbb{E} (T_{n,j}(x) T_{n,l}(x)) = \frac{1}{m} \sum_{j=1}^{m} \mathbb{E} (T_{n,j}^2(x)).
\end{align*}
We have,
\begin{align*}
\mathbb{E} (T_{n,j}^2(x)) &= \sum_{i=1}^{n} \sum_{k=1}^{n} m_{x,h}(t_i) m_{x,h}(t_k) \mathbb{E} \Big(  \big (\varepsilon_j (t_i )-\varepsilon_j(x) \big) \big (\varepsilon_j (t_k )-\varepsilon_j(x) \big) \Big)\\
& = \sum_{i=1}^{n} \sum_{k=1}^{n} m_{x,h}(t_i) m_{x,h}(t_k) \Big(R(t_i,t_k)-R(t_i,x)-R(x,t_k)+R(x,x)\Big).
\end{align*}
Note that $\mathbb{E} (T_{n,j}^2(x))$ does not depend on $j$ hence,
\begin{align}
\mathbb{E}(A_{m,n}^2(x))& = \sum_{i=1}^{n} \sum_{k=1}^{n} m_{x,h}(t_i) m_{x,h}(t_k) \Big(R(t_i,t_k)-R(t_i,x)-R(x,t_k)+R(x,x)\Big)\nonumber\\
& \overset{\Delta}{=} B_{n,1}(x)-B_{n,2}(x)-B_{n,3}(x)+B_{n,4}(x). \label{sumBn}
\end{align}
Using  Lemma  \ref{mxhlemma} and the approximation of a sum by an integral (see, for instance, Lemma 2 in Benelmadani et al. (2018a) \cite{Benelmadani and Benhenni and Louhichi}) we obtain,
\begin{align*}
B_{n,1}(x) & = \int_{x-h}^{x+h} \int_{x-h}^{x+h} \varphi_{x,h}(s)\varphi_{x,h}(t)R(s,t)~ds~dt+ O\Big(\frac{1}{nh}\Big)=\sigma_{x,h}^2+ O\Big(\frac{1}{nh}\Big).
\end{align*}
Using Equation \eqref{devesigma}  we obtain, \[B_{n,1}(x)  = R(x,x)-\frac{1}{2}\alpha(x)C_K h  + o(h)+ O\Big(\frac{1}{nh}\Big).\]
where $C_K = \int_{-1}^1 \int_{-1}^1 |u-v|K(u)K(v) du dv.$ Since $\underset{n \to \infty}{\lim}h=0$ and $\underset{n \to \infty}{\lim}nh=\infty$ then,
\begin{equation}\label{Bn1}
\underset{n \to \infty}{\lim} B_{n,1}(x) = R(x,x).
\end{equation}
Consider now the term $B_{n,2}(x)$. We obtain using Lemma  \ref{mxhlemma} and the approximation of a sum by an integral,
\begin{align*}
B_{n,2}(x) & = \int_{x-h}^{x+h} \int_{x-h}^{x+h} \varphi_{x,h}(s)\varphi_{x,h}(t)R(s,x)~ds~dt+ O\Big(\frac{1}{nh}\Big)\\
&= \int_{x-h}^{x+h}  \varphi_{x,h}(s)R(s,x)~ds+ O\Big(\frac{1}{nh}\Big)\\
& = \int_{-1}^{1}  K(s)R(x-hs,x)~ds+ O\Big(\frac{1}{nh}\Big)\\
& = \int_{-1}^{0}  K(s)R(x-hs,x)~ds+\int_{0}^{1}  K(s)R(x-hs,x)~ds+  O\Big(\frac{1}{nh}\Big).
\end{align*}
For $s \in ]-1,0[$, Taylor expansion of $R(\cdot,x)$ around $x$ yields, \[R(s,x)  = R(x-sh,x)-shR^{(1,0)}(x^+,x)+o(h).\]
Similarly for  $s \in ]0,1[$ we obtain,\[R(x-sh,x)  = R(x,x)-shR^{(1,0)}(x^-,x)+o(h).\]
Thus,
{\small
\begin{align*}
B_{n,2}(x) & = R(x,x)-hR^{(1,0)}(x^+,x)\int_{-1}^{0}s~k(s)~ds-hR^{(1,0)}(x^-,x)\int_{0}^{1}s~k(s)~ds+o(h)+ O\Big(\frac{1}{nh}\Big).
\end{align*}}
Hence, \begin{equation} \label{Bn2}
\underset{n \to \infty}{\lim} B_{n,2}(x) = R(x,x).
\end{equation}
Similarly, \begin{equation}\label{Bn3}
\underset{n \to \infty}{\lim} B_{n,3}(x) = R(x,x).
\end{equation}
It is easy to verify that,\begin{equation}\label{Bn4}
\underset{n \to \infty}{\lim} B_{n,4}(x) = R(x,x).
\end{equation}
Inserting \eqref{Bn1}, \eqref{Bn2}, \eqref{Bn3} and \eqref{Bn4} in \eqref{sumBn} yields, \[\underset{n,m \to \infty}{\lim} \mathbb{E}(A_{m,n}^2(x)) =0. \]
This concludes the proof of Theorem \ref{asymptotic normality theorem}. $\Box$
\subsection{Proof of Theorem \ref{compavariancegen}.} Let $x \in ]0,1[$.
On the one hand, we have from Proposition \ref{exactvarianceformula} and Remark \ref{remark4},
\begin{equation}
\Var~\hat{g}^{pro}_n (x)= \frac{\sigma_{x,h}^2}{m}-\frac{A }{12mn^2h}\frac{\alpha(x)}{f^2(x)}+ O\Big(\frac{1}{mn^3h^2}+\frac{1}{mn^2}\Big),\label{varproproof}
\end{equation}
where $A= \int_{-1}^1  K^2(t)\;dt$.
On the other hand,  it can be seen in Benelmadani et al. (2019) \cite{Benelmadani and Benhenni and Louhichi and Su} that,
\begin{align}
\Var \;\hat{g}_n^{GM}(x) & = \frac{\sigma_{x,h}^2}{m} + O\Big(\frac{1}{mn^2}+\frac{1}{mn^3h^2}\Big).\label{varmuproof}
\end{align}
Equations \eqref{varproproof} and \eqref{varmuproof} then yield,
\begin{align*}
mn^2h ~\Big(\Var\;\hat{g}_n^{GM}-\Var\;\hat{g}_n^{pro} \Big) & = \frac{A}{12}\frac{\alpha(x)}{f^2(x)}+ O\Big(h+\frac{1}{nh} \Big).
\end{align*}
Recall that $\alpha(x)>0$ and that $\frac{1}{f(x)}>0$. Since  $h \to 0$ and $nh \to \infty$ as $n,m \to \infty$ we obtain, \[ \underset{n, m \to \infty}{ \lim}~mn^2h ~\Big(\Var\;\hat{g}_n^{GM}(x)-\Var\;\hat{g}_n^{pro}(x) \Big)= \frac{A}{12}\frac{\alpha(x)}{f^2(x)} > 0.\]
This concludes the proof of Theorem \ref{compavariancegen}. $\Box$

\subsection{Proof of Theorem \ref{compavariance}.}
We have from the proof of Proposition \ref{biasofwinner} (Equation \eqref{espewinner}) for any $x \in ]0,1[$,
\begin{align}
\mathbb E (\hat{g}_{n,m}^{pro}(x))-g(x) & = I_h(x) -g(x) + O\Big(\frac{1}{n^2h}\Big), \label{biasproproof}
\end{align}
where, \[I_h(x) = \int_{x-h}^{x+h}\varphi_{x,h}(s) g(s) ~ds.\]
Hence, using  \eqref{varproproof} and \eqref{biasproproof}  we get for a positive density measure  $w$,
\begin{align}
\IMSE (\hat{g}_n^{pro}) & = \frac{1}{m}\int_{0}^{1}\sigma_{x,h}^2~w(x)~dx-\frac{A }{12mn^2h}\int_{0}^{1}\frac{\alpha(x)}{f^2(x)}~w(x) ~dx+ \int_{0}^{1}\big(I_h(x)-g(x)\big)^2~w(x)~dx\nonumber\\
&~~~~~~+ O\Big(\frac{1}{n^4h^2}+\frac{h}{n^2}+\frac{1}{mn^3h^2}+\frac{1}{mn^2}\Big). \label{firstpart}
\end{align}
It can be seen in Benelmadani et al. (2019) \cite{Benelmadani and Benhenni and Louhichi and Su} that,
\begin{align}
\mathbb E (\hat{g}_{n,m}^{GM}(x))-g(x) & = I_h(x) -g(x) + O\Big(\frac{1}{n^2h}\Big). \label{biasmuproof}
\end{align}
Using \eqref{varmuproof} and \eqref{biasmuproof} yield,
\begin{align}
\IMSE (\hat{g}_n^{GM}) & = \frac{1}{m}\int_{0}^{1}\sigma_{x,h}^2~w(x)~dx+ \int_{0}^{1}\big(I_h(x)-g(x)\big)^2~w(x)~dx~\nonumber\\&~~~~ + O\Big(\frac{1}{n^4h^2}+\frac{h}{n^2}+\frac{1}{mn^2}+\frac{1}{mn^3h^2}\Big). \label{secondpart}
\end{align}
Then, Equations \eqref{firstpart} and \eqref{secondpart} yield,
\begin{align*}
mn^2h ~\Big(\IMSE\;(\hat{g}_n^{GM})-\IMSE\;(\hat{g}_n^{pro}) \Big) & = \frac{A}{12}\int_{0}^{1}\frac{\alpha(x)}{f^2(x)}~w(x) ~dx+ O\Big(\frac{m}{n^2h}+mh^2+h+\frac{1}{nh} \Big).
\end{align*}
Since $\frac{m}{n}=O(1)$ and $mh^2 \to 0$  as $n,m \to \infty$ we obtain, \[ \underset{n, m \to \infty}{\lim}~mn^2h ~\Big(\IMSE\;(\hat{g}_n^{GM})-\IMSE\;(\hat{g}_n^{pro}) \Big) = \frac{A}{12}\int_{0}^{1}\frac{\alpha(x)}{f^2(x)}~w(x) ~dx > 0.\]
This concludes the proof of Theorem \ref{compavariance}. $\Box$
\section{Appendix}\label{RKHSAppendix}
Let $\varepsilon=(\varepsilon(t))_{t\in [0,1]}$ be a centered and a second order process of autocovariance R, such that $R$ is invertible when restricted to any finite set on $[0,1].$ Let $L(\varepsilon(t), t \in [0,1])$ be the set of all random variables which maybe be written as a linear combinations of $\varepsilon(t)$ for $t \in [0,1]$, i.e., the set of random variables of the form $\sum_{i=1}^{l} \alpha_i \varepsilon(t_i)$ for some positive integer $l$ and some constants $\alpha_i$, $t_i\in[0,1]$ for $i=1,\cdots,l$.
 Let also $L_2(\varepsilon)$ be the Hilbert space of all square integrable random variables in the linear manifold $L(\varepsilon(t), t \in [0,1])$, together with all random variables $U$ that are limits in $\mathbb L^2$ of  a sequence of random variables $U_n$ in $L(\varepsilon(t), t \in [0,1])$, i.e,  $U$ is such that, \[\exists\, (U_n)_{n\geq 0}\in L(\varepsilon(t), t \in [0,1]):\,\,\, \underset{n \to \infty}{\lim}\mathbb E ( (U_n-U)^2 )=0 .\]
\\
 Denote by $\mathcal F (\varepsilon)$ the family of functions $g$ on $[0,1]$ defined by,\[\mathcal F (\varepsilon) = \{g : [0,1]\to \mathbb R~\text{with}~g(\cdot)=\mathbb E (U\varepsilon(\cdot))~\text{ where }~U \in L_2(\varepsilon) \},\]
We note here that for every $g \in \mathcal F (\varepsilon)$, the associated $U$ is unique. It is easy to verify that $\mathcal F (\varepsilon)$ is a Hilbert space equipped with  the norm $||~||$ defined for $g \in \mathcal F (\varepsilon)$ by,
\[ ||g||^2= \mathbb E (U^2).\] In fact,
let $g\in \mathcal F (\varepsilon)$, i.e, $g(\cdot) = \mathbb E(U\varepsilon(\cdot))$ for some $U\in L_2(\epsilon)$. We have,
\begin{itemize}
	\item $||g||  = \sqrt {\mathbb E(U^2)}  \geq 0.$
	\item $||g||  = \sqrt {\mathbb E(U^2)} = 0 \Rightarrow U=0~\text{a.s.}~\Rightarrow g=0.$
	\item For $g\in \mathcal F (\varepsilon)$, i.e, $f(\cdot) = \mathbb E(V\varepsilon(\cdot))$ some $V\in L_2(\epsilon)$. We have, \begin{align*}
	||g+f||^2  & =  \mathbb E( (U+V)^2) =  \mathbb E (U^2)+\mathbb E (V^2)+2\mathbb E (UV)\\
	& \leq \mathbb E (U^2)+\mathbb E (V^2)+2\sqrt {\mathbb E (U^2)}\sqrt {\mathbb E (V^2)} = \Big( \sqrt {\mathbb E (U^2)}+\sqrt {\mathbb E (V^2)} \Big)^2.
	\end{align*}
\noindent Thus,
	$||g+f|| \leq \sqrt {\mathbb E (U^2)}+\sqrt {\mathbb E (V^2)} = ||g||+||f||.$
\end{itemize}
We now prove the completeness of $\mathcal F (\varepsilon)$. For this let $g_n(\cdot)=E(U_n \varepsilon (\cdot))$ be a Cauchy sequence in $\mathcal F (\varepsilon)$, i.e., \[\underset{n,m \to \infty}{\lim}||g_n-g_m||^2=0.\]
From the definition of the norm $||~||$  we obtain,
 \[\underset{n,m \to \infty}{\lim} \mathbb E ( (U_n-U_m)^2 )=\underset{n,m \to \infty}{\lim}||g_n-g_m||^2  =0.\]
This yields that $(U_n)_{n \geq 1}$ is a Cauchy sequence in $L_2(\varepsilon)$, which is a Hilbert space as proven by  \cite{Parzen} (see page 8 there). Thus it exists $U \in L_2(\varepsilon)$ such that, \[\underset{n \to\infty}{\lim} \mathbb E ( (U_n-U)^2 ) =0.\]
Taking $g(\cdot)=\mathbb E ( U\varepsilon(\cdot)),$ which is clearly an element of $\mathcal F (\varepsilon)$ gives, \[\underset{n \to \infty}{\lim}||g_n-g||^2 = \underset{n \to \infty}{\lim} \mathbb E ( (U_n-U)^2 ) =0.\]This concludes the proof of completness of $\mathcal F (\varepsilon)$. \\ 

The Hilbert space $\mathcal F (\varepsilon)$ can easily be identified as the Reproducing Kernel Hilbert Space associated to a reproducing kernel $R$  (with $R(s,t)=\mathbb E(\varepsilon(s)\varepsilon(t))$),
which is defined as follows.
\begin{definition} Parzen (1959)  \cite{Parzen}
	A Hilbert space $H$ is said to be a Reproducing Kernel Hilbert Space associated to a reproducing kernel (or function) $R$ (RKHS($R$)), if its members are functions on some set $T$, and if there is a kernel $R$ on $T \times T$ having the following two properties:
	\begin{equation}
	\begin{cases}
	R(\cdot,t) \in H ~~~~~~~~~~~\text{ for all} ~~ t \in T, \\
	\langle g,R(\cdot,t)\rangle=g(t) ~~~~~ \text{for all} ~~t \in T~~ \text{and} ~~g \in H,
	\end{cases}
	\label{propertiesRKHS}
	\end{equation}
	where $\langle \cdot,\cdot \rangle$ is the inner (or scalar) product in H.
\end{definition}
To prove this, we need to verify the properties given in \eqref{propertiesRKHS}. For  $t \in [0,1]$ we have,
\[R(s,t)=\mathbb E(\varepsilon(s)\varepsilon(t))~~~\text{ for all}~~ s \in [0,1].\]
Since $\varepsilon(s) \in L_2 (\varepsilon)$ then $R(\cdot,t) \in \mathcal F (\varepsilon)$ for any fixed $t\in [0,1]$.
Now let $g \in \mathcal F (\varepsilon)$, i.e., \[g(\cdot)=\mathbb E(U\varepsilon(\cdot))~~\text{for some}~~ U \in L_2(\varepsilon). \]Then,
\begin{align*}
\langle g,R(\cdot,t)\rangle& = \frac{1}{2} \big( ||g||^2 + ||R(\cdot,t)||^2 - ||g-R(\cdot,t)||^2 \big) = \frac{1}{2} \big( \mathbb E (U^2) + \mathbb E (\varepsilon(t)^2) -\mathbb E ((U- \varepsilon(t))^2) \big) \\
& = \frac{1}{2} \mathbb E ( 2 U \varepsilon(t)) = g(t).
\end{align*}
These properties together with the following theorem  yield that $\mathcal F (\varepsilon)$ is the RKHS($R$).
\begin{theorem}[E. H. Moor] Aronszajn (1944) \cite{Aronszajn}
	A symmetric non-negative Kernel R generates a unique Hilbert space. \label{unicity}
\end{theorem}

In the sequel, we take $R$ to be continuous on $[0,1]^2$ and we shall
consider the function of interest given  by \eqref{f_xh}. More generally, we consider the function $f$, defined for a continuous function $\varphi$ and $t \in [0,1]$, by\begin{equation} \label{f=int}
	f(t)= \int_0^1 R(s,t) \varphi (s)\;ds.
	\end{equation}
	
	\begin{lemma}\label{finFlemma}
	We have $f \in \mathcal F (\varepsilon)$, i.e., there exists $X \in L_2(\varepsilon)$ with,
	\begin{equation}\label{finF}
	f(\cdot) = \mathbb E(X \varepsilon(\cdot)).
	\end{equation}
	In addition, \[\|f\|^2=\mathbb E(X^2)=\int_{0}^{1}\int_{0}^{1}R(s,t)\varphi (s) \varphi (t)\;dt\;ds.\]
\end{lemma}
\textbf{Proof.}  Define, for a suitable partition $(x_{i,n})_{i =1,\cdots,n}$  of $[0,1]$, 
\[X_n = \sum_{i=1}^{n-1} (x_{i+1,n}-x_{i,n})\varphi(x_{i,n}) \varepsilon(x_{i,n}) \in L_2(\varepsilon),\]such that for any $t\in[0,1]$, \[f(t)= \underset{n \to \infty}{\lim} \sum_{i=1}^{n-1} (x_{i+1,n}-x_{i,n})\varphi(x_{i,n}) R(x_{i,n},t)=\underset{n \to \infty}{\lim} \mathbb E (X_n \varepsilon (t)) .\]
We shall prove that $(X_n)_n$ converges to a certain element of $\mathbb L^2$, i.e., \begin{equation}
\exists~ X \in \mathbb L^2 ~:~ \underset{n \to \infty}{\lim} \mathbb E \big( (X_n-X)^2\big)=0,\label{toprovethis}
\end{equation}
and by the definition of $L_2(\varepsilon)$ the limit in (\ref{toprovethis}) proves that $X$ is an element of $L_2(\varepsilon)$.
Now the proof (\ref{toprovethis}) is immediate, in fact it is easy to check that $(X_n)$ id a Cauchy sequence in  $\mathbb L^2$. By the completeness of $\mathbb L^2$, we deduce (\ref{toprovethis}).
In addition we have,
$\underset{n \to \infty}{\lim} \mathbb E (X_n \varepsilon(t))=\mathbb E (X \varepsilon(t)),$ this is due to the following inequality, \begin{align*}
\Big|\mathbb E (X_n \varepsilon (t))-\mathbb E (X \varepsilon (t)) \Big| & \leq \mathbb E \Big|(X_n -X)\varepsilon (t) \Big| \leq   \sqrt{\mathbb E ((X_n -X)^2)}\sqrt{\mathbb E (\varepsilon (t)^2)},
\end{align*}
and the fact that $\underset{n \to \infty}{\lim}\mathbb E ((X_n -X)^2)= 0$ and $\mathbb E (\varepsilon (t)^2) < \infty$. The proof of \eqref{finF} is concluded.
Finally, \begin{align*}
\mathbb E(X^2)& =\underset{n \to \infty}{\lim} \mathbb E(X_n^2) =\underset{n \to \infty}{\lim}\sum_{i=1}^{n}\sum_{j=1}^{n}(x_{i+1,n}-x_{i,n})(x_{j+1,n}-x_{j,n})\varphi(x_{i,n})\varphi(x_{j,n})R(x_{i,n},x_{j,n})\\
& = \int_{0}^{1}\int_{0}^{1}\varphi(t)\varphi(t)R(s,t)~ds~dt.
\end{align*}This concludes the proof of Lemma \ref{finFlemma}. $\Box$\\[0.3cm]
Now let $T_n = (t_1,t_2,\cdots,t_n)$ with $0 \leq t_1 < t_2 < \cdots < t_n \leq 1$ and let $V_{T_n}$ be the subspace of $\mathcal F (\varepsilon)$ spanned by the functions $R(\cdot,t)$ for $t \in T_n$, i.e., \[V_{T_n} = \{g : [0,1] \to \mathbb R~\text{ with }~g(\cdot)=\mathbb E(U\varepsilon(\cdot))~\text{ where }~ U \in L(\varepsilon(t),t \in T_n) \}.\]
Our task is to prove that if $R_{|T_n}=(R(t_i,t_j)_{1\leq i,j\leq n})$ is a non-singular matrix then $V_{T_n}$ is a closed subspace of $\mathcal F (\varepsilon)$. For this 
let, $(g_m)_{m \geq 1}$ be a sequence in $V_{T_n}$ converging to $g\in \mathcal F (\varepsilon)$. We shall prove that $g \in V_{T_n}$. Note that, \[g_m(t)=\mathbb E (U_m \varepsilon(t))~~\text{ with }~~U_m = \sum_{i=1}^{n}a_{i,m}\varepsilon(t_i),~~\text{ where }~(a_{i,m})_{m \geq 1} \in \mathbb R.\]
Since $g_m$  converges in $\mathcal F (\varepsilon)$ then it is a Cauchy sequence, i.e., \[\underset{m_1,m_2 \to \infty}{\lim~~~}||g_{m_1}-g_{m_2} ||^2=0.\]
By the definition of the norm on $\mathcal F (\varepsilon)$ we have,
\begin{align*}
&||g_{m_1}-g_{m_2} ||^2  =  \mathbb E((U_{m_1}-U_{m_2})^2)=  \mathbb E\Big(  \big(\sum_{i=1}^{n}(a_{i,m_1}-a_{i,m_2})\varepsilon(t_i) \big)^2\Big)\\
&~~ = \sum_{i=1}^{n}\sum_{j=1}^{n}(a_{i,m_1}-a_{i,m_2})(a_{j,m_1}-a_{j,m_2})R(t_i,t_j) = A_{m_1,m_2}' R_{|T_n}A_{m_1,m_2},
\end{align*}
where  $A_{m_1,m_2}' = (a_{1,m_1}-a_{1,m_2},\cdots,a_{n,m_1}-a_{n,m_2})'$. Thus,
\[\underset{m_1,m_2 \to \infty}{\lim~~~} A_{m_1,m_2}' R_{|T_n}A_{m_1,m_2}= 0.\]
Since $R_{|T_n}$ is a symmetric positive matrix, we obtain, \[\underset{m_1,m_2 \to \infty}{\lim~~~}A_{m_1,m_2}'=\underset{m_1,m_2 \to \infty}{\lim~~~}(a_{1,m_1}-a_{1,m_2},\cdots,a_{n,m_1}-a_{n,m_2})'=(0,\dots,0)',\]
which yields that $(a_{i,m})_m$ is a Cauchy sequence on $\mathbb R$ for all $i = 1, \cdots,n.$ Taking $a_i = \underset{m \to \infty}{\lim} a_{i,m}$ we obtain by the uniqueness of the limit, \[g(\cdot)=\mathbb E (U\varepsilon (\cdot))~~\text{with}~~U=\sum_{i=1}^n a_i \varepsilon(t_i),\]
which yields that $g \in V_{T_n}$. Hence $V_{T_n}$ is closed. $\Box$ \\[0.2cm]
Since $V_{T_n}$  is a closed subspace in the Hilbert space $\mathcal F (\varepsilon)$, one can define the orthogonal projection operator from $\mathcal F (\varepsilon)$ to $V_{T_n}$ which we note by $P_{|T_n}$, i.e.,  for every $f \in \mathcal F (\varepsilon)$, \[P_{|T_n} f = \underset{g ~\in V_{T_n}}{\argmin}~||f-g||.\]
Par definition of $P_{|T_n}$, we have for any $g \in V_{T_n}$
$$
\langle P_{|T_n}f- f, g\rangle=0.
$$
Now, for $t_i\in T_n$, $R(\cdot,t_i) \in V_{T_n}$. Hence, for every $i=1,\dots,n.$
\[\langle P_{|T_n}f- f, R(\cdot,t_i) \rangle =0\,\,\, {\mbox{or equivalently}}\,\,\,\langle P_{|T_n}f,R(\cdot,t_i)  \rangle = \langle f,R(\cdot,t_i)  \rangle.\]
The last equality, together with \eqref{propertiesRKHS}, gives that,
\begin{equation}
P_{|T_n}f(\cdot)=f(\cdot)\,\,\, {\mbox{on}}\,\,\, T_n.\label{appendixPtn=fxh}~~\Box
\end{equation}
\subsubsection*{Supplementary facts}
\begin{enumerate}
	\item[(F1)] \label{normdef-g_n} Let $f$ be defined by \eqref{f=int}. We shall prove that if  $g \in V_{T_n}$, i.e., if $g(\cdot )= \sum_{j=1}^n a_j R(t_j,\cdot)$ for some $a_i \in \mathbb R$, then \[||f -g||²=  \int_{0}^{1} \varphi (s)(f(s) -g(s))\;ds- \sum_{i=1}^{n} a_i ( f(t_{i}) -g(t_{i})).\]In fact,
	\begin{align*}
	||f -g||^2 & = \langle f -g,f -g \rangle = \langle f ,f -g\rangle -\langle g,f -g\rangle
	\end{align*}
	On the one hand, note that $f-g \in \mathcal F (\varepsilon)$ and by using \eqref{propertiesRKHS} we obtain,
	\begin{align}
	\langle g,f-g\rangle & = \sum_{i=1}^{n} a_i\langle R(t_i,\cdot) ,f -g\rangle= \sum_{i=1}^{n} a_i ( f(t_{i}) -g(t_{i})).
	\end{align}
	On the another hand,  Lemma \ref{finFlemma} and its proof yield that $f(\cdot)=\mathbb E (X \varepsilon(\cdot))$ where $X \in L_2 (\varepsilon)$ and that, \[\underset{l \to \infty}{\lim}\mathbb E (X_l - X)^2 = 0~~\text{where}~~~X_l= \sum_{j=1}^{l-1} (x_{j+1,l}-x_{j,l})\varphi_{x,h}(x_{j,l})\varepsilon(x_{j,l}), \]
	where $(x_{j,l})_{j =1,\cdots,l}$ is a suitable partition of $[0,1]$. Let $F_l(\cdot) =  \mathbb E(X_l \varepsilon(\cdot))$ which is an element of $\mathcal F (\varepsilon)$.  Clearly,
	\[\langle f ,f -g\rangle=\langle f-F_l ,f -g\rangle+\langle F_l ,f -g\rangle.\]
	We have,
	\begin{align*}
	|\langle f-F_l,f-g\rangle| &\leq  ||f-F_l ||~||f-g||\leq \sqrt{\mathbb E((X_l-X)^2)}||f-g||.
	\end{align*}
	Thus $\underset{l \to \infty}{\lim}\langle f-F_l,f-g\rangle = 0$.
	In addition,
	\begin{align*}
	&\langle F_l,f-g\rangle = \Big\langle \sum_{j=1}^{l-1} (x_{j+1,l}-x_{j,l}) \varphi (x_{j,l})R(x_{j,l},\cdot),f -g \Big\rangle \\
	& = \sum_{j=1}^{l-1} (x_{j+1,l}-x_{j,l}) \varphi (x_{j,l}) \langle  R(x_{j,l},\cdot),f -g\rangle = \sum_{j=1}^{l-1} (x_{j+1,l}-x_{j,l}) \varphi (x_{j,l})(f(x_{j,l}) -g(x_{j,l})).
	\end{align*}
	Hence,
	\begin{align*}
	\underset{l \to \infty}{\lim} \langle F_l,f-g\rangle & = \int_{0}^{1} \varphi (t)(f(t) -g(t))\;dt.
	\end{align*}
	Finally, \[ \langle f, f-g\rangle  = \int_{0}^{1} \varphi (t)(f(t) -g(t))\;dt.~\Box \]
	\item[(F2)] For $x \in [0,1]$, let $f_{x,h}$ be defined by \eqref{f_xh}. We shall prove that, \[m \Var(\hat{g}^{pro}_n (x)) =||P_{|T_n} f_{x,h}||^2. \]
	In fact, by the definition of the projection operator $P_{|T_n}$, we have $P_{|T_n}f_{x,h} \in V_{T_n}$ and for $t \in [0,1]$, \[P_{|T_n}f_{x,h}(t)= \sum_{i=1}^{n} a_i R(t_i,t)= \mathbb E(\sum_{i=1}^na_i\epsilon(t_i)\epsilon(t)) ~~\text{for some}~~a_i \in \mathbb R~\text{for}~i=1,\cdots,n,\] 
and then,
$$
||P_{|T_n}f_{x,h} ||^2=\mathbb E\left(\sum_{i=1}^na_i\epsilon(t_i)\right)^2=\sum_{i=1}^{n} a_i \sum_{j=1}^{n}a_jR(t_i,t_j)=\sum_{i=1}^{n} a_iP_{|T_n}f_{x,h}(t_i).
$$
Recall that ${m_{x,h}}_{|T_n}'={f_{x,h}}_{|T_n}'R_{|T_n}$ and using \eqref{appendixPtn=fxh} we obtain, \begin{equation}\label{mvar=|||PT_n|}
	P_{|T_n}f_{x,h}(t_i)=f_{x,h}(t_i)=\sum_{j=1}^{n} m_{x,h} (t_j)R(t_i,t_j).\end{equation}
	We have then, using \eqref{mvar=|||PT_n|},
\begin{align*}
	& ||P_{|T_n}f_{x,h} ||^2   =\sum_{i=1}^{n} a_i\sum_{j=1}^{n} m_{x,h} (t_j)R(t_i,t_j)
	 = \sum_{j=1}^{n} m_{x,h} (t_j)\sum_{i=1}^{n} a_i R(t_i,t_j) \\ & =\sum_{j=1}^{n} m_{x,h} (t_j)\sum_{i=1}^{n} m_{x,h} (t_i) R(t_i,t_j)  =m\Var(\hat{g}^{pro}_n (x)).~\Box
	\end{align*}
	\item[(F3)] \label{everyfunctioninFiscontnuous} We shall now prove prove that every function in $\mathcal F (\varepsilon)$ is continuous on $[0,1]$. In fact let $g \in \mathcal F (\varepsilon)$, i.e.,
	\[g(\cdot)=\mathbb E(U\varepsilon(\cdot))~~\text{for some}~~ U \in L_2(\varepsilon).\]
	For $s,t \in [0,1]$, (\ref{propertiesRKHS}) and Cauchy-Swartz inequality yields,
	\begin{align*}
	|g(t)-g(s)|& = |\langle R(\cdot,t),g \rangle-\langle R(\cdot,s),g \rangle| = |\langle R(\cdot,t)- R(\cdot,s) ,g \rangle|\\
	&  \leq ||R(\cdot,t)- R(\cdot,s)||~||g||= ||R(\cdot,t)- R(\cdot,s)||~\sqrt{\mathbb E (U^2)}.
	\end{align*}
	Since $\varepsilon$ is of second order process then $\mathbb E (U^2) < \infty$ and since $R$ is continuous on $[0,1]^2$ we obtain,
	\[\underset{s \to t}{\lim}~||R(\cdot,t)- R(\cdot,s)||^2= \underset{s \to t}{\lim} ~(R(t,t)+R(s,s)-2R(s,t) )=0,\]
	which yields that $\underset{s \to t}{\lim}~|g(t)-g(s)|=0$. Hence $g$ is continuous. $\Box$
	\item[(F4)]  \label{f"-g"}Suppose that $R$ verifies Assumptions $(A),(B)$ and $(C)$.  Let $f$ be defined by \eqref{f=int}. We shall prove that if $g \in V_{T_n}$, i.e., $g(\cdot )= \sum_{j=1}^n a_j R(t_j,\cdot)$ with $(a_i)_i \in \mathbb R$ then,
	\[f''(t) - g''(t^+) =-\alpha(t)\varphi(t)+\langle R^{(0,2)}(\cdot,t^+),f-g\rangle. \]
	In fact, we have, as in Equation \eqref{f''t},
	\begin{align*}
	f''(t) &= -\alpha(t) \varphi(t)+ \int_{0}^{1}R^{(0,2)}(s,t^+)\varphi(s)\;ds.
	\end{align*}
	In addition, we have clearly \[g''(t^+) =\sum_{j=1}^n a_j R^{(0,2)}(t_j,t^+). \] Thus,
	\[f''(t) - g''(t^+) =-\alpha(t)\varphi(t)+\int_{0}^{1}R^{(0,2)}(s,t^+)\varphi(s)\;ds-\sum_{j=1}^n a_j R^{(0,2)}(t_j,t^+). \]We have,
	\begin{align*}
	\langle R^{(0,2)}&(\cdot,t^+),f-g\rangle = \langle R^{(0,2)}(\cdot,t^+),f\rangle-\langle R^{(0,2)}(\cdot,t^+),g\rangle
	\end{align*}
	On the one hand, since by Assumption (C), $R^{(0,2)}(\cdot,t^+)$ is in $\mathcal F (\varepsilon)$ then \eqref{propertiesRKHS} yields, \begin{align}
	\langle R^{(0,2)}(\cdot,t^+),g\rangle = \sum_{j=1}^{n} a_j \langle R^{(0,2)}(\cdot,t^+),R(\cdot,t_j)\rangle  = \sum_{j=1}^{n} a_j R^{(0,2)}(t_j,t^+). \label{g''sum}
	\end{align}
	On the other hand, from Lemma \ref{finFlemma} we have $f(\cdot)=\mathbb E (X \varepsilon(\cdot))$ where $X \in L_2(\varepsilon)$ and,
	\[\underset{l \to \infty}{\lim}\mathbb E (X_l - X)^2 = 0~~\text{with}~~~X_l= \sum_{j=1}^{l-1} (x_{j+1,l}-x_{j,l})\varphi(x_{j,l})\varepsilon(x_{j,l}), \]
	where $(x_{j,l})_{j =1,\cdots,l}$ is a suitable partition of $[0,1]$. Let $F_l(\cdot) = \mathbb E  (X_l \varepsilon(\cdot)) \in \mathcal F (\varepsilon)$, we have,
	\begin{align}
	\langle R^{(0,2)}(\cdot,t^+),f\rangle &=\langle R^{(0,2)}(\cdot,t^+),f-F_l\rangle+\langle R^{(0,2)}(\cdot,t^+),F_l\rangle,\label{iterm}
	\end{align}
	and,
	\begin{align*}
	|\langle R^{(0,2)}(\cdot,t^+),f-F_l\rangle|& \leq ||R^{(0,2)}(\cdot,t^+) ||~||f-F_l ||= ||R^{(0,2)}(\cdot,t^+) ||\sqrt{\mathbb E((X_l-X)^2)}.
	\end{align*}
	The last bound together with  Assumption $(C)$ gives  $\underset{l \to \infty}{\lim} |\langle R^{(0,2)}(\cdot,t^+),f-F_l\rangle| =0,$
	in addition,
	\begin{align*}
	\langle R^{(0,2)}(\cdot,t^+),F_l\rangle & = \sum_{j=1}^{l-1} (x_{j+1,l}-x_{j,l})\varphi(x_{j,l})\langle R^{(0,2)}(\cdot,t^+),\varepsilon(x_{j,l})\rangle \\
	&= \sum_{j=1}^{l-1} (x_{j+1,l}-x_{j,l})\varphi(x_{j,l})R^{(0,2)}(x_{j,l},t^+).
	\end{align*}
	Thus, \begin{equation} \underset{l \to \infty}{\lim} \langle R^{(0,2)}(\cdot,t^+),F_l\rangle = \int_0^1 \varphi(s) R^{(0,2)}(s,t^+)~ds. \label{lim(R,Fl)} \end{equation}
	Finally, using \eqref{g''sum}, \eqref{iterm} and \eqref{lim(R,Fl)} yield, \[\langle R^{(0,2)}(\cdot,t^+),f-g\rangle=\int_0^1 \varphi(s) R^{(0,2)}(s,t^+)~ds -\sum_{j=1}^{n} a_j R^{(0,2)}(t_j,t^+).~\Box\]
\end{enumerate}

\end{document}